\documentclass[1p]{elsarticle}

\usepackage{lineno}
\usepackage[colorlinks,linkcolor={blue}]{hyperref}
\modulolinenumbers[5]

\usepackage{amsfonts}
\usepackage{graphicx}
\usepackage{amsmath}
\usepackage{graphicx}
\usepackage{amssymb}

\newtheorem{theorem}{Theorem}

\newtheorem{corollary}[theorem]{Corollary}

\newtheorem{definition}[theorem]{Definition}

\newtheorem{lemma}[theorem]{Lemma}

\newtheorem{proposition}[theorem]{Proposition}

\newtheorem{remark}[theorem]{Remark}
\newenvironment{proof}[1][Proof]{\textbf{#1.} }{\ \rule{0.5em}{0.5em}}

\newcommand{\comment}[1]{}
\newcommand{\cover}[1]{\stackrel{#1}{\Longrightarrow}}
\newcommand{\invcover}[1]{\stackrel{#1}{\Longleftarrow}}

\newcommand{\inter}{\mathrm{int}}
\newcommand{\id}{\mathrm{Id}}

\makeatletter
\usepackage{tikz}
\newcommand*\circled[2][1.6]{\tikz[baseline=(char.base)]{
    \node[shape=circle, draw, inner sep=1pt,
        minimum height={\f@size*#1},] (char) {\vphantom{WAH1g}#2};}}
\newcommand\NoStart{\circled[0.0]{$N_0$} }
\makeatother

\journal{Journal of Differential Equations}

\begin{document}
\begin{frontmatter}
\title{Oscillatory Motions and Parabolic Manifolds at Infinity in the Planar Circular Restricted Three Body Problem}

\author{Maciej J. Capi\'nski\footnote{M. C. has been partially supported by the NCN grant 2018/29/B/ST1/00109}}
\ead{maciej.capinski@agh.edu.pl}
\address{AGH University of Science and Technology, al. Mickiewicza 30, 30-059 Krak\'ow, Poland}

\author{Marcel Guardia\footnote{M. G. has received funding from the European Research Council (ERC) under the European Union’s Horizon 2020 research and innovation programme (grant agreement
No 757802).  M. G. is supported by the Catalan Institution for Research and Advanced Studies via an ICREA
Academia Prize 2019.}}
\ead{marcel.guardia@upc.edu}
\address{Universitat Polit\`ecnica de Catalunya, Departament de Matem\`atiques \& IMTECH, Pau Gargallo 14, Barcelona, Spain \& Centre de Recerca Matem\`atica,  Campus de Bellaterra, Edifici C, Barcelona, Spain}

\author{Pau Mart\'in \footnote{P. M. has been partially supported by the Spanish MINECO-FEDER Grant PGC2018-100928-B-I00 and the Catalan grant 2017SGR1049}}
\ead{p.martin@upc.edu}
\address{Universitat Polit\`ecnica de Catalunya, Departament de Matem\`atiques \& IMTECH, Pau Gargallo 14, Barcelona, Spain \& Centre de Recerca Matem\`atica,  Campus de Bellaterra, Edifici C, Barcelona, Spain}

\author{Tere Seara\footnote{T. S. has been also partly supported by the Spanish MINECO-FEDER Grant
PGC2018-098676-B-100 (AEI/FEDER/UE), the Catalan grant 2017SGR1049 and  by the Catalan Institution for Research and Advanced Studies via an ICREA
Academia Prize 2019.}}
\ead{tere.m-seara@upc.edu}
\address{Universitat Polit\`ecnica de Catalunya, Departament de Matem\`atiques \& IMTECH, Pau Gargallo 14, Barcelona, Spain \& Centre de Recerca Matem\`atica,  Campus de Bellaterra, Edifici C, Barcelona, Spain}

\author{Piotr Zgliczy\'nski\footnote{P. Z. has been partially supported by the NCN grant 2019/35/B/ST1/00655 }}
\ead{umzglicz@cyf-kr.edu.pl}
\address{Jagiellonian University, ul. prof. Stanis\l awa \L ojasiewicza 6,
30-348 Krak\'ow, Poland}

\begin{abstract}
Consider the Restricted Planar Circular 3 Body Problem with both realistic mass ratio and Jacobi constant for the Sun-Jupiter pair. We prove the existence of all possible combinations of past and future final motions. In particular, we obtain the existence of oscillatory motions. All the constructed trajectories cross the orbit of Jupiter but avoid close encounters with it. The proof relies on the method of correctly aligned windows and is computer assisted.

\end{abstract}

\end{frontmatter}




\section{Introduction}

The planar circular restricted three
body problem (PCR3BP) models the motion of a body of zero mass under the Newtonian graviational force exerted by two bodies of masses $\mu$ and $1-\mu$ which evolve in circular motion around their center of mass on the same plane. In rotating coordinates, if we denote by $q\in\mathbb{R}^2$ the position of the zero mass body and $p\in \mathbb{R}^2$ its associated momentum, the PCR3BP is a Hamiltonian system with respect to
\begin{equation}\label{def:HamCartesian}
H(q,p;\mu)= \frac{\|p\|^2}{2} -q_1p_2+q_2p_1- \frac{1-\mu}{\left\|q+\mu
      \right\|}-\frac{\mu}{\left\|q-(1-\mu)\right\|}.
\end{equation}

Since the Hamiltonian is autonomous, it is a first integral which correspond to the Jacobi constant in non-rotating coordinates. (Often the Jacobi constant is defined as $J=-2H$).

%
%
In the 1922, J. Chazy classified the possible final motions the massless body in the PCR3BP may possess \cite{Chazy22} (see also~\cite{ArnoldKN88}), that is, the  possible ``states'' that $q(t)$ may possess as $t\to\pm\infty$. They can be:
\begin{itemize}
 \item $H^\pm$ (hyperbolic):  $\|q(t)\|\rightarrow\infty$ and $\|\dot q(t)\|\rightarrow c>0$ as $t\rightarrow\pm\infty$.
\item $P^\pm$ (parabolic): $\|q(t)\|\rightarrow\infty$ and $\|\dot q(t)\|\rightarrow 0$ as $t\rightarrow\pm\infty$.
\item $B^\pm$ (bounded): $\limsup_{t\rightarrow \pm\infty}\|q\|<+\infty$.
\item $OS^\pm$ (oscillatory): $\limsup_{t\rightarrow \pm\infty}\|q\|=+\infty$ and  $\liminf_{t\rightarrow \pm\infty}\|q\|<+\infty$.
\end{itemize}
Examples of all these types of motion, except the oscillatory ones,
were already known by Chazy.

Oscillatory motions were proven to exist for the first time by Sitnikov \cite{Sitnikov60} in the 1960's for the nowadays called \emph{Sitnikov model}, which is a symmetric restricted spatial 3 body problem. Moreover he proved that one can construct orbits with any prescribed past and future final motions.

The approach by Sitnikov and by most of the subsequent references (see below) to construct oscillatory motions strongly rely on perturbative methods (see for instance \cite{LlibreS80, GuardiaMS16}). As a consequence, the motions obtained are either confined to ``small'' regions of the phase space or only exist in certain narrow ranges of some parameters.

The purpose of this paper is to develop techniques to prove such behaviors in \emph{non-perturbative} regimes. These techniques will rely on Computer Assisted Proofs. This will allow to deal with physical ranges of parameters and regions of the phase space
(that is  regions  ``close to'' the orbits of the Sun and Jupiter).

We consider the PCR3BP and apply these techniques to prove the existence of any combination of past and future motions, including the oscillatory ones, for some  explicitly
given values of the mass parameter $\mu $ and  energy level $H$ (equivalently for a given value of the  Jacobi constant). More concretely, we consider $\mu =0.0009537$ which corresponds to the mass ratio for the pair Sun-Jupiter and $H=-1$.
%

\begin{theorem}
\label{th:main} Consider the PCR3BP, that is, \eqref{def:HamCartesian} with  $\mu =0.0009537$. Then,
\[
 X^+\cap Y^-\cap \{H=-1\}\neq\emptyset,
\]
where $X,Y=H,P,B,OS$.

In particular, for $r_0:=0.5002$, 
\begin{itemize}
\item There exist trajectories $(q(t),p(t))$ such that
\begin{equation}\label{def:oscillatory}
\limsup_{t\to\pm\infty}\|q(t)\|=+\infty \qquad \text{ and }\qquad \liminf_{t\to\pm\infty}\|q(t)\| \le r_0.
\end{equation}
\item 
For every sufficiently large $K\gg 1$, there exists a periodic orbit $(q(t),p(t))$  such that
\[
\sup_{t\in\mathbb{R}}\|q(t)\|\geq K \qquad \text{ and }\qquad \inf_{t\in\mathbb{R}}\|q(t)\| \le r_0.
\]
\end{itemize}
\end{theorem}

This theorem gives the possibility of combining any past and future final motions at a given energy level and with a realistic mass ratio.
Moreover, Item 1 in the theorem implies that there exist oscillatory motions which reach points which are closer to the Sun than Jupiter.
That is, these oscillatory orbits cross the orbit of Jupiter (but stay away from collision with it). Item 2 of the theorem gives the existence of  periodic orbits of the PCR3BP (in the rotating frame) which encircle Jupiter and go very far from the primaries.
%
%
%
%

%
%

This result focuses on the energy level $H=-1$ since it is far from the limit cases
where oscillatory motions can be proven analytically
(see \cite{LlibreS80, GuardiaMS16}, where the value of $-H$ needs to be taken sufficiently large).
Our methodology can be
applied also to different energy levels and there is nothing
special about  $H=-1$. (In fact, from our  proof it
follows that there are orbits with oscillatory motions at any energy level sufficiently close to  $H=-1$ .)
%

The analysis of final motions, and in particular of oscillatory motions, has drawn considerable attention in the last decades since the pioneering work by  Sitnikov \cite{Sitnikov60}. In 1968,  Alekseev extended the results of Sitnikov constructing all possible combinations of future and past final motions (and thus oscillatory motions) for the full three body problem  assuming the third mass is small enough \cite{Alekseev68}.

Later, J. Moser \cite{Moser01} gave a new proof which related the existence of oscillatory motions to symbolic dynamics.
His approach  has been very influential and has been applied to different Restricted 3 Body Problems \cite{SimoL80,LlibreS80,GuardiaMS16,GuardiaPSV20} and, roughly speaking, it is also applied in the present paper.
Moeckel has also proved the existence of oscillatory motions via symbolic dynamics for the three body problem relying on dynamics close to triple collision \cite{Moeckel07} (and therefore for sufficiently small total angular momentum). Oscillatory motions have been also constructed relying on other techniques closer to those of Arnold diffusion \cite{GuardiaSMS17,Seara20}.

Concerning the PCR3BP, \cite{SimoL80} gives the existence of oscillatory motions and symbolic dynamics for arbirtarily large Jacobi constant
assuming the mass ratio to be exponentially small with respect to the Jacobi constant. The paper \cite{GuardiaMS16} proved the result for any mass ratio and large enough Jacobi constant. As mentioned before, both references strongly rely on perturbative methods and only apply to nearly integrable regimes.
Note that, in particular, the last result is \emph{non-perturbative} with respect to the mass but requires large Jacobi constant which implies that the orbit are extremely far from the orbits of the Sun and Jupiter.
\begin{figure}
\begin{center}
\includegraphics[height=5.5cm]{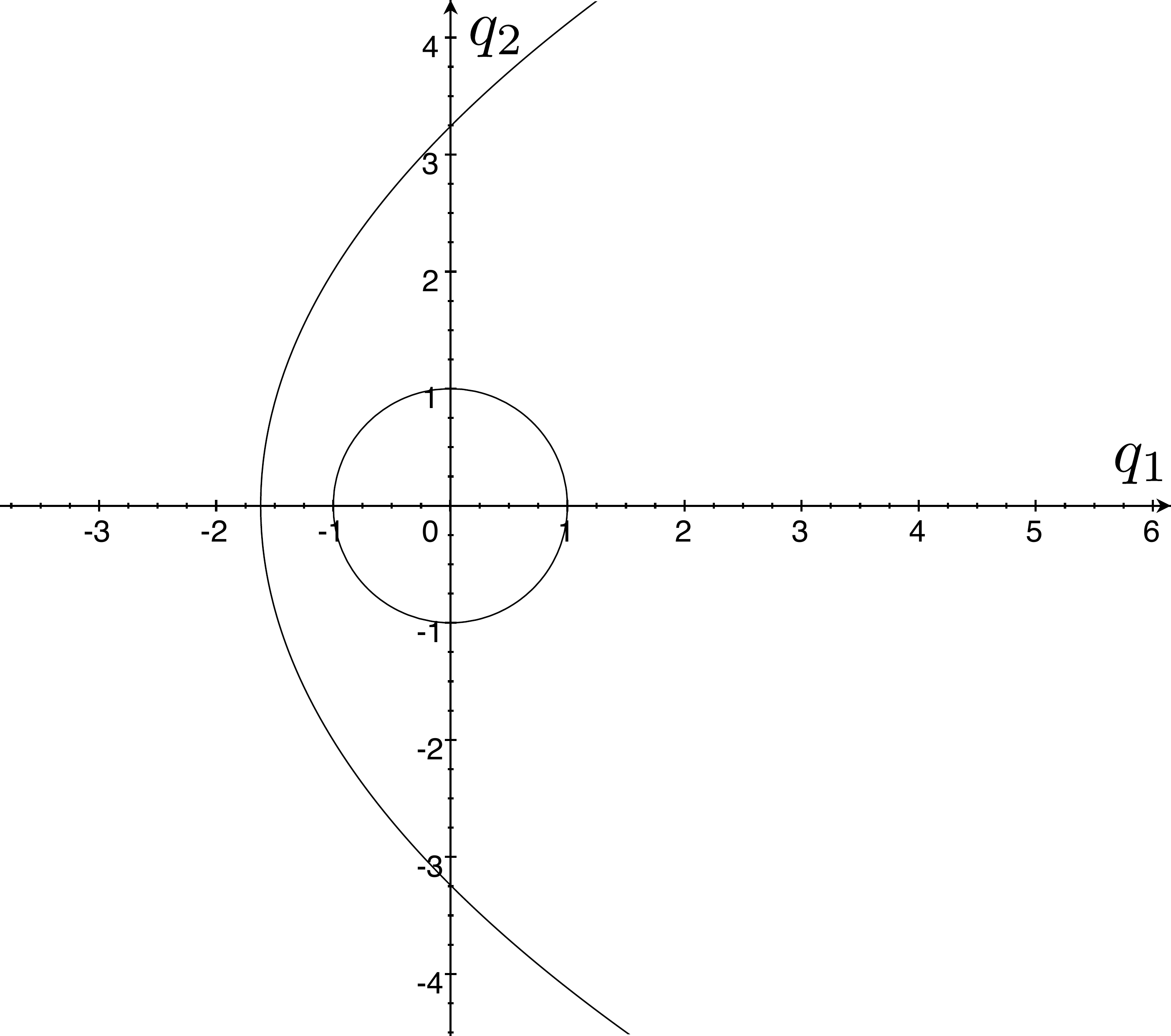}\qquad\includegraphics[height=5.5cm]{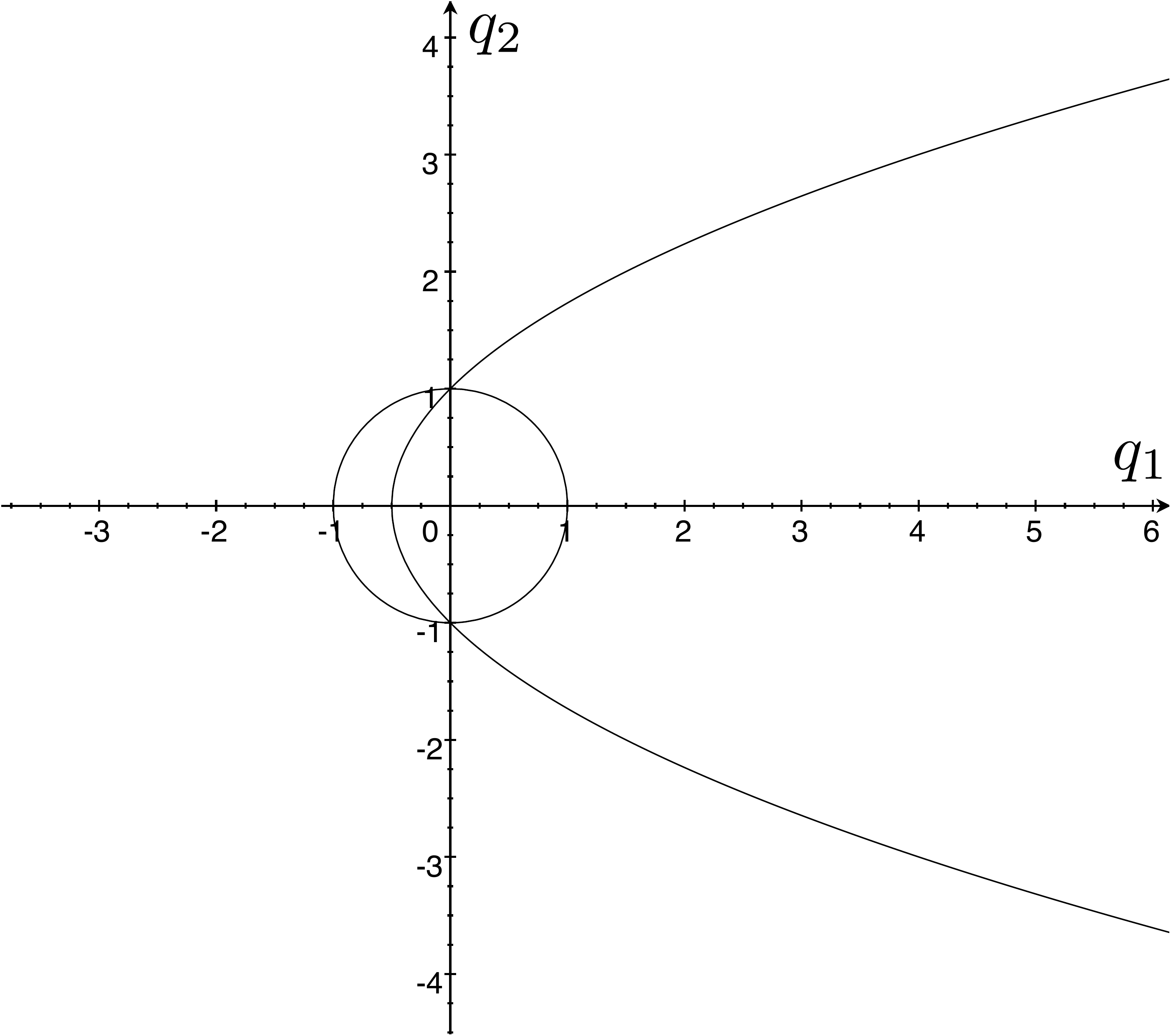}
\end{center}
\caption{Parabolic motion from infinity for $H=-1.8$ on the left (which is the case considered in \cite{GalanteK11,GalanteK10b,GalanteK10c}), and on the right for $H=-1$ as is considered in the current paper. The circle represents the path of Jupiter and the Sun is positioned in the origin. The plots are for the limit case $\mu=0$, but are very close to the trajectories in the Jupiter-Sun system.}\label{fig:energies}
\end{figure}

As far as the authors know, the only papers which deal with realistic values of  both the mass ratio and energy level   are \cite{GalanteK11,GalanteK10b,GalanteK10c} which are also computer assisted. Relying on completely different techniques from those of \cite{Moser01},  Kaloshin and Galante construct trajectories whose initial conditions are ``within the range of the Solar system'' and become oscillatory as $t\to+\infty$. These orbits have energy $H\leq -1.52$ (the needed conditions are rigorously verified with the aid of computers for $H=-1.8$). See Figure \ref{fig:energies} for the difference between their choice of energy and ours. 

Not only the methods by Kaloshin and Galante are very different from ours, but also the orbits they construct are very different from those constructed in the present paper. In particular, their orbits undergo a drastic change in eccentricity whereas they stay away from the orbits of the primaries (note that the condition $H\leq -1.52$ implies that the outer Hill region is disconnected from the inner ones). On the contrary, the oscillatory trajectories in the present paper have rather high eccentricity but can cross the trajectory of Jupiter.
%

\subsection{The Moser approach and its implementation}
Let us finish the introduction by explaining the approach that Moser developed to prove the existence of oscillatory motions for the Sitnikov problem and how his ideas are adapted in the present paper to prove Theorem \ref{th:main}.

The Sitnikov problem is a Hamiltonian system of one a half degrees of freedom (one degree of freedom plus periodic time dependence). Let us denote it by $\mathcal{H}=\mathcal{H}(q,\dot q, t)$ (its particular form is now not important).
If one performs a  suitable change of coordinates and considers the stroboscopic Poincar\'{e} map, the ``parabolic infinity'' $q=+\infty$, $\dot q=0$, can be seen as a fixed point. The linearization of the stroboscopic map at this point is degenerate (equal to the identity).

The first step of Moser's approach is to prove that, even if the fixed point is degenerate, it  possesses one-dimensional stable and unstable invariant manifolds,  which correspond to  the parabolic orbits $P^\pm$. For the Sitnikov and the PCR3BP this fact had been proven previously in \cite{McGehee73}.

The second step is to prove that these invariant manifolds intersect transversally. This is the step which crucially relies on classical perturbative techniques such as Melnikov Theory \cite{Melnikov63} (as in \cite{Moser01, SimoL80}) or singular perturbative techniques to deal with exponentially small phenomena (as in \cite{GuardiaMS16}).

If the fixed point at infinity was hyperbolic, a standard adaptation of Smale Theorem \cite{Smale65} (based on the classical Lambda Lemma) would lead to symbolic dynamics and oscillatory motions. However, since it is degenerate one needs to analyze carefully the dynamics close to the fixed point by a specific ``parabolic Lambda lemma''.


In this paper, relying on the just explained Moser approach, we develop techniques which can be implemented in a computer to prove the existence of oscillatory motions.

First, in Section \ref{sec:weak-nhims}, we prove the existence of the local invariant manifolds and obtain quantitative estimates of its graph parameterizations. The approach, by using cone-shaped isolating blocks, is reminiscent to that of \cite{McGehee73}.

Then, in Sections~\ref{sec:eqinfty} and~\ref{sec:InvManPCR3B}, these techniques are applied to the PCR3BP. First, in Section~\ref{sec:eqinfty}, we perform several changes of coordinates to the PCR3BP so that it fits into the framework of Section~\ref{sec:weak-nhims}. Then, in Section~\ref{sec:InvManPCR3B}, we obtain estimates of the local invariant manifolds. In this section, we also extend them by the flow. This extension, computer assisted, is done in a way that one obtains fine rigorous estimates for the global invariant manifolds.

The fact that the infinity fixed point is degenerate implies that this computer implementation is by no means standard. Indeed, the dynamics in the invariant manifolds is extremely slow (its decay to the fixed point is polynomial in time instead of exponential).

The invariant manifolds intersect thanks to the Hamiltonian structure, and moreover, one can easily locate (some of) the intersections thanks to the reversibility of the PCR3BP with respect to the involution
\[
(q_1,q_2,p_1,p_2)\to (q_1,-q_2,-p_1,p_2).
\]
Our method does not require proving that the invariant manifolds intersect transversally. However, the method to construct oscillatory motions explained in the next paragraph implicitly relies on the fact that these invariant manifolds have intersections which are topologically transverse.

Finally, in Section~\ref{sec:oscillatory-proof}, we use the methods of correctly aligned windows (covering relations) \cite{MR2060531,Easton75,Easton78} to construct the motions described in Theorem~\ref{th:main}. More precisely, we construct a sequence of windows which go from a small neighborhood of the fixed point at infinity to a neighborhood of one of the intersections between the stable and unstable invariant manifolds. Relying on the analysis of the local dynamics close to infinity in Sections~\ref{sec:eqinfty} and~\ref{sec:InvManPCR3B} and integrating with rigorous numerics the flow of the PCR3BP, we show that these windows are correctly aligned. Different choices of sequence of windows lead to different final motions. If one chooses a sequence such that (some of) the windows get closer and closer to the invariant manifolds of infinity, the orbits ``hitting'' this sequence of windows are oscillatory. On the contrary, if one chooses the windows uniformly away  from the invariant manifolds (for instance one can take a fixed loop of correctly aligned windows), the corresponding orbits are bounded.
Moreover, we also show that orbits passing through the edges of some of our windows lead to hyperbolic motions and that orbits reaching the parabolic stable/unstable manifolds lead to parabolic motions. 
From our topological construction it follows that we can link all these types of motions in forward and backward time. 
%
%
%
%


\section{Preliminaries}

First we introduce some notation, which will be used throughout the paper. We
write $B_{k}$ for a open unit ball in $\mathbb{R}^{k},$ centered at zero, under
some norm of our choice. (In our application we will use the max norm, but
many of the arguments can be made norm independent.) We will write
$\overline{B_{k}}$ for the closure of $B_{k}$.

For a given norm $\|\cdot\|$ on $\mathbb{R}^n$ and for a matrix $A\in \mathbb{R}^{n \times n}$ we define%
\begin{equation}\label{def:matrixconst}
\begin{split}
m\left(A\right)   &  =\min_{p\in\mathbb{R}^{n},\left\Vert p\right\Vert
=1}\left\Vert Ap\right\Vert ,\\
l\left(A\right)   &  =\lim_{h\rightarrow0^{+}}\frac{\left\Vert
\id+hA\right\Vert -\left\Vert \id\right\Vert }{h},\\
m_{l}\left(A\right)   &  =\lim_{h\rightarrow0^{+}}\frac{m\left(
\id+hA\right)  -\left\Vert \id\right\Vert }{h}.
\end{split}
\end{equation}
The $l(A)$ is the logarithmic norm of $A$ \cite{MR0100966,MR0145662}. It is known that $l(A)$ is a convex
function. The number $m(A)$ is useful to us since it can be used to obtain the lower
bound $\left\Vert Ap\right\Vert \geq m\left(A\right)  \left\Vert
p\right\Vert $. The number $m_{l}\left(A\right)$ can be thought of as a
`lower bound version' of the logarithmic norm.

If $s >0$, then 
\begin{equation*}
  l(sA)=s l(A), \quad m_l(sA)=s m_l(A).
\end{equation*}

\begin{lemma}[\cite{MR3567489}]
\label{lem:norms-computation}We have%
\[
m\left(  A\right)  =\left\{
\begin{array}
[c]{lll}%
\frac{1}{\left\Vert A^{-1}\right\Vert } &  & \text{if }\det A\neq0,\\
0 &  & \text{otherwise,}%
\end{array}
\right.  \qquad m_{l}\left(  A\right)  =-l\left(  -A\right)  .
\]

\end{lemma}

We now give two technical lemmas that allow us to obtain upper and lower
bounds on $l$ and $m_{l}$, respectively, for a weighted average of a family of matrices.

\begin{lemma}
\label{lem:norm-integral-1}Let $h:\left[  0,1\right]  \rightarrow
\mathbb{R}_{+}$ be a continuous function and let $A:\left[  0,1\right]
\rightarrow\mathbb{R}^{n\times n}$ be a measurable function such that
$l\left(  A\left(  s\right)  \right)  \leq L$ for $s\in\left[  0,1\right]  .$
Then%
\[
l\left(  \int_{0}^{1}h\left(  s\right)  A\left(  s\right)  ds\right)  \leq
L\int_{0}^{1}h\left(  s\right)  ds.
\]

\end{lemma}

\begin{proof}
From Jensen's inequality applied to the convex function $l$ we obtain
\[
l\left(  \int_{0}^{1}h(s)A(s)ds\right)  \leq\int_{0}^{1}l(h(s)A(s))ds=\int
_{0}^{1}h(s)l(A(s))ds\leq L\int_{0}^{1}h(s)ds,
\]
as required.
\end{proof}

\begin{lemma}
\label{lem:norm-integral-2}Let $h:\left[  0,1\right]  \rightarrow
\mathbb{R}_{+}$ be a continuous function and let $A:\left[  0,1\right]
\rightarrow\mathbb{R}^{n\times n}$ be a measurable function such that
$m_{l}\left(  A\left(  s\right)  \right)  \geq M$ for $s\in\left[  0,1\right]
.$ Then%
\[
m_{l}\left(  \int_{0}^{1}h\left(  s\right)  A\left(  s\right)  ds\right)  \geq
M\int_{0}^{1}h\left(  s\right)  ds.
\]

\end{lemma}

\begin{proof}
Since $A\mapsto l\left(  -A\right)  $ is convex, $m_{l}(A)=-l(-A)$ is concave,
so from Jensen's inequality%
\begin{align*}
m_{l}\left(  \int_{0}^{1}h(s)A(s)ds\right)   &  \geq\int_{0}^{1}m_{l}
(h(s)A(s))ds\\
&  =\int_{0}^{1}h(s)m_{l}(A(s))ds\geq M\int_{0}^{1}h(s)ds,
\end{align*}
as required.
\end{proof} 


\section{Topologically hyperbolic
manifolds\label{sec:weak-nhims}}
Let $\Lambda $ be a smooth compact $c$-dimensional manifold without
boundary. We will consider a vector field 
\begin{equation*}
F:\mathbb{R}^{u}\times \mathbb{R}^{s}\times \Lambda \rightarrow \mathbb{R}%
^{u}\times \mathbb{R}^{s}\times T\Lambda
\end{equation*}
(here $T\Lambda $ is the tangent space) and an ODE%
\begin{equation}
p^{\prime }=F(p).  \label{eq:our-ode}
\end{equation}%
We shall write $\Phi _{t}$ for the flow induced by (\ref{eq:our-ode}). We
shall assume that 
\begin{equation*}
\tilde{\Lambda}=\left\{ \left( 0,0\right) \right\} \times \Lambda
\end{equation*}%
is invariant under the flow. We will be in the context where for $p\in 
\tilde{\Lambda}$ the derivative $DF(p)$ can be zero.

The objective of this section will be to provide sufficient conditions for
the existence and the construction of stable and unstable manifolds of $%
\tilde{\Lambda}$.

Let us introduce the following notation for coordinates: $x\in \mathbb{R}%
^{u},$ $y\in \mathbb{R}^{s},$ $\lambda \in \Lambda $. The coordinate $x$ is
towards the expanding direction, $y$ is towards the contracting direction
and $\lambda $ is the center direction. We do not need to assume though that
the coordinates $x$ and $y$ are perfectly aligned with the unstable and
stable bundles of our system, respectively. A `rough alignment' will turn
out good enough, provided that the conditions needed for our construction
are fulfiled.

Let $L\in (0,1]$ be a fixed constant, let $\beta _{u}\subset \left\{
1,\ldots ,u\right\} ,$ $\beta _{s}\subset \left\{ 1,\ldots ,s\right\} $ (the
sets $\beta _{u},\beta _{s}$ can be empty) and consider the following sets 
\begin{eqnarray*}
S^{u} &=&S_{L}^{u}=\left\{ \left( x,y,\lambda \right) :\lambda \in \Lambda
,\left\Vert y\right\Vert <L\left\Vert x\right\Vert ,\left\Vert x\right\Vert
<1,x_{i}>0\text{ for }i\in \beta _{u}\right\} , \\
S^{s} &=&S_{L}^{s}=\left\{ \left( x,y,\lambda \right) :\lambda \in \Lambda
,\left\Vert x\right\Vert <L\left\Vert y\right\Vert ,\left\Vert y\right\Vert
<1,y_{i}>0\text{ for }i\in \beta _{s}\right\} .
\end{eqnarray*}

We define%
\begin{eqnarray*}
S_{-}^{u} &=&\left\{ \left( x,y,\lambda \right) \in \overline{S^{u}}%
:\left\Vert x\right\Vert =1\text{ or }x_{i}=0\text{ for some }i\in \beta
_{u}\right\} , \\
S_{-}^{s} &=&\left\{ \left( x,y,\lambda \right) \in \overline{S^{s}}%
:\left\Vert y\right\Vert =1\text{ or }y_{i}=0\text{ for some }i\in \beta
_{s}\right\} .
\end{eqnarray*}%
We shall refer to $S_{-}^{u},S_{-}^{s}$ as \emph{exit sets}. We also define%
\begin{eqnarray*}
S_{+}^{u} &=&\left\{ \left( x,y,\lambda \right) \in \overline{S^{u}}%
:\left\Vert y\right\Vert =L\left\Vert x\right\Vert \text{ and }\left\Vert
x\right\Vert <1\right\} , \\
S_{+}^{s} &=&\left\{ \left( x,y,\lambda \right) \in \overline{S^{s}}%
:\left\Vert x\right\Vert =L\left\Vert y\right\Vert \text{ and }\left\Vert
y\right\Vert <1\right\} .
\end{eqnarray*}%
We shall refer to $S_{+}^{u},S_{+}^{s}$ as \emph{entry sets}. Note that 
\begin{equation*}
\partial S^{u}=S_{-}^{u}\cup S_{+}^{u}\qquad \text{and\qquad }\partial
S^{s}=S_{-}^{s}\cup S_{+}^{s}.
\end{equation*}

\begin{definition}
\label{def:Scu-sector} We say that $S^{u}$ is an \emph{unstable sector} if
for every $p\in S^{u}$ the forward trajectory leaves $S^{u}$ through $%
S_{-}^{u}$ and enters through $S_{+}^{u}$. More precisely, if the following
two conditions are satisfied:

\begin{enumerate}
\item If $p\in S_{-}^{u}$ then $\Phi _{\lbrack 0,t]}\left( p\right) \notin
S^{u}$ for some $t>0$,

\item If $p\in S_{+}^{u}$ then $\Phi _{(0,t]}\left( p\right) \in S^{u}$ for
some $t>0$.
\end{enumerate}
\end{definition}

\begin{definition}
\label{def:Scs-sector} We say that $S^{s}$ is a \emph{stable sector}, if it
is a unstable sector for the flow with reversed time$.$
\end{definition}

The sets $S^{u}$ and $S^{s}$ will provide bounds for the domains in which
the manifolds are positioned. To simplify the statements, we will sometimes
refer to $S^{u}$ and $S^{s}$ as sectors.

\begin{remark}
Depending on the particular system the sets $\beta _{u},\beta _{s}$ can be
empty. We consider them since in the equations of the PRC3BP at infinity,
some of the coordinates will only have physical meaning when they are
greater or equal to zero.
\end{remark}

\begin{definition}
\label{def:parab-manifold}We will define the unstable and stable sets as 
\begin{align*}
W^{u}(\tilde{\Lambda})& =\left\{ p\in S^{u}:\Phi _{t}\left( p\right) \in
S^{u}\text{ for }t\in (-\infty ,0]\text{ and }\lim_{t\rightarrow -\infty }%
\mathrm{dist}\left( \Phi _{t}\left( p\right) ,\tilde{\Lambda}\right)
=0\right\} , \\
W^{s}(\tilde{\Lambda})& =\left\{ p\in S^{s}:\Phi _{t}\left( p\right) \in
S^{s}\text{ for }t\in \lbrack 0,+\infty )\text{ and }\lim_{t\rightarrow
+\infty }\mathrm{dist}\left( \Phi _{t}\left( p\right) ,\tilde{\Lambda}%
\right) =0\right\} ,
\end{align*}%
respectively.
\end{definition}

In our work we will present tools which will allow us to establish the
existence of unstable and stable sets which are graphs of Lipschitz
functions 
\begin{eqnarray*}
w^{u} &:&\pi _{x,\lambda }S^{u}\rightarrow S^{u}, \\
w^{s} &:&\pi _{y,\lambda }S^{s}\rightarrow S^{s},
\end{eqnarray*}%
\begin{equation*}
W^{u}(\tilde{\Lambda})=\mathrm{graph}\left( w^{u}\right) =w^{u}\left( \pi
_{x,\lambda }S^{u}\right) ,\qquad W^{s}(\tilde{\Lambda})=\mathrm{graph}%
\left( w^{s}\right) =w^{s}\left( \pi _{y,\lambda }S^{s}\right) ,
\end{equation*}%
where $w^{u},$ $w^{s}$ satisfy 
\begin{equation*}
\pi _{x,\lambda }w^{u}\left( x,\lambda \right) =\left( x,\lambda \right)
,\qquad \pi _{y,\lambda }w^{s}\left( y,\lambda \right) =\left( y,\lambda
\right) .
\end{equation*}%
This will in particular mean that $W^{u}(\tilde{\Lambda})$ and $W^{s}(\tilde{%
\Lambda})$ are Lipschitz manifolds.

We now define what we will mean by saying that $\tilde{\Lambda}$ is a \emph{%
topologically hyperbolic} manifold.

\begin{definition}
Assume that the unstable and stable sets are manifolds. If we have a
neighbourhood $U$ of $\tilde{\Lambda}$ in which all points whose forward
trajectories remain in $U$ are in $W^{s}(\tilde{\Lambda})$, and all points
from $U$ whose backward trajectories remain in $U$ are contained in $W^{u}(%
\tilde{\Lambda})$, then we call $\tilde{\Lambda}$ a topologically hyperbolic
manifold.
\end{definition}

We will be working under the assumption that in $S^{u}$ and $S^{s}$ we can
factor out suitable terms from the derivative of $F$. From now on let us
focus on the sector $S^{u}$ within which we will establish the existence of
the manifold $W^{u}(\tilde{\Lambda})$. (The results for $W^{s}(\tilde{\Lambda%
})$ within $S^{s}$ will follow by reversing the sign of the vector field,
and swapping the roles of the coordinates $x,y$.)

For the factorisation of the suitable terms in $S^{u}$ we will assume that
there exist functions $h:S^{u}\rightarrow \mathbb{R}$ and $%
G:S^{u}\rightarrow \mathcal{L}\left( \mathbb{R}^{u}\times \mathbb{R}%
^{s}\times T\Lambda \right) $, (here $\mathcal{L}(X)$ stands for the space
of Linear operators on $X$), such that 
\begin{equation}
h\left( x,y,\lambda \right) >0,\qquad \text{for all }\left( x,y,\lambda
\right) \in S^{u},  \label{eq:h>0}
\end{equation}%
and 
\begin{equation}
DF\left( x,y,\lambda \right) =h(x,y,\lambda )G\left( x,y,\lambda \right) .
\label{eq:factorised-derivative}
\end{equation}

\begin{remark}
Note that we allow $h=0$ on $\tilde\Lambda$. We also note that the case of (%
\ref{eq:our-ode}--\ref{eq:factorised-derivative}) is fundamentally different
from considering the case where we have an ODE with a vector field $F\left(
p\right) =h\left( p\right) g(p)$ and where $h>0$ and $g$ has a NHIM. The
latter case is trivial since the NHIM for $g$ and its associated stable and
unstable manifolds become invariant manifolds for $F$ by a simple rescaling
of time. 
\end{remark}

Our objective will be to impose some normally-hyperbolic-type conditions on $%
G$ in (\ref{eq:factorised-derivative}), from which we will be able to deduce
the existence $W^{u}(\tilde{\Lambda})$ in $S^{u}$. Our methods will lead to
establishing the existence of the function $w^{u}$, which will be Lipschitz.
They can be applied in a more general context, but to simplify the arguments
we restrict to the case where $\Lambda $ is a $c$-dimensional torus $\Lambda
=\mathbb{S}^{c}=\left( \mathbb{R}/\left( \mbox{mod}\,2\pi \right) \right)
^{c}$. Then we are in a convenient situation, since we have a covering%
\begin{equation}
\varphi :\mathbb{R}^{c}\rightarrow \left( \mathbb{R}/\left( \mbox{mod}\,2\pi
\right) \right) ^{c},  \label{eq:torust-lift}
\end{equation}%
which gives us local charts as restrictions of $\varphi $ to balls, provided
that the radius of such balls is smaller than $\pi $.


\subsection{Cone conditions and outflowing along cones}

Let $L_{u},L_{cu},L_{s},L_{cs}>0$%
\begin{equation}
Q_{u},Q_{cu},Q_{s},Q_{cs}:\mathbb{R}^{u+s+c}\rightarrow \mathbb{R},
\label{eq:cones}
\end{equation}%
defined as%
\begin{align*}
Q_{u}\left( x,y,\lambda \right) & :=L_{u}\left\Vert x\right\Vert -\left\Vert
\left( y,\lambda \right) \right\Vert , \\
Q_{cu}\left( x,y,\lambda \right) & :=L_{cu}\left\Vert \left( x,\lambda
\right) \right\Vert -\left\Vert y\right\Vert , \\
Q_{s}\left( x,y,\lambda \right) & :=L_{s}\left\Vert y\right\Vert -\left\Vert
\left( x,\lambda \right) \right\Vert , \\
Q_{cs}\left( x,y,\lambda \right) & :=L_{cs}\left\Vert \left( y,\lambda
\right) \right\Vert -\left\Vert x\right\Vert ,
\end{align*}%

\begin{figure}[tbp]
\begin{center}
\includegraphics[height=4cm]{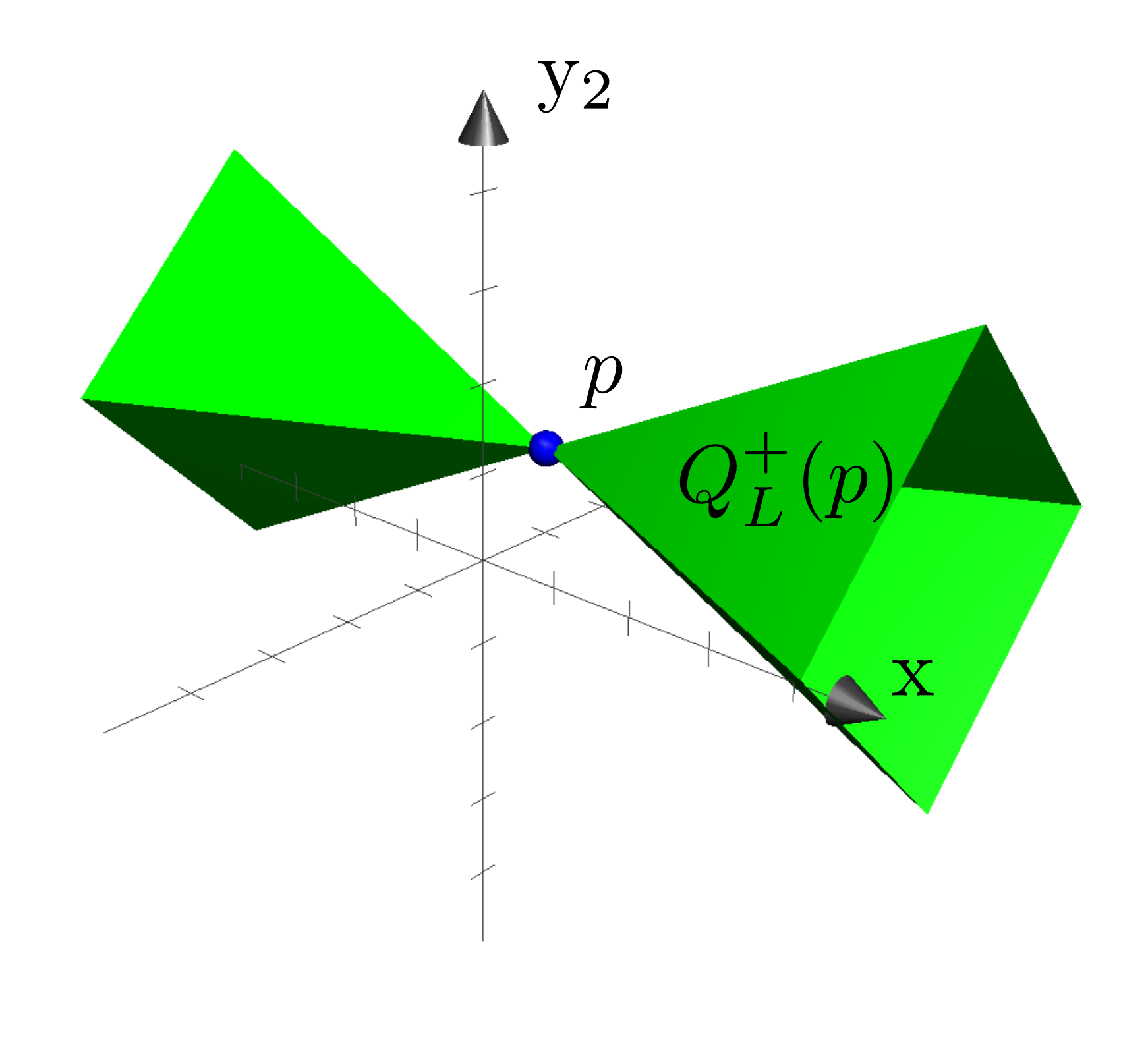}
\end{center}
\caption{A cone attached at $p=(0,1,1)$, in the case when $n_{1}=1$, $%
n_{2}=2 $, $L=\frac{1}{2}$, $\Vert \mathrm{x_{1}}\Vert =|\mathrm{x_{1}}|$
and $\Vert \mathrm{y}\Vert =\Vert (\mathrm{y}_{1},\mathrm{y}_{2})\Vert =|%
\mathrm{y}_{1}|+|\mathrm{y}_{2}|$.}
\label{fig:cone}
\end{figure}

We slightly abuse notations by referring to $Q_{u},Q_{cu},Q_{s}$ and $Q_{cs}$
as cones. We do so since for any point $p\in \mathbb{R}^{u+s+c}$ the sets
\begin{equation*}
Q_{\kappa }^{+}\left( p\right) :=\left\{ q:Q_{\kappa }\left( p-q\right) \geq
0\right\} ,\qquad \kappa \in \left\{ u,cu,s,cs\right\}
\end{equation*}%
are cones centered at $p$. (See Figure \ref{fig:cone}.) This means that $%
Q_{u},Q_{cu},Q_{s}$ and $Q_{cs}$ define cones that can be attached to any
point $p\in \mathbb{R}^{u+s+c}$.

We will assume that $L_{s},L_{u}\in \left( 0,\pi \right) $. We do so for
convenience: We are working in the simplified setting where $\Lambda $ is a
torus $\left(\mathbb{R}/(2\pi)\right)^c$. When $L_{u}\in \left( 0,\pi \right) $ and $p\in S^{cs}$, then the set
$Q_{u}^{+}\left( p\right) \cap S^{cs}$ is contained in a single chart, since
for any $q\in Q_{u}^{+}\left( p\right) $ with $\left\Vert \pi
_{x}q\right\Vert \leq 1$ and $\left\Vert \pi _{y}q\right\Vert \leq 1$ we
will have%
\begin{equation*}
\left\Vert \pi _{\theta }\left( q-p\right) \right\Vert \leq \left\Vert \pi
_{y,\theta }\left( q-p\right) \right\Vert \leq L_{u}\left\Vert \pi
_{x}\left( q-p\right) \right\Vert \leq L_{u}\left( \left\Vert \pi
_{x}q\right\Vert +\left\Vert \pi _{x}p\right\Vert \right) \leq 2L_{u}<2\pi .
\end{equation*}%
A mirror argument can be made that for $p\in S^{cu}$ the set $%
Q_{s}^{+}\left( p\right) \cap S^{cu}$ is also contained in a single chart.

We note though that for a given point $p\in S^{cu}$ the set $%
Q_{cu}^{+}\left( p\right) $ is only locally defined in a neighbourhood of $p$%
, which is small enough to be contained in a single chart. The same is for $%
Q_{cs}^{+}\left( p\right)$.

\begin{remark}
\label{rem:points-in-same-chart}Whenever we write $Q_{cu}\left( p-q\right) $
or $Q_{cs}\left( p-q\right) $ we implicitly assume that $q$ and $p$ are in
some common local chart.
\end{remark}

\begin{definition}
Let $U\subset \mathbb{R}^{u+s+c}$ and let $\kappa \in \left\{ u,cu\right\} $%
. We say that a flow $\Phi _{t}$ satisfies (forward) $Q_{\kappa }$-cone
conditions in $U$ if for every $p_{1},p_{2}\in U$ satisfying $Q_{\kappa
}\left( p_{1}-p_{2}\right) \geq 0$ the fact that $\Phi _{\left[ 0,t\right]
}\left( p_{i}\right) \subset U$, for both $i=1,2$ and some $t>0$, implies
that%
\begin{equation}
Q_{\kappa }\left( \Phi _{t}(p_{1})-\Phi _{t}(p_{2})\right) \geq 0.  \notag
\end{equation}
\end{definition}

\begin{definition}\label{defin:backwardcone}
Let $U\subset \mathbb{R}^{u+s+c}$ and let $\kappa \in \left\{ s,cs\right\} $%
. We say that a flow $\Phi _{t}$ satisfies (backward) $Q_{\kappa }$-cone
conditions in $U$ if for every $p_{1},p_{2}\in U$ satisfying $Q_{\kappa
}\left( p_{1}-p_{2}\right) \geq 0$ the fact that $\Phi _{\left[ -t,0\right]
}\left( p_{i}\right) \subset U$, for both $i=1,2$ and some $t>0$, implies
that%
\begin{equation}
Q_{\kappa }\left( \Phi _{-t}(p_{1})-\Phi _{-t}(p_{2})\right) \geq 0.  \notag
\end{equation}
\end{definition}


\begin{definition}
We say that $\Phi _{t}$ is (forward) outflowing from $S^{cs}$ along $Q_{u}$
if for every $p_{1},p_{2}\in S^{cs}$ satisfying $Q_{u}\left(
p_{1}-p_{2}\right) \geq 0$ there exists a $t>0$ such that%
\begin{equation*}
\Phi _{t}(p_{i})\notin S^{cs}\qquad \text{for some }i\in \left\{ 1,2\right\}
.
\end{equation*}
\end{definition}

\begin{definition}\label{defin:backoutflowing}
We say that $\Phi _{t}$ is (backward) outflowing from $S^{cu}$ along $Q_{s}$
if for every $p_{1},p_{2}\in S^{cu}$ satisfying $Q_{s}\left(
p_{1}-p_{2}\right) \geq 0$ there exists a $t>0$ such that%
\begin{equation*}
\Phi _{-t}(p_{i})\notin S^{cu}\qquad \text{for some }i\in \left\{
1,2\right\} .
\end{equation*}
\end{definition}

Intuitively, if $\Phi_{t}$ satisfies cone conditions then any two points
which are aligned by the cones will flow to points, which are also aligned
by the cones. The outflowing condition states that at least one of two such
points will eventually flow out of the considered set.


\subsection{Construction of stable and unstable manifolds}

The aim of this section is to prove the following two theorems.

\begin{theorem}
\label{th:Wcu-bound}Let $S^{u}$ be a sector (see Definition \ref%
{def:Scu-sector}) and denote by $\Phi _{t}$ the flow induced by $F$. Assume
that:

\begin{enumerate}
\item The flow $\Phi _{t}$ satisfies forward cone conditions for $Q_{cu}$ in 
$S^{u}$ and backward cone conditions for $Q_{s}$ in $S^{u}$.

\item Every forward trajectory starting from a point in the sector $S^{u}$
must exit the sector. 

\item The flow $\Phi _{t}$ is backward outflowing from $S^{u}$ along $Q_{s}$.
\end{enumerate}

Then the unstable manifold $W^{u}(\tilde{\Lambda})$ is contained in $%
S^{u}$. Moreover, $w^{u}$ is Lipschitz, with Lipschitz constant $L_{cu}$.
(The $L_{cu}$ is the constant in the cone $Q_{cu}$; see (\ref{eq:cones})).
\end{theorem}

\begin{theorem}
\label{th:Wcs-bound}Let $S^{s}$ be a sector (see Definition \ref%
{def:Scs-sector}) and denote by $\Phi _{t}$ the flow induced by $F$. Assume
that:

\begin{enumerate}
\item The flow $\Phi _{t}$ satisfies backward cone conditions for $Q_{cs}$
in $S^{s}$ and forward cone conditions for $Q_{u}$ in $S^{s}$.

\item Every backward trajectory starting from a point in the sector $S^{s}$
must exit the sector. 

\item The flow $\Phi _{t}$ is forward outflowing from $S^{s}$ along $Q_{u}$.
\end{enumerate}

Then the stable manifold $W^{s}(\tilde{\Lambda})$ is contained in $%
S^{s}$. Moreover, $w^{s}$ is Lipschitz, with the Lipschitz constant $L_{cs}$%
. (The $L_{cs}$ is the constant in the cone $Q_{cu}$; see (\ref{eq:cones}).)
\end{theorem}

We will focus on proving Theorem \ref{th:Wcu-bound}, since Theorem \ref%
{th:Wcs-bound} follows from Theorem \ref{th:Wcu-bound} by reversing the sign
of the vector field and swapping the roles of the coordinates $x$ and $y$.

Before we prove Theorem \ref{th:Wcu-bound}, we need some additional notions
and technical lemmas.

To simplify the notation, throughout this section let us write here 
\begin{equation*}
\mathrm{x}=(x,\lambda ).
\end{equation*}

\begin{definition}
We say that $h:\pi _{\mathrm{x}}\overline{S^{u}}\rightarrow \overline{S^{u}}$
is a center-horizontal disc satisfying $Q_{cu}$ cone condition if 
\begin{equation}
\pi _{\mathrm{x}}h=\id_{\mathrm{x}},  \label{eq:h-projection-cond}
\end{equation}%
and for every $\mathrm{x}_{1},\mathrm{x}_{2}\in \pi _{\mathrm{x}}\overline{%
S^{u}}$, such that $h\left( \mathrm{x}_{1}\right) $, $h\left( \mathrm{x}%
_{2}\right) $ lie in a single chart, 
\begin{equation}
Q_{cu}\left( h\left( \mathrm{x}_{1}\right) -h\left( \mathrm{x}_{2}\right)
\right) \geq 0.  \label{eq:h-cone-cond}
\end{equation}
\end{definition}

For discs $h$ as defined above we will write%
\begin{equation*}
\mathrm{graph}\left( h\right) :=h\left( \pi _{\mathrm{x}}\overline{S^{u}}%
\right) .
\end{equation*}

\begin{remark}
In our proof of Theorem \ref{th:Wcu-bound} we will show that there exists a
center-horizontal disc $w^{u}$ satisfying $Q_{cu}$ cone condition such that%
\begin{equation*}
W^{u}(\tilde{\Lambda})=\mathrm{graph}\left( w^{u}\right) .
\end{equation*}
\end{remark}

The lemma below will be the main building block for the construction of $%
W^{u}(\tilde{\Lambda})$ in the proof of Theorem \ref{th:Wcu-bound}.

\begin{lemma}
\label{lem:ch-disc-propagation}Assume that $S^{u}$ is an unstable sector and $\Phi _{t}$ satisfies forward cone
conditions for $Q_{cu}$, then there exists $T>0$ such that for every
center-horizontal disc satisfying $Q_{cu}$ cone condition $h:\pi _{\mathrm{x}%
}S^{u}\rightarrow S^{u}$ there exists a center-horizontal disc satisfying $%
Q_{cu}$ cone conditions $h^{\prime }:\pi _{\mathrm{x}}S^{u}\rightarrow S^{u}$
such that%
\begin{equation*}
\Phi _{T}\left( \mathrm{graph}\left( h\right) \right) \cap S^{u}=\mathrm{%
graph}\left( h^{\prime }\right) .
\end{equation*}%
Moreover, if for $q\in \mathrm{graph}\left( h\right) $ we have $\Phi
_{T}\left( q\right) \in S^{u}$, then $q\in \mathrm{graph}\left( h^{\prime
}\right) $ and $\Phi _{t}\left( q\right) \in S^{u}$ for $t\in \left[ 0,T%
\right] $.
\end{lemma}

\begin{proof}
Let $r\in \left( 0,\pi \right) $ be fixed. By the continuity of the flow
with respect to time and initial conditions, for every $q_{1}\in \overline{%
S^{u}}$ there exists a $T>0$ such that for all $t\in \left[ 0,T\right] $ and
all $q_{2}$ such that $\| q_1-q_2\|=r$ we have 
\begin{equation}
\frac{r}{2}<\left\Vert \Phi _{t}\left( q_{2}\right) -q_{1}\right\Vert <\pi .
\label{eq:points-not-far}
\end{equation}
By compactness of $\overline{S^{u}}$ the $T$ can be chosen to be independent
of the choice of $q_{1}$. Since the points in $S_{-}^{u}$ exit the set $%
S^{u}$ (see Condition 1 from Definition \ref{def:Scu-sector}), and $S_{-}^{u}$ is compact, we can choose $T$ small enough so that
in addition to (\ref{eq:points-not-far}),
%
%

\begin{equation}
\Phi _{t}\left( q\right) \notin S^{u},\qquad \text{for every }t\in \left[
0,2T\right] \text{ and }q\in S_{-}^{u}.  \label{eq:points-exit}
\end{equation}%
Condition (\ref{eq:points-exit}) ensures that if we exit the set $S^{u}$
then we can not return to it in a time shorter than $2T$. This implies that
if $q\in S^{u}$ and $\Phi _{T}\left( q\right) \in S^{u}$, then $\Phi
_{t}\left( q\right) \in S^{u}$ for all $t\in \left[ 0,T\right] $.

Let us introduce the following notation. We will write $\mathrm{x}=\left(
x,\lambda \right) $ for a point in $\mathbb{R}^{u}\times \Lambda $. Let us
also introduce the following set 
\begin{equation*}
\mathcal{D}_{u}:=\pi _{\mathrm{x}}S^{u}
\end{equation*}%
Observe that $\partial \mathcal{D}_{u}:=\pi _{\mathrm{x}}\partial S^{u}=\pi
_{\mathrm{x}}S_{-}^{u}$.

For a given center-horizontal disc $h$ satisfying the $Q_{cu}$ cone
condition and fixed $t\in \mathbb{R}$, let us define $g_{t}:\overline{%
\mathcal{D}_{u}}\rightarrow \mathbb{R}^{u}\times \Lambda $%
\begin{equation*}
g_{t}\left( \mathrm{x}\right) =\pi _{\mathrm{x}}\Phi _{t}\circ h(\mathrm{x}).
\end{equation*}%
(The function $g_{t}$ depends on the choice of $h$.) For small $t$ the
function $g_{t}\mathrm{\ }$is close to identity. From (\ref%
{eq:points-not-far}--\ref{eq:points-exit}), for every $t\in \left[ 0,T\right]%
, $ 
\begin{align}
& \left. \frac{r}{2}<\left\Vert g_{t}\left( \mathrm{x}_{2}\right) -\mathrm{x}%
_{1}\right\Vert <\pi \qquad \text{for all }\mathrm{x}_{1}\in \overline{%
\mathcal{D}_{u}}\text{ and }\mathrm{x}_{2}\in \partial B\left( \mathrm{x}%
_{1},r\right) \cap \overline{\mathcal{D}_{u}},\right.
\label{eq:points-on-boundary-of-U-1} \\
& \left. g_{t}\left( \mathrm{x}_{3}\right) \notin \mathcal{D}_{u}\qquad 
\text{for all }\mathrm{x}_{3}\in \partial \mathcal{D}_{u}.\right.
\label{eq:points-on-boundary-of-U-2}
\end{align}%
Note that the choice of $T$ is independent from $h$.


We will show that for every $\mathrm{x}_{1}\in\mathcal{D}_{u}$ there exists
a point $\mathrm{x}_{1}^{\prime}\in\mathcal{D}_{u}$ such that $g_{T}\left( 
\mathrm{x}_{1}^{\prime}\right) =\mathrm{x}_{1}$. We will prove this by using
the local Brouwer degree (see Appendix). To do so, we will construct a
homotopy from $g_{T}$ to the identity map with some good properties. Let $%
U=B\left( \mathrm{x}_{1},r\right) \cap\mathcal{D}_{u}$. 
Let $H:\left[ 0,1\right] \times U\rightarrow\mathbb{R}^{u+c}$ be a homotopy
chosen as 
\begin{equation*}
H\left( \alpha,\mathrm{x}\right) =g_{\alpha T}\left( \mathrm{x}\right) .
\end{equation*}
Note that 
\begin{equation*}
H\left( 0,\mathrm{x}\right) =\mathrm{x},\qquad H\left( 1,\mathrm{x}\right)
=g_{T}\left( \mathrm{x}\right) .
\end{equation*}
From (\ref{eq:points-on-boundary-of-U-1}--\ref{eq:points-on-boundary-of-U-2}%
) we see that for $\mathrm{x}_{2}\in\partial B\left( \mathrm{x}_{1},r\right)
\cap\overline{\mathcal{D}_{u}}$ and $\mathrm{x}_{3}\in\partial\mathcal{D}%
_{cu}$%
\begin{equation}
H\left( \alpha,\mathrm{x}_{i}\right) \neq\mathrm{x}_{1},\qquad\text{for }%
i=2,3.  \label{eq:x1-not-on-bd}
\end{equation}
Since $\partial U=\left( \partial B\left( \mathrm{x}_{1},r\right) \cap%
\overline{\mathcal{D}_{u}}\right) \cup(\overline{B\left( \mathrm{x}%
_{1},r\right) }\cap\partial\mathcal{D}_{u})$, (\ref{eq:x1-not-on-bd})
implies that for every $\mathrm{x}\in\partial U$ and every $\alpha\in\left[
0,1\right] $%
\begin{equation*}
\mathrm{x}_{1}\notin H\left( \left[ 0,1\right] ,\partial U\right) .
\end{equation*}
By the homotopy property of the Brouwer degree%
\begin{align*}
\deg\left( g_{T},U,\mathrm{x}_{1}\right) & =\deg\left( H\left(
1,\cdot\right) ,U,\mathrm{x}_{1}\right) =\deg\left( H\left( 0,\cdot \right)
,U,\mathrm{x}_{1}\right) \\
& =\deg\left( \id,U,\mathrm{x}_{1}\right) =1,
\end{align*}
hence by the solution property of the Brouwer degree we see that there
exists an $\mathrm{x}_{1}^{\prime}\in U\subset\mathcal{D}_{u}$ such that 
\begin{equation*}
g_{T}\left( \mathrm{x}_{1}^{\prime}\right) =\mathrm{x}_{1}.
\end{equation*}

Above we have shown that $\mathcal{D}_{u}\subset g_{T}\left( \mathcal{D}%
_{cu}\right) ,$ hence from the continuity of $g_{T}$ and compactness of $%
\overline{\mathcal{D}_{u}}$, $g_{T}(\overline{\mathcal{D}_{u}})$ is
compact and $\overline{\mathcal{D}_{u}}\subset g_{T}(\overline{\mathcal{D}%
_{cu}})$.

Now we will show that $g_{T}:\overline{\mathcal{D}_{u}}\rightarrow \mathbb{R%
}^{u}\times \Lambda $ is injective. By choosing small $T$ the map $g_{T}$ is
close to identity. It is therefore enough to show that for $\mathrm{x}%
_{1}\neq \mathrm{x}_{2}$ close enough so that $h(\mathrm{x}_{1})-h(\mathrm{x}%
_{2})$ are in the same chart, we will have $g_{T}\left( \mathrm{x}%
_{1}\right) \neq $ $g_{T}\left( \mathrm{x}_{2}\right).$ For $\mathrm{x}%
_{1}\neq \mathrm{x}_{2}$ we know that 
\begin{equation*}
Q_{cu}\left( h(\mathrm{x}_{1})-h(\mathrm{x}_{2})\right) \geq 0,
\end{equation*}%
so, by the fact that $\Phi _{t}$ satisfies cone conditions for $Q_{cu}$ we
also have 
\begin{equation*}
0\leq Q_{cu}\left( \Phi _{T}\left( h(\mathrm{x}_{1})\right) -\Phi _{T}\left(
h(\mathrm{x}_{2})\right) \right) ,
\end{equation*}%
which implies%
\begin{align}
\left\Vert \pi _{y}\left[ \Phi _{T}\left( h(\mathrm{x}_{1})\right) -\Phi
_{T}\left( h(\mathrm{x}_{2})\right) \right] \right\Vert & \leq
L_{cu}\left\Vert \pi _{x}\left[ \Phi _{T}\left( h(\mathrm{x}_{1})\right)
-\Phi _{T}\left( h(\mathrm{x}_{2})\right) \right] \right\Vert
\label{eq:h-prime-hor} \\
& =L_{cu}\left\Vert g_{T}\left( \mathrm{x}_{1}\right) -g_{T}\left( \mathrm{x}%
_{2}\right) \right\Vert .  \notag
\end{align}%
If $\pi _{y}\left[ \Phi _{T}\left( h(\mathrm{x}_{1})\right) -\Phi _{T}\left(
h(\mathrm{x}_{2})\right) \right] \neq 0,$ then above implies that $%
g_{T}\left( \mathrm{x}_{1}\right) \neq g_{T}\left( \mathrm{x}_{2}\right) $.
If $\pi _{y}\left[ \Phi _{T}\left( h(\mathrm{x}_{1})\right) -\Phi _{T}\left(
h(\mathrm{x}_{2})\right) \right] =0$ then we see that since by uniqueness of
solutions of ODEs $\Phi _{T}\left( h(\mathrm{x}_{1})\right) \neq \Phi
_{T}\left( h(\mathrm{x}_{2})\right) $ we must also have $g_{T}\left( \mathrm{%
x}_{1}\right) \neq g_{T}\left( \mathrm{x}_{2}\right) $.

Above argument shows that also $g_{t}$ is injective, for every $t\in \left[
0,T\right] .$

Let us now define%
\begin{equation}
h^{\prime }\left( \mathrm{x}\right) :=\Phi _{T}\left( h(g_{T}^{-1}\left( 
\mathrm{x}\right) )\right) .  \label{eq:def-h'}
\end{equation}%
By definition of $h^{\prime }$ we see that $\pi _{\mathrm{x}}h^{\prime
}\left( \mathrm{x}\right) =\pi _{\mathrm{x}}\Phi _{T}\left(
h(g_{T}^{-1}\left( \mathrm{x}\right) )\right) =g_{T}\left( g_{T}^{-1}\left( 
\mathrm{x}\right) \right) =\mathrm{x}$. From (\ref{eq:h-prime-hor}) we see
that%
\begin{align*}
\left\Vert \pi _{y}\left[ h^{\prime }\left( \mathrm{x}_{1}\right) -h^{\prime
}\left( \mathrm{x}_{2}\right) \right] \right\Vert & =\left\Vert \pi _{y}%
\left[ \Phi _{T}\left( h(g_{T}^{-1}\left( \mathrm{x}_{1}\right) )\right)
-\Phi _{T}\left( h(g_{T}^{-1}\left( \mathrm{x}_{2}\right) )\right) \right]
\right\Vert \\
& \leq L_{cu}\left\Vert \pi _{x}\left[ \Phi _{T}\left( h(g_{T}^{-1}\left( 
\mathrm{x}_{1}\right) )\right) -\Phi _{T}\left( h(g_{T}^{-1}\left( \mathrm{x}%
_{2}\right) )\right) \right] \right\Vert \\
& =L_{cu}\left\Vert \mathrm{x}_{1}-\mathrm{x}_{2}\right\Vert.
\end{align*}%
Hence $Q_{cu}\left( h^{\prime }\left( \mathrm{x}_{1}\right) -h^{\prime
}\left( \mathrm{x}_{2}\right) \right) \geq 0$, so $h^{\prime }$ so satisfies 
$Q_{cu}$ cone condition.

What remains is to show that $\mathrm{graph}(h^{\prime})\subset \overline{%
S^{u}}.\ $Points can exit $\overline{S^{u}}$ along a forward trajectory only
through the set $S_{-}^{u}$. Once they exit, due to our choice of $T$ they
can not come back to $\overline{S^{u}}$. Moreover, $g_{T}\left( \mathrm{x}%
\right) $ is injective. This means that $h^{\prime }$ contains only points $%
\Phi _{T}\circ h(\mathrm{x})$ for which $\Phi _{\left[ 0,T\right] }\circ h(%
\mathrm{x})\subset \overline{S^{u}},$ hence $h^{\prime }\subset \overline{%
S^{u}}.$ As required.
\end{proof}

We are now ready to prove Theorem \ref{th:Wcu-bound}.

\begin{proof}[Proof of Theorem \protect\ref{th:Wcu-bound}]
Let $T$ be as obtained in Lemma~\ref{lem:ch-disc-propagation}.

For a given center-horizontal disc $h$ satisfying $Q_{cu}$ cone conditions,
by Lemma \ref{lem:ch-disc-propagation} we know that there exists the disc $%
h^{\prime }$. Let us introduce the notation $\mathcal{G}\left( h\right)
:=h^{\prime }$.

We shall now construct $w^{u}.$ The idea is to take $h_{0}\left( \mathrm{x}%
\right) =\left( \pi _{x}\mathrm{x},0,\pi _{\lambda }\mathrm{x}\right) $,
inductively define $h_{n+1}:=\mathcal{G}\left( h_{n}\right) $ for $n\geq 0$,
and show that $h_{n+1}$ converges to $w^{u}.$

Consider $h_{n}$ as defined above. Let $\mathrm{x}\in \pi_{\mathrm{x}} S^{u}$
be fixed. By compactness of $\{ p\in \overline{S^{u}}:$$\pi_{\mathrm{x}}p=%
\mathrm{x}\} $ and the fact that $\pi _{\mathrm{x}}h_{n}\left( \mathrm{x}%
\right) =\mathrm{x}$, there exists a convergent subsequence $%
\lim_{k\rightarrow \infty }h_{n_{k}}\left( \mathrm{x}\right) =q.$ We will
show that such $q$ has to be unique. Suppose that for two subsequences $%
m_{k} $ and $n_{k}$ we have $\lim_{k\rightarrow \infty }h_{m_{k}}\left( 
\mathrm{x}\right) =q_{1}$ and $\lim_{k\rightarrow \infty }h_{n_{k}}\left( 
\mathrm{x}\right) =q_{2}$.

Let us fix $t>0$. We have $\Phi _{-t}(q_{1})=\lim_{k\rightarrow \infty }\Phi
_{-t}(h_{m_{k}}\left( \mathrm{x}\right) )$ and $\Phi _{-t}(h_{m_{k}}\left( 
\mathrm{x}\right) )\in S^{u}$ for $m_{k}$ large enough ($t<m_{k}\cdot T$).
Hence $\Phi _{-t}(q_{1})\in S^{u}$. For $q_{2}$, by analogous argument, we
also obtain $\Phi _{-t}(q_{2})\in S^{u}$.

Since $\pi _{\mathrm{x}}q_{1}=\mathrm{x}=\pi _{\mathrm{x}}q_{2}$, we see
that $Q_{s}(q_{1}-q_{2})\geq 0$, and since the $\Phi _{-t}$ is outflowing
along $Q_{s}$ we see that $q_{1}$ has to be equal to $q_{2}$, otherwise one
of them would exit $S^{u}$.

We now define $w^{u}\left( \mathrm{x}\right) =\lim_{n\rightarrow \infty
}h_{n}\left( \mathrm{x}\right) $. The conditions (\ref{eq:h-projection-cond}%
--\ref{eq:h-cone-cond}) are preserved by passing to the limit, which
concludes the construction of $w^{u}$. In particular, the graph of $w^{u}$
is $L_{cu}$--Lipschitz.

By construction, the center horizontal disc $\mathrm{graph}\left(
w^{u}\right) $ consists of points, whose backward trajectories remain in $%
S^{u}$.

By repeating the above argument leading to $q_{1}=q_{2}$ we can easily prove
that any point whose backward trajectory remains in $S^{u}$ has to be in $%
\mathrm{graph}\left( w^{u}\right) $. Moreover, since the forward trajectory
starting from every point from $S^{u}$ must exit the sector, a backward
trajectory which remains in $S^{u}$ must accumulate on a limit set contained
in the boundary $\partial S^{u}$. Note that $\partial S^{u}=S^{u}_-\cup
S^{u}_+\cup \tilde\Lambda$. All points from $S^{u}_-$ exit $\overline{S^{u}}$
and all points from $S^{u}_+$ enter $S^{u}$. This means that any backward
trajectory which remains in $S^{u}$ must converge to $\tilde \Lambda$.

This concludes the proof.
\end{proof}


\subsection{Validation of cone and outflowing conditions based on
contraction\-/ex\-pan\-sion rates\label{sec:rates-verification}}

In this section we introduce `rates' of contraction and expansion associated
to a matrix and show how they can be used to validate cone conditions and
outflowing conditions. The conditions follow from estimates on the matrices $%
G$ appearing in (\ref{eq:factorised-derivative}). The tools presented in this
section make the abstract Theorems \ref{th:Wcu-bound} and \ref{th:Wcs-bound}
a practical tool for establishing the existence of the invariant manifolds.

Once again, we focus on the case of the unstable manifold, since the
stable manifold follows from changing the sign of the vector field
and swapping the roles of the coordinates $x,y$.

Throughout this section we keep the notation
\begin{equation*}
\mathrm{x}=\left( x,\lambda \right) .
\end{equation*}

Assume that (\ref{eq:factorised-derivative}) is satisfied, i.e. that%
\begin{equation*}
DF\left( \mathrm{x},y\right) =h(\mathrm{x},y)G\left( \mathrm{x},y\right)
,\qquad \text{for }(\mathrm{x},y)\in S^{u},
\end{equation*}%
where
\begin{equation}
h\left( \mathrm{x},y\right) >0\qquad \text{for }(\mathrm{x},y)\in S^{u}.
\label{eq:positive-h-2}
\end{equation}

Consider the matrix $G$ of the form%
\begin{equation*}
G(p)=\left(
\begin{array}{cc}
G_{\mathrm{xx}}(p) & G_{\mathrm{x}y}(p) \\
G_{y\mathrm{x}}(p) & G_{yy}(p)%
\end{array}%
\right) \qquad \text{for }p\in \overline{S^{u}},
\end{equation*}%
where $G_{\mathrm{xx}},$ $G_{\mathrm{x}y},$ $G_{y\mathrm{x}}$ and $G_{yy}$
are $\left( c+u\right) \times \left( c+u\right) $, $\left( c+u\right) \times
s,$ $s\times \left( c+u\right) $ and $s\times s$ matrices, respectively. Let
us define the following constants (see \eqref{def:matrixconst}),%
\begin{equation*}
\begin{array}{lll}
m_{l}\left( G_{\mathrm{xx}}\right) :=\inf_{p\in \overline{S^{u}}%
}m_{l}\left( G_{\mathrm{xx}}(p)\right) , &  & l\left( G_{yy}\right)
:=\sup_{p\in \overline{S^{u}}}l\left( G_{yy}(p)\right) ,\medskip \\
\left\Vert G_{\mathrm{x}y}\right\Vert :=\sup_{p\in \overline{S^{u}}%
}\left\Vert G_{\mathrm{x}y}\left( p\right) \right\Vert , &  & \left\Vert G_{y%
\mathrm{x}}\right\Vert :=\sup_{p\in \overline{S^{u}}}\left\Vert G_{y\mathrm{%
x}}\left( p\right) \right\Vert .%
\end{array}%
\end{equation*}

\begin{definition}
Let $\xi _{cu},\mu _{s}\in \mathbb{R}$ be defined as
\begin{equation}\label{eq:xi-def} 
\begin{aligned}
\xi _{cu}& :=m_{l}\left( G_{\mathrm{xx}}\right) -L_{cu}\left\Vert G_{\mathrm{%
x}y}\right\Vert ,  
\\
\mu _{s}& :=l\left( G_{yy}\right) +\frac{1}{L_{cu}}\left\Vert G_{y\mathrm{x}%
}\right\Vert . 
\end{aligned}
\end{equation}
We refer to $\xi _{cu}$ as the \emph{expansion rate} of $G$ and to $\mu _{s}$
as the \emph{contraction rate} of $G$.
\end{definition}

The lemma below provides a tool for validating cone conditions based on the
expansion and contraction rates.

\begin{lemma}
\label{lem:cone-verif}If the constants $\xi _{cu},\mu _{s}$ defined by \eqref{eq:xi-def} have the property
\begin{equation}
\mu _{s}<\xi _{cu},  \label{eq:rate-cond}
\end{equation}%
then the flow induced by (\ref{eq:our-ode}) satisfies $Q_{cu}$ cone
conditions in $S^{u}$.  
\end{lemma}

\begin{proof}
By Lemma 8 from \cite{MR3567489} we know that $m(I+tA)=1+tm_{l}(A)+O(t^{2})$
and the bound $O(t^{2})$ is uniform if we consider matrices $A$ in some
compact set. In our case this compact set is given as $\{DF(p) , p \in \overline{S^{u}}\}$.

 Let $p_{1}\neq p_{2}$, $p_{i}=\left( \mathrm{x}%
_{i},y_{i}\right) $ for $i=1,2$, be such that $Q_{cu}\left(
p_{1}-p_{2}\right) \geq 0.$ Since $\left\Vert y_{1}-y_{2}\right\Vert \leq
L_{cu}\left\Vert \mathrm{x}_{1}-\mathrm{x}_{2}\right\Vert $ and $h\geq 0$, 
by \eqref{def:matrixconst} and
Lemma \ref{lem:norm-integral-2},\textbf{\ }for $t>0$%


\begin{equation}\label{eq:xi-expansion-temp} 
\begin{aligned}
\lefteqn{ \frac{1}{t}\left( \left\Vert \pi _{\mathrm{x}}\left( \Phi _{t}\left(
p_{1}\right) -\Phi _{t}\left( p_{2}\right) \right) \right\Vert -\left\Vert
\mathrm{x}_{1}-\mathrm{x}_{2}\right\Vert \right)}\\
& =\frac{1}{t}\left( \left\Vert \left( \mathrm{x}_{1}-\mathrm{x}_{2}\right)
+t\pi _{\mathrm{x}}\left( F\left( p_{1}\right) -F\left( p_{2}\right) \right)
+O\left( t^{2}\right) \right\Vert -\left\Vert \mathrm{x}_{1}-\mathrm{x}%
_{2}\right\Vert \right)   \\
& =\left\Vert \frac{1}{t}\left( \id_{\mathrm{x}}+t\int_{0}^{1}\pi _{\mathrm{x}%
}\frac{\partial F}{\partial \mathrm{x}}\left( p_{1}+s\left(
p_{1}-p_{2}\right) \right) ds\right) \left( \mathrm{x}_{1}-\mathrm{x}%
_{2}\right) \right.   \\
& \quad \left. +\int_{0}^{1}\pi _{\mathrm{x}}\frac{\partial F}{\partial y}%
\left( p_{1}+s\left( p_{1}-p_{2}\right) \right) ds\left( y_{1}-y_{2}\right)
\right\Vert -\frac{1}{t}\left\Vert \mathrm{x}_{1}-\mathrm{x}_{2}\right\Vert
+O\left( t\right) \\
& \geq \frac{1}{t}\left( m\left( \id_{\mathrm{x}}+t\int_{0}^{1}\pi _{\mathrm{x%
}}\frac{\partial F}{\partial \mathrm{x}}\left( p_{1}+s\left(
p_{1}-p_{2}\right) \right) ds\right) -1\right) \left\Vert \mathrm{x}_{1}-%
\mathrm{x}_{2}\right\Vert \\
& \quad -\int_{0}^{1}\left\Vert \pi _{\mathrm{x}}\frac{\partial F}{\partial y%
}\left( p_{1}+s\left( p_{1}-p_{2}\right) \right) \right\Vert
dsL_{cu}\left\Vert \mathrm{x}_{1}-\mathrm{x}_{2}\right\Vert +O\left( t\right)\\
& =m_{l}\left( \int_{0}^{1}\pi _{\mathrm{x}}\frac{\partial F}{\partial
\mathrm{x}}\left( p_{1}+s\left( p_{1}-p_{2}\right) \right) ds\right)
\left\Vert \mathrm{x}_{1}-\mathrm{x}_{2}\right\Vert +O\left( t\right) \\
& \quad -L_{cu}\int_{0}^{1}\left\Vert \pi _{\mathrm{x}}\frac{\partial F}{%
\partial y}\left( p_{1}+s\left( p_{1}-p_{2}\right) \right) \right\Vert
ds\left\Vert \mathrm{x}_{1}-\mathrm{x}_{2}\right\Vert +O\left( t\right) \\
& =m_{l}\left( \int_{0}^{1}(h\cdot G_{\mathrm{xx}})\left( p_{1}+s\left( p_{1}-p_{2}\right) \right) ds\right)
\left\Vert \mathrm{x}_{1}-\mathrm{x}_{2}\right\Vert +O\left( t\right) \\
& \quad -L_{cu}\int_{0}^{1}\left\Vert (h\cdot G_{\mathrm{x}y})\left( p_{1}+s\left( p_{1}-p_{2}\right) \right) \right\Vert
ds\left\Vert \mathrm{x}_{1}-\mathrm{x}_{2}\right\Vert +O\left( t\right)\\
& \geq \left( m_{l}\left( G_{\mathrm{xx}}\right) -L_{cu}\left\Vert G_{%
\mathrm{x}y}\right\Vert \right) \int_{0}^{1}h\left( p_{1}+s\left(
p_{1}-p_{2}\right) \right) ds\left\Vert \mathrm{x}_{1}-\mathrm{x}%
_{2}\right\Vert +O\left( t\right) \\
& = \xi _{cu}\int_{0}^{1}h\left( p_{1}+s\left( p_{1}-p_{2}\right) \right)
ds\left\Vert \mathrm{x}_{1}-\mathrm{x}_{2}\right\Vert +O\left( t\right).
\end{aligned}%
\end{equation}
So, letting $t\to 0$,%
\begin{equation}
D_{-}\left\Vert \pi _{\mathrm{x}}\left( \Phi _{t}\left( p_{1}\right) -\Phi
_{t}\left( p_{2}\right) \right) \right\Vert |_{t=0}\geq \xi
_{cu}\int_{0}^{1}h\left( p_{1}+s\left( p_{1}-p_{2}\right) \right)
ds\left\Vert \mathrm{x}_{1}-\mathrm{x}_{2}\right\Vert,
\label{eq:change-in-x}
\end{equation}%
where
\begin{equation*}
D_{-}f(t_0)=\liminf_{t\rightarrow 0^{+}}\frac{f(t_0+t)-f(t_0)}{t}
\end{equation*}%
is the lower Dini derivative of $f$.

By Lemma 7 from \cite{MR3567489}, we know that $\left\Vert I+tA\right\Vert
=1+tl(A)+O(t^{2})$ and the bound $O(t^{2})$ is uniform is we consider
matrices $A$ in some compact set. By using the fact that $h\geq 0$
together with Lemma \ref{lem:norm-integral-1}, for $t>0$,%
\begin{align*}
& \frac{1}{t}\left( \left\Vert \pi _{y}\left( \Phi _{t}\left( p_{1}\right)
-\Phi _{t}\left( p_{2}\right) \right) \right\Vert -\left\Vert
y_{1}-y_{2}\right\Vert \right) \\
& =\frac{1}{t}\left( \left\Vert \left( y_{1}-y_{2}\right) +t\pi _{y}\left(
F\left( p_{1}\right) -F\left( p_{2}\right) \right) +O\left( t^{2}\right)
\right\Vert -\left\Vert y_{1}-y_{2}\right\Vert \right) \\
& =\left\Vert \frac{1}{t}\left( \id_{y}+t\int_{0}^{1}\pi _{y}\frac{\partial F%
}{\partial y}\left( p_{1}+s\left( p_{1}-p_{2}\right) \right) ds\right)
\left( y_{1}-y_{2}\right) \right. \\
& \quad \left. +\int_{0}^{1}\pi _{y}\frac{\partial F}{\partial \mathrm{x}}%
\left( p_{1}+s\left( p_{1}-p_{2}\right) \right) ds\left( \mathrm{x}_{1}-%
\mathrm{x}_{2}\right) \right\Vert -\frac{1}{t}\left\Vert
y_{1}-y_{2}\right\Vert +O\left( t\right) \\
& \leq \frac{1}{t}\left( \left\Vert \left( \id_{y}+t\int_{0}^{1}\pi _{y}\frac{%
\partial F}{\partial y}\left( p_{1}+s\left( p_{1}-p_{2}\right) \right)
ds\right) \right\Vert -1\right) L_{cu}\left\Vert \mathrm{x}_{1}-\mathrm{x}%
_{2}\right\Vert \\
& \quad +\int_{0}^{1}\left\Vert \pi _{y}\frac{\partial F}{\partial \mathrm{x}%
}\left( p_{1}+s\left( p_{1}-p_{2}\right) \right) \right\Vert ds\left\Vert
\mathrm{x}_{1}-\mathrm{x}_{2}\right\Vert +O\left( t\right) \\
& =l\left( \int_{0}^{1}\pi _{y}\frac{\partial F}{\partial y}\left(
p_{1}+s\left( p_{1}-p_{2}\right) \right) ds\right) L_{cu}\left\Vert \mathrm{x%
}_{1}-\mathrm{x}_{2}\right\Vert +O\left( t\right) \\
& \quad +\int_{0}^{1}\left\Vert \pi _{y}\frac{\partial F}{\partial \mathrm{x}%
}\left( p_{1}+s\left( p_{1}-p_{2}\right) \right) \right\Vert ds\left\Vert
\mathrm{x}_{1}-\mathrm{x}_{2}\right\Vert +O\left( t\right) \\
& \leq L_{cu}\left[ l\left( G_{yy}\right) +\frac{1}{L_{cu}}\left\Vert G_{y%
\mathrm{x}}\right\Vert \right] \left\Vert \mathrm{x}_{1}-\mathrm{x}%
_{2}\right\Vert \int_{0}^{1}h\left( p_{1}+s\left( p_{1}-p_{2}\right) \right)
ds+O\left( t\right) \\
&= L_{cu}\mu _{s}\int_{0}^{1}h\left( p_{1}+s\left( p_{1}-p_{2}\right)
\right) ds\left\Vert \mathrm{x}_{1}-\mathrm{x}_{2}\right\Vert +O\left(
t\right).
\end{align*}%
So, letting $t\to 0$,%
\begin{equation}
D^{+}\left\Vert \pi _{y}\left( \Phi _{t}\left( p_{1}\right) -\Phi _{t}\left(
p_{2}\right) \right) \right\Vert |_{t=0}\leq L_{cu}\mu
_{s}\int_{0}^{1}h\left( p_{1}+s\left( p_{1}-p_{2}\right) \right)
ds\left\Vert \mathrm{x}_{1}-\mathrm{x}_{2}\right\Vert ,
\label{eq:change-in-y}
\end{equation}%
where
\begin{equation*}
D_{+}f(t_0)=\limsup_{t\rightarrow 0^{+}}\frac{f(t_0+t)-f(t_0)}{t}
\end{equation*}%
is the upper Dini derivative of $f$.

These estimates of the Dini derivatives imply
\begin{equation}\label{def:DiniBound}
D_{-}\left( L_{cu}\Vert \mathrm{x}_{1}(t)-\mathrm{x}_{2}(t)\Vert -\Vert
y_{1}(t)-y_{2}(t)\Vert \right) _{|t=0}>0.
\end{equation}%
Indeed, using (\ref{eq:positive-h-2}), (\ref%
{eq:rate-cond}), (\ref{eq:change-in-x}), (\ref%
{eq:change-in-y}) we obtain
\begin{align*}
& D_{-}\left( L_{cu}\Vert \mathrm{x}_{1}(t)-\mathrm{x}_{2}(t)\Vert -\Vert
y_{1}(t)-y_{2}(t)\Vert \right) _{|t=0} \\
& \geq L_{cu}D_{-}\left\Vert \pi _{\mathrm{x}}\left( \Phi _{t}\left(
p_{1}\right) -\Phi _{t}\left( p_{2}\right) \right) \right\Vert
|_{t=0}-D^{+}\left\Vert \pi _{y}\left( \Phi _{t}\left( p_{1}\right) -\Phi
_{t}\left( p_{2}\right) \right) \right\Vert |_{t=0} \\
& \geq L_{cu}\xi _{cu}\int_{0}^{1}h\left( p_{1}+s\left( p_{1}-p_{2}\right)
\right) ds\left\Vert \mathrm{x}_{1}-\mathrm{x}_{2}\right\Vert \\
& \quad -L_{cu}\mu _{s}\int_{0}^{1}h\left( p_{1}+s\left( p_{1}-p_{2}\right)
\right) ds\left\Vert \mathrm{x}_{1}-\mathrm{x}_{2}\right\Vert \\
& =(\xi _{cu}-\mu _{s})L_{cu}\int_{0}^{1}h\left( p_{1}+s\left(
p_{1}-p_{2}\right) \right) ds\left\Vert \mathrm{x}_{1}-\mathrm{x}%
_{2}\right\Vert \\
& >0.
\end{align*}
Finally, if we assume that $L_{cu}\Vert \mathrm{x}_{1}-\mathrm{x}_{2}\Vert -\Vert
y_{1}-y_{2}\Vert \geq 0$ (i.e. we are also possibly on the boundary of the
cone), the inequality \eqref{def:DiniBound} implies $Q_{cu}$ forward cone conditions.
We have thus proven that the flow satisfies $Q_{cu}$ cone conditions  on $S^{u}$.
\end{proof}

We now discuss how to validate backward $Q_{s}$ cone condition (see Definition \ref{defin:backwardcone}) and
outflowing from $S^{u}$ condition (see Definition \ref{defin:backoutflowing}). To this end, we define the following two constants%
\[
\begin{split}
\xi _{s} &:=m_{l}\left( -G_{yy}\right) -L_{s}\left\Vert G_{y\mathrm{x}%
}\right\Vert , \\
\mu _{cu} &:=l\left( -G_{\mathrm{xx}}\right) +\frac{1}{L_{s}}\left\Vert G_{%
\mathrm{x}y}\right\Vert .
\end{split}
\]
\begin{lemma}
\label{lem:outflow-verif}If
\begin{equation*}
\mu _{cu}<\xi _{s}\qquad \text{and \qquad }\xi _{s}>0,
\end{equation*}%
then we have backward cone condition for $Q_{s}$ in $S^{u}$ and the
backward outflowing condition from $S^{u}$ along $Q_{s}$.
\end{lemma}

\begin{proof}
If we reverse the sign in the vector field (which then induces the flow with
reversed time), then $\xi _{s}$ plays the role of an expansion rate, and $%
\mu _{cu}$ the role of the contraction rate. This means that from Lemma \ref%
{lem:cone-verif} we obtain backward cone conditions for $Q_{s}$ in $S^{u}$.

We now turn to proving the outflowing from $S^{u}$ along $Q_{s}$ condition. From a
mirror derivation to (\ref{eq:xi-expansion-temp}) for $p_{1}=\left( \mathrm{x%
}_{1},y_{1}\right) $ and $p_{2}=\left( \mathrm{x}_{2},y_{2}\right) $ 
such that $\left\Vert \mathrm{x}_{1}-\mathrm{x}_{2}\right\Vert \leq
L_{s}\left\Vert y_{1}-y_{2}\right\Vert $ and for $t>0$ we obtain
\begin{align*}
& \frac{1}{t}\left( \left\Vert \pi _{y}\left( \Phi _{-t}\left( p_{1}\right)
-\Phi _{-t}\left( p_{2}\right) \right) \right\Vert -\left\Vert
y_{1}-y_{2}\right\Vert \right)  \\
& =\frac{1}{t}\left( \left\Vert \left( y_{1}-y_{2}\right) +t\pi _{y}\left(
-F\left( p_{1}\right) +F\left( p_{2}\right) \right) +O\left( t^{2}\right)
\right\Vert -\left\Vert y_{1}-y_{2}\right\Vert \right)  \\
& =\left\Vert \frac{1}{t}\left( \id_{y}+t\int_{0}^{1}\pi _{y}\frac{-\partial F%
}{\partial y}\left( p_{1}+s\left( p_{1}-p_{2}\right) \right) ds\right)
\left( y_{1}-y_{2}\right) \right.  \\
& \quad \left. +\int_{0}^{1}\pi _{y}\frac{-\partial F}{\partial \mathrm{x}}%
\left( p_{1}+s\left( p_{1}-p_{2}\right) \right) ds\left( \mathrm{x}_{1}-%
\mathrm{x}_{2}\right) \right\Vert -\frac{1}{t}\left\Vert
y_{1}-y_{2}\right\Vert +O\left( t\right)  \\
& \geq \frac{1}{t}\left( m\left( \id_{y}+t\int_{0}^{1}\pi _{y}\frac{-\partial
F}{\partial y}\left( p_{1}+s\left( p_{1}-p_{2}\right) \right) ds\right)
-1\right) \left\Vert y_{1}-y_{2}\right\Vert  \\
& \quad -\int_{0}^{1}\left\Vert \pi _{y}\frac{-\partial F}{\partial \mathrm{x%
}}\left( p_{1}+s\left( p_{1}-p_{2}\right) \right) \right\Vert
dsL_{s}\left\Vert y_{1}-y_{2}\right\Vert +O\left( t\right)  \\
& =m_{l}\left( \int_{0}^{1}\pi _{y}\frac{-\partial F}{\partial y}\left(
p_{1}+s\left( p_{1}-p_{2}\right) \right) ds\right) \left\Vert
y_{1}-y_{2}\right\Vert +O\left( t\right)  \\
& \quad -L_{s}\int_{0}^{1}\left\Vert \pi _{y}\frac{-\partial F}{\partial
\mathrm{x}}\left( p_{1}+s\left( p_{1}-p_{2}\right) \right) \right\Vert
ds\left\Vert y_{1}-y_{2}\right\Vert +O\left( t\right)  \\
& \geq \left( m_{l}\left( -G_{yy}\right) -L_{s}\left\Vert G_{y\mathrm{x}%
}\right\Vert \right) \int_{0}^{1}h\left( p_{1}+s\left( p_{1}-p_{2}\right)
\right) ds\left\Vert y_{1}-y_{2}\right\Vert +O\left( t\right)  \\
& = \xi _{s}\int_{0}^{1}h\left( p_{1}+s\left( p_{1}-p_{2}\right) \right)
ds\left\Vert y_{1}-y_{2}\right\Vert +O\left( t\right) .
\end{align*}%


Now, since for $p\in S^{u}$ we have $h(p)>0$ and $\xi _{s}>0$, for
sufficiently small $t>0$ we obtain
\begin{equation}
\left\Vert \pi _{y}\left( \Phi _{-t}\left( p_{1}\right) -\Phi _{-t}\left(
p_{2}\right) \right) \right\Vert >\left\Vert y_{1}-y_{2}\right\Vert .
\label{eq:change-in-y-2}
\end{equation}
Recall that we are assuming that  $p_{1}\neq p_{2}$ satisfy $Q_{s}\left( p_{1}-p_{2}\right) \geq 0$. Up to now, we have proven that, for $t>0$ and  as long as $\Phi _{-t}(p_{1}),\Phi _{-t}(p_{2})\in S^{u}$, the following is satisfied:
\begin{itemize}
 \item Since $Q_{s}^{+}\left( p\right) \cap S^{u}
$ is in a single chart, for every $p\in S^{u}$ and $t>0$, 
\[Q_{s}\left(
\Phi _{-t}(p_{1})-\Phi _{-t}(p_{2})\right) \geq 0.\]
\item The map $t\rightarrow \left\Vert \pi _{y}\left( \Phi
_{-t}\left( p_{1}\right) -\Phi _{-t}\left( p_{2}\right) \right) \right\Vert $
is strictly increasing. 
\end{itemize}
We will show that this implies that for some $t>0$
either $\Phi _{-t}\left( p_{1}\right) $ or $\Phi _{-t}\left( p_{2}\right) $
must exit $S^{u}$.

If this was not the case, by compactness of $\overline{S^{u}}$ we would
have a sequence $t_{n}\rightarrow \infty $ and a point $p_{1}^{\ast
}=\lim_{n\rightarrow \infty }\Phi _{-t_{n}}\left( p_{1}\right) $. From
compactness, we can choose a subsequence $t_{n_{k}}$ such that $p_{2}^{\ast
}=\lim_{k\rightarrow \infty }\Phi _{-t_{n_{k}}}\left( p_{2}\right) $. Since by the $Q_s$-backward cone condition $
\Phi _{-t_{n_{k}}}\left( p_{1}\right) \in Q_{s}^{+}(\Phi _{-t_{n_{k}}}\left(
p_{2}\right) )$, by passing to the limit, we obtain $p_{1}^{\ast }\in
Q_{s}^{+}(p_{2}^{\ast })$. Also%
\begin{equation*}
\left\Vert \pi _{y}\left( \Phi _{-t_{n_{k}}}\left( p_{1}\right) -\Phi
_{-t_{n_{k}}}\left( p_{2}\right) \right) \right\Vert \overset{k\rightarrow
\infty }{\rightarrow }\left\Vert \pi _{y}\left( p_{1}^{\ast }-p_{2}^{\ast
}\right) \right\Vert ,
\end{equation*}%
and by the strict monotonicity of $t \mapsto \left\Vert \pi _{y}\left( \Phi _{-t}\left( p_{1}\right) -\Phi
_{-t}\left( p_{2}\right) \right) \right\Vert$
 \begin{equation}
 \left\Vert \pi _{y}\left( \Phi _{-t}\left( p_{1}\right) -\Phi
_{-t}\left( p_{2}\right) \right) \right\Vert <\left\Vert \pi _{y}\left(
p_{1}^{\ast }-p_{2}^{\ast }\right) \right\Vert, \quad \mbox{for all $t>0$}.
\end{equation} 
Observe that from the above condition it follows that $\pi_y p_{1}^{\ast }\neq \pi_y p_{2}^{\ast }$. Hence at least one of the points
$p_1^\ast$, $p_2^\ast$ must not belong to $\Lambda$, so that (\ref{eq:change-in-y-2}) holds and the  function $t \mapsto \left\Vert \pi _{y}\left( \Phi _{-t}\left( p_{1}^\ast\right) -\Phi
_{-t}\left( p_{2}^\ast\right) \right) \right\Vert$ is increasing.

Taking any $s>0$. We have,
\[
\left\Vert \pi _{y}\left( p_{1}^{\ast }-p_{2}^{\ast }\right) \right\Vert
>\left\Vert \pi _{y}\left( \Phi _{-t_{n_{k}}-s}\left( p_{1}\right) -\Phi
_{-t_{n_{k}}-s}\left( p_{2}\right) \right) \right\Vert.
\]
Taking the limit $k\to\infty$, we obtain a contradiction:
\[
\left\Vert \pi _{y}\left( p_{1}^{\ast }-p_{2}^{\ast }\right) \right\Vert
>\left\Vert \pi _{y}\left( \Phi
_{-s}\left( p_{1}^{\ast }\right) -\Phi _{-s}\left( p_{2}^{\ast }\right)
\right) \right\Vert >\left\Vert \pi _{y}\left( p_{1}^{\ast }-p_{2}^{\ast
}\right) \right\Vert .
\]
This finishes our proof.
\end{proof}


\section{Description of the PCR3BP at infinity}

\label{sec:eqinfty}

In this section we introduce the Planar Circular Restricted 3-Body Problem
and present several sets of coordinates that will be useful in our construction.

\subsection{Equations at infinity\label{subsec:eqatinfty}}

Let $(r,\alpha)$ be polar coordinates in the plane and let $(y,G)$ be their
symplectic conjugate momenta, i.e. $y=\dot{r}$ is the momentum in the radial
direction and $G=r^{2}\dot{\alpha}$ is the angular momentum. Then, the
Hamiltonian for the planar circular restricted three body problem (PCR3BP)
in the inertial frame, where the primaries are rotating, takes the form
\begin{equation}
\mathcal{H}(r,\alpha,y,G,t)=\frac{1}{2}\left( \frac{G^{2}}{r^{2}}+y^{2}\right)
-U(r,\alpha-t),  \label{eq:Ham}
\end{equation}
where
\begin{equation*}
U(r,\phi)=\frac{1-\mu}{\sqrt{r^{2}+2\mu r\cos\phi+\mu^{2}}}+\frac{\mu}{\sqrt{%
r^{2}-2\left( 1-\mu\right) r\cos\phi+\left( 1-\mu\right) ^{2}}}
\end{equation*}
is the Newtonian potential describing the interaction of the massless body with the
primaries, which move on circular orbits. 
In the rotating coordinate frame $\phi=\alpha-t$, the Hamiltonian \eqref{eq:Ham} becomes the Hamiltonian $H$ in \eqref{def:HamCartesian}
in polar coordinates, that is
\begin{equation}\label{eq:Jacobi-integral}
H(r,\phi,y,G)=\frac{1}{2}\left( \frac{G^{2}}{r^{2}}+y^{2}\right)
-U(r,\phi)-G.
\end{equation}

Since we want to study the invariant manifolds of infinity, we consider the
McGehee coordinates $(x,y,\phi,G)$ where
\begin{equation}\label{def:McGehee}
r=\frac{2}{x^{2}},\qquad x>0.
\end{equation}
Taking%
\begin{equation}
\mathcal{U}\left( x,\phi\right) =U\left( 2x^{-2},\phi\right) ,
\label{eq:put-mathc-U}
\end{equation}
we obtain the following ODE%
\begin{align}
\dot{x} & =-\frac{1}{4}x^{3}y,  \label{eq:xy-ode} \\
\dot{y} & =\frac{1}{8}x^{6}G^{2}-\frac{x^{3}}{4}\frac{\partial\mathcal{U}}{%
\partial x},  \notag \\
\dot{\phi} & =\frac{1}{4}x^{4}G-1,  \notag \\
\dot{G} & =\frac{\partial\mathcal{U}}{\partial\phi},  \notag
\end{align}
which is reversible with respect to the involution
\begin{equation}
\mathcal{S}\left( x,y,\phi,G\right) =\left( x,-y,-\phi,G\right).
\label{eq:symmetry}
\end{equation}
That is,  the flow $\Phi_{t}$ induced by (\ref{eq:xy-ode}) satisfies
\begin{equation}
\Phi_{t}\circ\mathcal{S}=\mathcal{S}\circ\Phi_{-t}.  \label{eq:symmetry-flow}
\end{equation}

Let us use the following notation $\mathcal{O}_{k}(x)=\mathcal{O}(|x|^{k})$.
Since $
\left( 1+x\right) ^{-1/2}=1-\frac{1}{2}x+\frac{3}{8}x^{2}+\mathcal{O}_{3}(x)$, 
one has
\comment{
\begin{eqnarray*}
U(r,\phi ) &=&\frac{1-\mu }{r}\left( 1+2\frac{\mu }{r}\cos \phi +\left(
\frac{\mu }{r}\right) ^{2}\right) ^{-1/2}+\frac{\mu }{r}\left( 1-2\frac{\left( 1-\mu \right) }{r}\cos \phi +\left( \frac{1-\mu }{r}\right)
^{2}\right) ^{-1/2} \\
&=&\frac{1}{r}-\frac{1-\mu }{r}\frac{1}{2}\left( 2\frac{\mu }{r}\cos \phi
+\left( \frac{\mu }{r}\right) ^{2}\right) -\frac{\mu }{r}\frac{1}{2}\left( -2\frac{\left( 1-\mu \right) }{r}\cos \phi +\left( \frac{1-\mu }{r}\right)
^{2}\right) \\
&&+\frac{1-\mu }{r}\frac{3}{4}\left( 2\frac{\mu }{r}\cos \phi +\left( \frac{\mu }{r}\right) ^{2}\right) ^{2}+\frac{\mu }{r}\frac{3}{4}\left( -2\frac{\left( 1-\mu \right) }{r}\cos \phi +\left( \frac{1-\mu }{r}\right)
^{2}\right) ^{2} \\
&&+\mathcal{O}_1\left(\frac{\mu(1-\mu)}{r^4}\right) \\
&=&\frac{1}{r}-\frac{1-\mu }{r}\frac{1}{2}\left( \frac{\mu }{r}\right) ^{2}-\frac{\mu }{r}\frac{1}{2}\left( \frac{1-\mu }{r}\right) ^{2} \\
&&+\frac{1-\mu }{r}\frac{3}{4}\left( 2\frac{\mu }{r}\cos \phi \right) ^{2}+\frac{\mu }{r}\frac{3}{4}\left( 2\frac{\left( 1-\mu \right) }{r}\cos \phi
\right) ^{2}+\mathcal{O}_{1}\left( \frac{\mu(1-\mu) }{r^{4}}\right) \\
&=&\frac{1}{r}-\frac{\left( 1-\mu \right) \mu }{2}\frac{1}{r^{3}} \\
&&+\left( 1-\mu \right) \mu 3\cos ^{2}\phi \frac{1}{r^{3}}+\mathcal{O}_{1}\left( \frac{\mu(1-\mu) }{r^{4}}\right) \\
&=&\frac{1}{r}-\frac{\mu (1-\mu )}{2}\left( 1-6\cos ^{2}(\phi )\right) \frac{1}{r^{3}}+\mathcal{O}_{1}\left( \frac{\mu(1-\mu) }{r^{4}}\right)
\end{eqnarray*}
} 
\begin{align*}
U(r,\phi) =\,&\frac{1-\mu}{r}\left( 1+2\frac{\mu}{r}\cos\phi+\left( \frac{\mu%
}{r}\right) ^{2}\right) ^{-1/2} \\
& +\frac{\mu}{r}\left( 1-2\frac{\left( 1-\mu\right) }{r}\cos\phi+\left(
\frac{1-\mu}{r}\right) ^{2}\right) ^{-1/2} \\
 =\,&\frac{1}{r}\left( 1-\frac{\mu(1-\mu)}{2}\left( 1-3\cos^{2}\phi\right)
\frac{1}{r^{2}}+\mathcal{O}_{1}\left( \frac{\mu(1-\mu)}{r^{3}}\right)
\right).
\end{align*}
Hence
\begin{equation}
\mathcal{U}(x,\phi)=U\left( 2/x^{2},\phi\right) =\frac{x^{2}}{2}\left( 1-%
\frac{\mu(1-\mu)}{2}\left( 1-3\cos^{2}\phi\right) \frac{x^{4}}{4}+\mathcal{%
O}\left( \mu x^{6}\right) \right) .  \label{eq:U-x2}
\end{equation}
Therefore,
\[
\frac{\partial\mathcal{U}}{\partial x}  =x+\mathcal{O}_{5}(x), \qquad
\frac{\partial\mathcal{U}}{\partial\phi}  =\beta(\phi)x^{6}+\mathcal{O}_{8}(x)
\]
where
\[
 \beta(\phi)=\frac{3\mu(1-\mu)}{8}\cos\phi\sin\phi.
\]
This means that in (\ref{eq:xy-ode}%
) in the equation for $\dot{y}$ the dominant term near $x=y=0$ will be $-%
\frac{x^{3}}{4}\frac{\partial\mathcal{U}}{\partial x}=-\frac{x^{4}}{4}+%
\mathcal{O}_{6}(x)$.

The manifold at infinity
\begin{equation*}
\Lambda=\{(x,y,\phi,G)\in\mathbb{R}\times\mathbb{R}\times\mathbb{T}\times%
\mathbb{R}\ |\ x=y=0\},
\end{equation*}
is invariant and is foliated by periodic orbits
\begin{equation}\label{def:LambdaInfty}
\Lambda_{I}=\Lambda\cap\left\{ G=-I\right\} .
\end{equation}
Observe that $H\equiv 0$ on $\Lambda$. 

\begin{remark}
In the definition of $\Lambda_{I}$ we introduce the minus on the right hand
side because then the periodic orbit $\Lambda_{I}$ will belong to the energy level  $H=I$ (see (\ref{eq:Jacobi-integral})).
\end{remark}

\subsection{Invariant sector}

The system (\ref{eq:xy-ode}) can be written as
\begin{align}
\dot{x} & =-\frac{1}{4}x^{3}y  \label{eq:dot-x1} \\
\dot{y} & =-\frac{1}{4}x^{4}+x^{6}\mathcal{O}_{1}  \notag \\
\dot{\phi} & =\frac{1}{4}x^{4}G-1  \notag \\
\dot{G} & =\beta(\phi)x^{6}+x^{8}\mathcal{O}_{1},  \notag
\end{align}
where $\beta$ and all $O_{1}$ functions above are $2\pi$-periodic in $\phi$.

To straighten the lowest order terms, we make the change
\begin{equation}\label{def:linearchange}
\begin{aligned}
q & =\frac{1}{2}(x-y), \\
p & =\frac{1}{2}(x+y), \\
\theta & =\phi+Gy.
\end{aligned}
\end{equation}
The coordinate change $\theta=\phi+Gy$, might appear artificial, as we are
adding $Gy$ to an angle $\theta$, however the system is still $2\pi$-periodic in $%
\theta$, which means that we can treat this variable as an angle.

Note that we are only interested in the region $p+q \geq 0$ (see \eqref{def:McGehee}). We do assume this fact throughout this section without mentioning it.

Then we have the new system, which has the form
\begin{align}
\dot{q} & =\frac{1}{4}(q+p)^{3}\left( q+(q+p)^{3}\mathcal{O}_{0}\right) ,
\label{eq:ode-diag} \\
\dot{p} & =-\frac{1}{4}(q+p)^{3}\left( p+(q+p)^{3}\mathcal{O}_{0}\right) ,
\notag \\
\dot{G} & =(p+q)^{6}\mathcal{O}_{0},  \notag \\
\dot{\theta} & =(p+q)^{6}\mathcal{O}_{0}-1.  \notag
\end{align}

Let $F$ denote the vector field on the right hand side of (\ref{eq:ode-diag}%
). We note that the derivative of $F$ is of the form%
\begin{equation}
DF=(p+q)^{3}\left( \frac{1}{4}\left(
\begin{array}{cccc}
1+\frac{3q}{q+p} & \frac{3q}{q+p} & 0 & 0 \\
-\frac{3p}{q+p} & -1-\frac{3p}{q+p} & 0 & 0 \\
0 & 0 & 0 & 0 \\
0 & 0 & 0 & 0%
\end{array}
\right) +\mathcal{O}_{1}\right) .  \label{eq:deriv-not-diag}
\end{equation}
This means that we can factor out the term
\begin{equation*}
h(q,p,\theta,G)=\left( p+q\right) ^{3},
\end{equation*}
in front of the derivative of the vector field.

On the level set of  the Hamiltonian $H=I$, we have $G=G(x,y,\theta,I)$.
From (\ref{eq:Jacobi-integral}) we know that%
\begin{equation}
I=\frac{x^{4}G^{2}}{8}+\frac{1}{2}y^{2}-\mathcal{U}(x,\phi)-G.
\label{eq:Jacobi-G}
\end{equation}
Since the equation (\ref{eq:Jacobi-G}) is quadratic in $G$, it can be explicitly solved for $G(x,y,\theta,I)$. The
analytical formula has a singularity in the denominator for $x \to 0$, hence for the rigorous numerical computation of $G$ and its
derivatives  it is convenient to use a different approach.
Let us recall that by (\ref{eq:U-x2}) $\mathcal{U}(x,\phi)=x^{2}\mathcal{O}%
_{0}$, hence for $x$ close to zero $G$ can be solved for as $G\approx-I$.

Let $I$ be fixed and let $\varrho_{I}:\mathbb{R}^{3}\times\mathbb{S}%
^{1}\rightarrow\mathbb{R}$ be defined as%
\begin{equation*}
\varrho_{I}(q,p,G,\theta):=\frac{\left( p+q\right) ^{4}G^{2}}{8}+\frac{1}{2}%
\left( p-q\right) ^{2}-\mathcal{U}(p+q,\theta-G\left( p-q\right) )-G-I.
\end{equation*}
We consider $G_{I}=G_{I}(q,p,\theta)$ to be the solution of $\varrho
_{I}(q,p,G_{I},\theta)=0$, which satisfies $G_{I}(0,0,\theta)=-I$. The
lemma below is a tool which we use in our computer assisted proof to establish
that such $G_{I}$ is well defined and to validate explicit bounds for its
values.

\begin{lemma}
\label{lem:implicit-G-Newton}Let $\mathbf{G}$ be a closed interval and let $%
G_{0}\in\mathrm{int}\mathbf{G}$. Let $q,p,\theta$ be fixed. If%
\begin{equation*}
G_{0}-\left( \frac{\partial\varrho_{I}}{\partial G}(q,p,\mathbf{G}%
,\theta)\right) ^{-1}\varrho_{I}(q,p,G_{0},\theta)\subset\mathbf{G,}
\end{equation*}
then there esists a $G_{I}\left( q,p,\theta\right) \in\mathbf{G}$ such that $%
\varrho_{I}(q,p,G_{I}\left( q,p,\theta\right) ,\theta)=0.$
\end{lemma}

\begin{proof}
The result follows from the interval Newton method \cite[Theorem 13.2]{MR2652784}. 
\end{proof}


\begin{corollary}
\label{cor:implicit-G-Newton} Lemma \ref{lem:implicit-G-Newton} works under
the implicit assumption that $0\notin\frac{\partial\varrho_{I}}{\partial G}%
(q,p,\mathbf{G},\theta)$, so from the implicit function theorem we also
obtain that $G_{I}\left( q,p,\theta\right) \in C^{1}$ and the bounds on its
derivatives as%
\begin{align*}
\frac{\partial G_{I}}{\partial q}\left( q,p,\theta\right) & \in-\left( \frac{%
\partial\varrho_{I}}{\partial G}\left( q,p,\mathbf{G},\theta\right) \right)
^{-1}\frac{\partial\varrho_{I}}{\partial q}\left( q,p,\mathbf{G}%
,\theta\right) , \\
\frac{\partial G_{I}}{\partial p}\left( q,p,\theta\right) & \in-\left( \frac{%
\partial\varrho_{I}}{\partial G}\left( q,p,\mathbf{G},\theta\right) \right)
^{-1}\frac{\partial\varrho_{I}}{\partial p}\left( q,p,\mathbf{G}%
,\theta\right) , \\
\frac{\partial G_{I}}{\partial\theta}\left( q,p,\theta\right) & \in-\left(
\frac{\partial\varrho_{I}}{\partial G}\left( q,p,\mathbf{G},\theta\right)
\right) ^{-1}\frac{\partial\varrho_{I}}{\partial\theta }\left( q,p,\mathbf{G}%
,\theta\right) .
\end{align*}
\end{corollary}

We see that in a neighborhood of $\Lambda $ the coordinate $q$ is ``expanding'', $p$ is
``contracting'', and $\theta ,G$ are center coordinates. This means that we can
expect the set%
\begin{equation}
S^{u}=S_{I,L,R}^{u}:=\left\{ \left( q,p,G_{I}\left( q,p,\theta \right)
,\theta \right) :q\in \left( 0,R\right) ,\left\vert p\right\vert <Lq,\theta
\in \mathbb{S}^{1}\right\} ,  \label{eq:Scu-3bp}
\end{equation}%
to be an invariant sector (see Definition \ref{def:Scu-sector}) for suitably
chosen $R,L>0$.

\begin{remark}
We fix the energy level $H=I$, and for such fixed value our
system becomes three dimensional. We therefore treat the sector $%
S_{I,L,R}^{u}$ as a subset of a three dimensional space, with coordinates $%
q,p,\theta $.
\end{remark}

The following lemma gives conditions to prove the existence of an unstable
sector  (according to Definition \ref{def:Scu-sector}).

\begin{lemma}
\label{th:invariant-sector}Let $F=\left( F_{1},F_{2},F_{3},F_{4}\right) :%
\mathbb{R}^{3}\mathbb{\times \mathbb{S}}^{1}\mathbb{\rightarrow R}^{4}$
stand for the vector field on the right hand side of (\ref{eq:ode-diag}). If
for every
\begin{align*}
\left\langle \left( F_{1},F_{2}\right) (z),\left( L,-1\right)
\right\rangle & >0\qquad \text{for all }z\in \left\{ p=Lq\right\} \cap
\overline{S_{I,L,R}^{u}}, \\
\left\langle \left( F_{1},F_{2}\right) (z),\left( L,1\right)
\right\rangle & >0\qquad \text{for all }z\in \left\{ p=-Lq\right\} \cap
\overline{S_{I,L,R}^{u}}, \\
F_{1}(z)&>0 \qquad \text{for all }z\in \left\{ q=R\right\} \cap
\overline{S_{I,L,R}^{u}},
\end{align*}
then $S_{I,L,R}^{u}$ is an unstable sector.
\end{lemma}

\begin{proof}
The assumptions imply that the flow can not exit $S_{I,L,R}^{u}$ through $%
\left\vert p\right\vert =Lq$ and must exit through $\{q=R\}$.
\end{proof}

We define a sector%
\begin{equation}
S^{s}=S_{I,L,R}^{s}:=\left\{ \left( q,p,G_{I}\left( q,p,\theta \right)
,\theta \right) :p\in \left( 0,R\right) ,\left\vert q\right\vert <Lp,\theta
\in \mathbb{S}^{1}\right\} .  \label{eq:Scs-3bp}
\end{equation}

The proposition below shows that points exiting a neighbourhood of $\Lambda _{I}$ (within $p+q>0$)
must do so through the unstable sector. This is a technical result,
which will be useful in our construction for the proof of oscillatory
motions in section \ref{sec:oscillatory-proof}.

\begin{figure}
\begin{center}
\includegraphics[height=4cm]{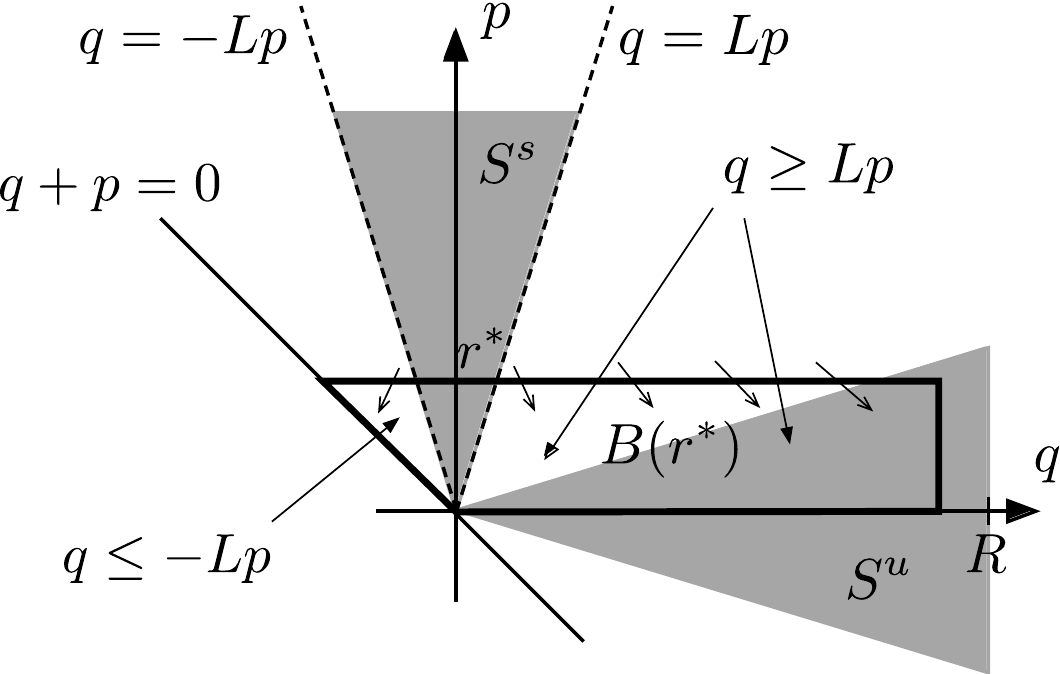}
\end{center}
\caption{The set $B(r^*)$ from Proposition \ref{th:orbits-enter-Scu}.\label{fig:Br}}
\end{figure}

\begin{proposition}
\label{th:orbits-enter-Scu}
Assume that $S^{u}$ and $S^{s}$ defined in
(\ref{eq:Scu-3bp}) and (\ref{eq:Scs-3bp}) are unstable invariant and stable
sectors, respectively. Let
\[
B(r)=\{(q,p,G_{I}\left(  q,p,\theta\right)  ,\theta)\,|\,-r<q<\sqrt
{r},0<p<r,\theta\in\mathbb{S}^{1},q+p>0\}.
\]
There exists $r^{\ast}>0$, such that:

\begin{enumerate}
\item For every
\[
z_{0}\in\left(  B(r^{\ast})\cap\left\{  q>0\right\}  \right)  \setminus\left(
S^{u}\cup S^{s}\right)
\]
there exists $T=T(z_{0})>0$ such that $z(T)\in S^{u}$, where $z(t)$ is a
solution of (\ref{eq:ode-diag}) with initial condition $z(0)=z_{0}$.

\item For every
\[
z_{0}\in\left(  B(r^{\ast})\cap\left\{  q<0\right\}  \right)  \setminus S^{s}%
\]
and the solution $z(t)$ starting from $z_{0}$ we have%
\[
\lim_{t\rightarrow\infty}\pi_{q,p}z\left(  t\right)  =\left(  q^{\ast}%
,p^{\ast}\right)  ,
\]
where $p^{\ast}\in\left(  -\pi_{q}z_{0},r^{\ast}\right)  $ and $q^{\ast
}+p^{\ast}=0$.  
\end{enumerate}
\end{proposition}

\begin{proof}
Below we will choose $r^{\ast}$ sufficiently small so that $r^{\ast}<1,$
$\sqrt{r^{\ast}}<R$ and (see Figure \ref{fig:Br})
\[
(q=\sqrt{r^{\ast}},p=r^{\ast})\in\pi_{q,p}S^{u}=\pi_{q,p}S_{I,L,R}^{u}.
\]
Observe that we have freedom to decrease $r^*$ and still the above condition will be satisfied.

Let $M$ be such that for $\mathcal{O}_{0}$ in (\ref{eq:ode-diag}) we have a
bound $|\mathcal{O}_{0}|<M$ for some macroscopic neighborhood of $(0,0)$. We
choose $r^{\ast}$ small enough so that $B(r^{\ast})$ lies in this neighbourhood.

We start by showing that if $r^{\ast}$ is also small enough so that $r^{\ast
}<\frac{1}{64M^{2}}$ then for every point from $\overline{B\left(  r^{\ast
}\right)  }\cap\left\{  p=r^{\ast}\right\}  $ and any $t$ holds
\begin{equation}
\frac{dp}{dt}<0.\label{eq:p-prime-negative}%
\end{equation}
In view of (\ref{eq:ode-diag}) and the positivity of $p+q$, in order to show
(\ref{eq:p-prime-negative}) it is enough to check that $p+(p+q)^{3}%
\mathcal{O}_{0}>0$. Since $\left\vert q\right\vert \leq\sqrt{r^{\ast}}$ and
$r^{\ast}<\frac{1}{64M^{2}}$ we obtain
\[
p+(p+q)^{3}\mathcal{O}_{0}\geq r^{\ast}-M(r^{\ast}+\sqrt{r^{\ast}}%
)^{3}>r^{\ast}-M(2\sqrt{r^{\ast}})^{3}>0.
\]

We can see from (\ref{eq:p-prime-negative}) that a trajectory can not exit the
set $B(r^{\ast})\setminus\left(  S^{u}\cup S^{d}\right)  $ through $\left\{
p=r^{\ast}\right\}  $. (See Figure \ref{fig:Br})

We now show that if we choose $r^{\ast}<M^{-1}\left(  1+\frac{1}{L}\right)
^{-3}$ then for every point from  $(B(r^\ast) \cap \{q>0\}) \setminus (S^u \cup S^s) \subset $
 $B(r^{\ast})\cap\left\{  q\geq Lp\right\}  $
and any $t$ holds
\begin{equation}
\frac{dq}{dt}>0.\label{eq:q-prime-positive}%
\end{equation}
For this it is enough to show that $q+(p+q)^{3}\mathcal{O}_{0}>0$. Since
$q\geq Lp$ and $q\leq\sqrt{r^{\ast}}<\sqrt{M^{-1}\left(  1+\frac{1}{L}\right)
^{-3}}$ we see that
\[
q+(p+q)^{3}\mathcal{O}_{0}\geq q-M(q+\frac{1}{L}q)^{3}=q\left(  1-q^{2}%
M\left(  1+\frac{1}{L}\right)  ^{3}\right)  >0.
\]

We are ready to prove the first claim. Let us fix $z_{0}\in ( B(r^{\ast
})\cap \{q>0\})\setminus\left(  S^{u}\cup S^{s}\right)  $. Let
\[
c\left(  z_{0}\right)  =\min\left\{  \frac{dq}{dt}(z): z=(q,p,G,\theta) \in B(r^\ast), \mbox{ such that} \quad q\in\left[  \pi_{q}%
z_{0},\sqrt{r^{\ast}}\right]  ,q\geq Lp\right\}  >0.
\]
For a trajectory $z\left(  t\right)  $ starting from $z_{0}$, for $t>0$ we
will therefore have%
\[
\pi_{q}z\left(  t\right)  >\pi_{q}z_{0}+c\left(  z_{0}\right)  t,
\]
as long as $z\left(  t\right)  \in B\left(  r^{\ast}\right)  $. Since $z(t)$
can not pass through $\left\{  p=r^{\ast}\right\}  $ we will have $z\left(
T(z_{0})\right)  \in S^{u}$ for some $T(z_{0})<\sqrt{r^{\ast}}/c\left(
z_{0}\right)  $.

To prove the second claim, observe that using mirror arguments to the proof of
(\ref{eq:q-prime-positive}) we obtain that for every point from $B\left(
r^{\ast}\right) \cap \{q <0\} \setminus S^{s}=B\left(  r^{\ast}\right)  \cap\left\{
q\leq-Lp\right\}  $   we have%
\[
\frac{dq}{dt}<0.
\]

For every point $z_{0}\in\left(  B(r^{\ast})\cap\left\{  q<0\right\}  \right)
\setminus S^{s}$ the trajectory $z\left(  t\right)  $ which starts from
$z_{0}$ can not exit $B\left(  r^{\ast}\right)  $ through $\left\{  p=r^{\ast
}\right\}  $. Moreover it cannot exit through the line $x=p+q=0$ as it consist of the fixed points. Since $q(t)$ is strictly increasing and bounded, it converges to $q^\ast \in [-r^\ast,\pi_q z_0]$ and it is easy to see that then the trajectory must converge to the line $x=p+q=0$.  We have obtained
\[
\lim_{t\rightarrow\infty}\pi_{q,p}z\left(  t\right)  =\left(  q^{\ast}%
,p^{\ast}\right)
\]
for some $\left(  q^{\ast},p^{\ast}\right)  \in\left\{  q+p=0\right\}  $.
Therefore $p^{\ast}=-q^{\ast}\in\left(  -\pi_{q}z_{0},r^{\ast}\right) $. This finishes our proof.
\end{proof}

\subsection{Bringing the factorised derivative at infinity to diagonal form}

\label{subsec:diag}

In the dominant part of matrix $DF$ given by (\ref{eq:deriv-not-diag}) in
the limit of $\frac{|p|}{q}\rightarrow0$ (i.e. very tight sector) we have $%
\frac{3q}{p+q}\rightarrow3$ and $\frac{3p}{p+q}\rightarrow0$. Therefore,  we obtain
there is an off-diagonal ``non-small'' term corresponding to the entry $\frac{\partial F_{q}}{%
\partial p}$. The presence of such term is undesirable, because it makes the
verification of cone conditions for the parabolic invariant manifold harder.

In this section we discuss a simple change of coordinates of
the form $u=q+\left( 1-b\right) p$, for some $b\in\mathbb{R}$, which leaves the
remaining coordinates $p,\theta,G$ unchanged. This change will take the derivative
of the factorised vector field (\ref{eq:deriv-not-diag}) at $\Lambda$ to
diagonal form.

By taking
\begin{equation}\label{def:qtou}
u=q+\left( 1-b\right) p
\end{equation}
we obtain the ODE%
\begin{equation}\label{eq:final-ode-3bp}
\begin{split}
\dot{u} & =\frac{1}{4}(u+bp)^{3}\left( u-2\left( 1-b\right) p+(u+bp)^{3}%
\mathcal{O}_{0}\right), \\
\dot{p} & =-\frac{1}{4}(u+bp)^{3}\left( p+(u+bp)^{3}\mathcal{O}_{0}\right) ,
\\
\dot{G} & =\left( u+bp\right) ^{6}\mathcal{O}_{0},  \\
\dot{\theta} & =\left( u+bp\right) ^{6}\mathcal{O}_{0}-1.
\end{split}
\end{equation}
Let us denote the vector field on the right hand side in (\ref%
{eq:final-ode-3bp}) by $F=\left( F_{u},F_{p},F_{G},F_{\theta}\right) $.

It is easy to see that
\begin{equation*}
\frac{\partial F_{u}}{\partial p}=\frac{1}{4}(u+bp)^{3}\left( W+(u+bp)^{2}%
\mathcal{O}_{0}\right) ,
\end{equation*}
where
\begin{align*}
W & =\frac{3b}{u+bp}\left( u-2(1-b)p\right) -2(1-b) \\
& =3b\frac{u}{u+bp}-6b(1-b)\frac{p}{u+bp}-2(1-b) \\
& =3b\left( \frac{u}{u+bp}-1\right) -6b(1-b)\frac{p}{u+bp}+(5b-2) \\
& =3b\frac{-bp}{u+bp}-6b(1-b)\frac{p}{u+bp}+(5b-2) \\
& =3b(b-2)\frac{p}{u+bp}+(5b-2).
\end{align*}
So if we take $b=2/5$, then we get rid of the first
term in the bracket and can factor out $(u+bp)^{3}$ in front of the
derivative of the vector field, obtaining%
\begin{equation} \label{eq:factorised-derivative-3bp}
\begin{aligned}
& DF\left( u,p,G,\theta\right) = \\
&\frac{ (u+bp)^{3}}{4}\left(\left(
\begin{array}{cccc}
4+3\left( b-2\right) \frac{p}{u+bp} & 3b(b-2)\frac{p}{u+bp} & 0 & 0 \\
-3\frac{p}{u+bp} & -1-3\frac{p}{u+bp} & 0 & 0 \\
0 & 0 & 0 & 0 \\
0 & 0 & 0 & 0%
\end{array}
\right) +\mathcal{O}_{2}\left( u+bp\right) \right) .
\end{aligned}
\end{equation}
with $b=2/5$.

Since for points from a sector $|p|\leq L_{1}q$
\begin{equation*}
\left\vert \frac{p}{u+bp}\right\vert =\left\vert \frac{p}{q+p}\right\vert
\leq\frac{L_{1}}{1-L_{1}},
\end{equation*}
we see that for small $L_{1}$ the factorised derivative of the vector field
will be close $\mathrm{diag}\left( 4,-1,0,0\right) ,$ when computed at
points from the sector.


\section{Bounds on the unstable manifold at infinity in the PCR3BP}\label{sec:InvManPCR3B}

Let us fix $I$ and consider the manifold $\{H=I\}$ (see (\ref{eq:Jacobi-integral})) and the periodic orbit $\Lambda_{I}$ introduced in  \eqref{def:LambdaInfty}. Analogously to \eqref{eq:Scu-3bp},  we consider a
sector $S_{I,L,R}^{u}$ in the coordinates $(u,p,\theta)$.

We will work in a
setting where, by Lemma \ref{lem:implicit-G-Newton}, for $(u,p,G,\theta )\in \{J=I\}$, we have $G=G_{I}(u,p,\theta)$.
This means that $(u,p,\theta)$ uniquely define the point $(u,p,G_{I}(u,p,\theta),\theta)$.

Our aim will be to apply Theorem \ref{th:Wcu-bound} to establish the existence of a
center-horizontal disc $w^{u}:\overline{B_{u}}\times \mathbb{S}%
^{1}\rightarrow S_{I,L,R}^{u}$, for $\overline{B_{u}}=\left[ 0,R\right]
\subset \mathbb{R}$, such that the unstable manifold $W_{\Lambda
_{I}}^{u}$ of $\Lambda _{I}$  for the flow of (\ref{eq:final-ode-3bp})
projected on to the $u,p,\theta $ coordinates is a graph of $w^{u}$%
\begin{equation*}
W_{\Lambda _{I}}^{u}=\mathrm{graph}\left( w^{u}\right) .
\end{equation*}

Recalling the change of coordinates \eqref{def:qtou} and that $b=2/5$, we define
\begin{equation*}
\mathrm{G}_{I}\left( u,p,\theta\right) :=G_{I}\left( u-\left( 1-b\right)
p,p,\theta\right) .
\end{equation*}

From now on, for points belonging to the sector $S_{I,L,R}^{u}$, we denote by $\widetilde\Phi
_{t}\left( u,p,\theta \right) $ the projection onto the $(u,p,\theta)$ coordinates of the
flow assoacaiated to \eqref{eq:final-ode-3bp} with initial condition at the point $\left(
u,p,{\mathrm{G}}_{I}\left( u,p,\theta \right) ,\theta \right) $.

Let $P$ and $M$ be the following matrices
\begin{equation*}
P:=\left(
\begin{array}{cccc}
1 & 0 & 0 & 0 \\
0 & 1 & 0 & 0 \\
0 & 0 & 0 & 1%
\end{array}%
\right) ,\qquad M:=\left(
\begin{array}{ccc}
1 & 0 & 0 \\
0 & 1 & 0 \\
\frac{\partial {\mathrm{G}}_{I}}{\partial u} & \frac{\partial {\mathrm{G}}%
_{I}}{\partial p} & \frac{\partial {\mathrm{G}}_{I}}{\partial \theta } \\
0 & 0 & 1%
\end{array}%
\right).
\end{equation*}%
On the invariant surface $J=I$, the coordinates are $(u,p,\theta )$ and the vector field is $\tilde{F}(u,p,\theta)=P F(u,p,G_I(u,p,\theta),\theta)$. Hence
the derivative of the vector field is
\begin{equation}
D\tilde{F}(u,p,\theta ):=P\,DF(u,p,\mathrm{G}_{I}\left( u,p,\theta \right)
,\theta )\,M.  \label{eq:derivative-in-reduced-coord}
\end{equation}%
From $D\tilde{F}(u,p,\theta )$ we can factorise the term $h\left( u,p,\theta
\right) =(u+bp)^{3}$. Let us use the notation $\mathbf{G}$ for a $3\times 3$
interval matrix, which is an interval enclosure of the factorised $D\tilde{F}
$. In other words, for every $\left( u,p,\theta \right) \in \pi _{u,p,\theta
}S_{I,L,R}^{u}$ let
\begin{equation}
D\tilde{F}\left( u,p,\theta \right) \in (u+bp)^{3}\mathbf{G}%
=(u+bp)^{3}\left(
\begin{array}{ccc}
\mathbf{G}_{uu} & \mathbf{G}_{up} & \mathbf{G}_{u\theta } \\
\mathbf{G}_{pu} & \mathbf{G}_{pp} & \mathbf{G}_{p\theta } \\
\mathbf{G}_{\theta u} & \mathbf{G}_{\theta p} & \mathbf{G}_{\theta \theta }%
\end{array}%
\right) .  \label{eq:factorised-reduced-deriv}
\end{equation}

\begin{remark}
The $\mathbf{G}$ can be obtained by validating assumptions of Lemma \ref%
{lem:implicit-G-Newton} to obtain bounds on $\mathrm{G}_{I}\left( u,p,\theta
\right) $, bounds on the  derivatives of $\mathrm{G}_{I}$ via Corollary %
\ref{cor:implicit-G-Newton}, and using these boundes for computing an interval
enclosure of $(u+bp)^{-3}P\,DF(z)\,M$ for all $z\in S_{I,L,R}^{u}$.
\end{remark}

We now formulate our main result, which is our tool for obtaining the bounds
on the unstable manifold of $\Lambda _{I}$. First we consider the
following notation. As in \eqref{eq:cones} and \eqref{eq:xi-def}, we define cones%
\begin{eqnarray*}
Q_{cu}\left( u,p,\theta \right)  &=&L_{cu}\left\Vert \left( u,\theta \right)
\right\Vert -\left\Vert p\right\Vert , \\
Q_{s}\left( u,p,\theta \right)  &=&L_{s}\left\Vert p\right\Vert -\left\Vert
\left( u,\theta \right) \right\Vert ,
\end{eqnarray*}%
and consider constants $\xi_{cu},\mu_s,\xi_s,\mu_{cu}\in \mathbb{R}$ satisfying
\begin{align*}
\xi _{cu}& \le m_{l}\left( \left(
\begin{array}{cc}
\mathbf{G}_{uu} & \mathbf{G}_{u\theta } \\
\mathbf{G}_{\theta u} & \mathbf{G}_{\theta \theta }%
\end{array}%
\right) \right) -L_{cu}\left\Vert \left(
\begin{array}{c}
\mathbf{G}_{up} \\
\mathbf{G}_{\theta p}%
\end{array}%
\right) \right\Vert , \\
\mu _{s}& \ge l\left( \mathbf{G}_{pp}\right) +\frac{1}{L_{cu}}\left\Vert \left(
\begin{array}{cc}
\mathbf{G}_{pu} & \mathbf{G}_{p\theta }%
\end{array}%
\right) \right\Vert , \\
\xi _{s}& \le m_{l}\left( -\mathbf{G}_{pp}\right) -L_{s}\left\Vert \left(
\begin{array}{cc}
\mathbf{G}_{pu} & \mathbf{G}_{p\theta }%
\end{array}%
\right) \right\Vert , \\
\mu _{cu}& \ge l\left( -\left(
\begin{array}{cc}
\mathbf{G}_{uu} & \mathbf{G}_{u\theta } \\
\mathbf{G}_{\theta u} & \mathbf{G}_{\theta \theta }%
\end{array}%
\right) \right) +\frac{1}{L_{s}}\left\Vert \left(
\begin{array}{c}
\mathbf{G}_{up} \\
\mathbf{G}_{\theta p}%
\end{array}%
\right) \right\Vert .
\end{align*}

The choice of coordinates for $Q_{cu}$ is motivated by the fact that the
coordinates $u,\theta $ are center unstable and $p$ is a stable coordinate
for the flow $\widetilde\Phi _{t}$. The constants $\xi _{cu}$ and $\mu _{s}$ will be
used to validate $Q_{cu}$ cone conditions (see Lemma \ref%
{lem:cone-verif}). The choice of coordinates for $Q_{s}$ is due to the fact
that $p$ is the unstable and $u,\theta $ are center stable coordinates for $%
\widetilde\Phi _{-t}$. The constants $\xi _{s},\mu _{cu}$ and Lemma \ref%
{lem:outflow-verif} will be used to validate the backward $Q_{s}$ cone
conditions and backward outflowing from $S^{u}$ along $Q_{s}$.

\begin{theorem}
\label{th:manifold-in-sector-3bp}Let $I,R,L_{cu},L_{s},L\in \mathbb{R}$ be
fixed and such that $R,L_{cu}>0$, $L\in (0,1)$ and $L_{s}\in (0,\pi) $. Let $%
S_{I,L,R}^{u}$ be an unstable  sector for \eqref{eq:final-ode-3bp}  in the energy level $J=I$. Assume also that every forward trajectory starting from it must exit the sector and
\begin{align}
\xi _{cu}& >\mu _{s},  \label{eq:xi-ch-mu-ch-cond} \\
\xi _{s}& >\mu _{cu},  \label{eq:xi-v-mu-v-cond} \\
\xi _{s}& >0,  \label{eq:xi-v-cond}
\end{align}%
Then, the unstable manifold $W_{\Lambda _{I}}^{u}$ is a graph of a $%
Q_{cu}$ center horizontal disc $w^{u}:\pi _{u,\theta
}S_{I,L,R}^{u}\rightarrow S_{I,L,R}^{u}$, which satisfies $\pi _{u,\theta
}w^{u}=\mathrm{id},$ and
\begin{equation*}
\left\vert \pi _{p}\left( w^{u}\left( u_{1},\theta _{1}\right)
-w^{u}\left( u_{2},\theta _{2}\right) \right) \right\vert \leq
L_{cu}\left\Vert \left( u_{1},\theta _{1}\right) -\left( u_{2},\theta
_{2}\right) \right\Vert .
\end{equation*}
\end{theorem}

\begin{proof}
By Lemma \ref{lem:cone-verif} and (\ref{eq:xi-ch-mu-ch-cond}), the flow $\widetilde\Phi_t$  induced by $\widetilde F$ satisfies forward cone conditions for $Q_{cu}$ in $%
S^{u}.$ By Lemma \ref{lem:outflow-verif} and (\ref{eq:xi-v-mu-v-cond}--\ref%
{eq:xi-v-cond}) the flow satisfies backward cone conditions for $Q_{s}$ in $%
S^{u}$ and is backward outflowing from $S^{u}$ along $Q_{s}$. The result
follows from Theorem \ref{th:Wcu-bound}.
\end{proof}

We have used Theorem \ref{th:manifold-in-sector-3bp} to validate the
following result.

\begin{theorem}
\label{th:manifold-for-given-I}Let
\begin{equation*}
I=-1,\qquad L=4\cdot 10^{-9},\qquad R=10^{-4},\qquad L_{cu}=L_{s}=10^{-5}.
\end{equation*}%
Then $S_{I,L,R}^{u}$ is an unstable sector and the unstable
manifold $W_{\Lambda _{I}}^{u}$ is a graph of a $Q_{cu}$ center horizontal
disc in $S_{I,L,R}^{u}$.
\end{theorem}

\begin{proof}
The proof follows by computer assisted validation of the assumptions of Theorems %
\ref{th:invariant-sector} and \ref{th:manifold-in-sector-3bp}.

In our computer assisted validation, we obtain the following bound of the factorised derivative term $\mathbf{G}$ from (\ref{eq:factorised-reduced-deriv}),
\begin{align}& \mathbf{G}= \notag \\
& \left(
\begin{array}{lll}
[0.9999999, 1] & [-2.561\mbox{\tt e-08}, 1.921\mbox{\tt e-09}] & [-3.142\mbox{\tt e-13}, 4.011\mbox{\tt e-21}] \\
{[-3.004\mbox{\tt e-09}, 1.511\mbox{\tt e-07}]} & [-0.25, -0.2499999] & [-1.003\mbox{\tt e-20}, 7.854\mbox{\tt e-13}] \\
{[-9.751\mbox{\tt e-12}, 1.593\mbox{\tt e-15}]} & [-3.903\mbox{\tt e-12}, 1.050\mbox{\tt e-15}] &[-1.044\mbox{\tt e-20}, 2.531\mbox{\tt e-20}]
\end{array}\right), \notag
\end{align}
which results in
\[\begin{array}{lll}
\xi_{cu}=-5.288\mbox{\tt e-12}, & \qquad & \mu_{cu}=0.00257, \\
\mu_s=-0.2348, & & \xi_s=0.2499999.
\end{array}
\]

Note that $S^{u}_{I,L,R}\cup S^{u}_{I,L,R\,-} \subset  \{x>0,y<0\}$ (in the ``original'' $(x,y)$ coordinates, see \eqref{def:linearchange}). This,  by (\ref{eq:dot-x1}), implies that  $\dot x>0$. This means that every forward trajectory starting from a point in the set $S^{u}_{I,L,R}$ must exit the set.

The computer assisted proof takes a fraction of a second, running on a standard laptop.
\end{proof}

\subsection{Extending the unstable manifold}

The goal of this section is to extend the unstable manifold beyond $%
\overline{\pi _{u,p,\theta }S_{I,L,R}^{u}}$. We use an argument based on
the properties of the vector field, to establish that the unstable manifold of the periodic orbit $\Lambda_I$ (see \eqref{def:LambdaInfty})
streches away from  $\Lambda_I$  along the corresponding unstable manifold of the two body problem. In the case when $\mu $ is small, the unstable manifold of
the two body problem proves to be a sufficiently good approximation.

We start with the description of the system for $\mu=0$. In the case of the two
body problem the Hamiltonian is given by
\begin{equation}
H(r,\alpha,y,G)=\frac{1}{2}\left( \frac{G^{2}}{r^{2}}+y^{2}\right) -\frac {1%
}{r},   \label{eq:ham2b}
\end{equation}
(compare with (\ref{eq:Ham})) and the equations of motion are
\begin{align*}
\dot{r} & =y, \\
\dot{y} & =\frac{G^{2}}{r^{3}}-\frac{1}{r^{2}},
\\
\dot{\alpha} & =\frac{G}{r^{2}}, \\
\dot{G} & =0.
\end{align*}

The equations for $r,y$ form a closed system (with $G$ being a parameter)
and we will focus just on them. We see that
\begin{equation}\label{eq:2b-r}
\begin{aligned}
\dot{r} & =y,
\\
 \dot{y} &=\frac{G^{2}}{r^{3}}-\frac{1}{r^{2}}.
\end{aligned}
\end{equation}
Let us define an effective potential $W$ for (\ref{eq:2b-r})
\begin{equation}  \label{eq:W2b}
W(r)=\frac{G^{2}}{2 r^{2}} - \frac{1}{r}.
\end{equation}
and an effective hamiltonian for (\ref{eq:2b-r})
\begin{equation}
H(r,y)=\frac{y^{2}}{2} + W(r).   \label{eq:2b-ham}
\end{equation}

We are interested in the parabolic solution (it has $y\rightarrow0$ for
large $r$), which is a solution with $H(r,y)=0$. We have
\begin{equation}
y=\pm\sqrt{\frac{2}{r}-\frac{G^{2}}{r^{2}}}=\pm x\sqrt{1-\left( \frac{Gx}{2}%
\right) ^{2}}.   \label{eq:approx-2bp}
\end{equation}
where in the last equality we have used McGehee coordinates \eqref{def:McGehee} $r=2/x^2$.

The solution with $+$ is the stable manifold of the point at infinity
(because $\dot{r}=y>0$ hence it will grow to infinity) and with $-$ is the
unstable manifold. Since, for $r\gg 1$ one has $G\approx-I$,
we expect
\begin{equation}
\label{def:gammaI}
\gamma_{I}\left( x\right) :=-x\sqrt{1-\left( \frac{Ix}{2}\right) ^{2}}
\end{equation}
to be a good approximation of the unstable manifold in the energy level $%
H=I $, when $\mu$ is small.

\begin{figure}[ptb]
\begin{center}
\includegraphics[width=6cm]{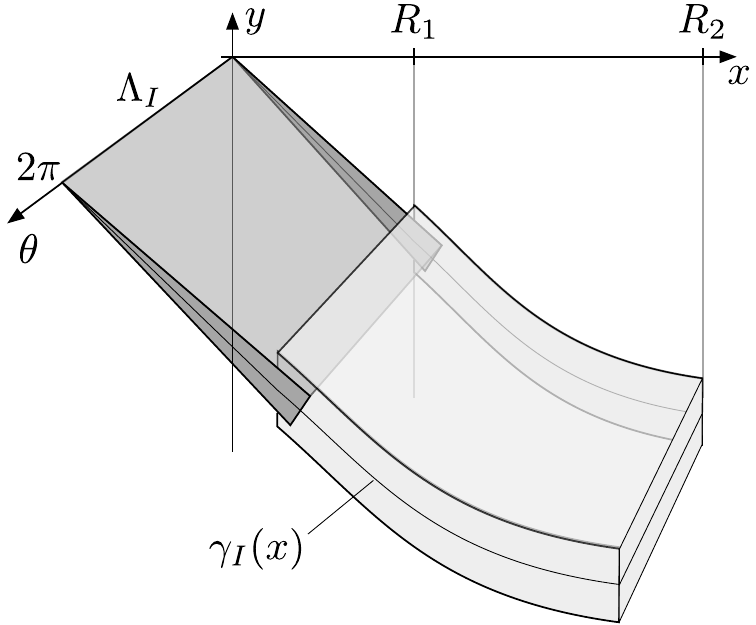}
\end{center}
\caption{The unstable sector $S^{u}_{I,L,R}$ ind dark grey, and its extension $%
S_{I,R_{1},r_{2},\protect\varepsilon}$ in light grey.}
\label{fig:extension}
\end{figure}

Let $I$ be fixed. We  now extend the sector $S_{I,L,R}^{u}$ further
away from zero. Let us consider the change of coordinates from $\left(
x,y,G,\phi \right) $ to $\left( x,\eta,G,\phi \right) $ defined as%
\begin{equation*}
y=\gamma _{I}\left( x\right) +\eta.
\end{equation*}%
Let $\mathcal{G}_{I}\left( x,\eta,\phi \right) $ be the solution for $G$ of the
quadratic equation%
\begin{equation}
I=\frac{x^{4}G^{2}}{8}+\frac{1}{2}\left( \gamma _{I}\left( x\right)
+\eta\right) ^{2}-\mathcal{U}(x,\phi )-G,  \label{eq:Jacobi-xy}
\end{equation}%
for given $x,\eta,\phi $. Let us now define the following extension of a
sector. Let $R_{1},R_{2},\varepsilon \in \mathbb{R}$ be such that $%
0<R_{1}<R_{2}$, $\varepsilon >0$, and let (see Figure \ref{fig:extension})
\begin{equation}
\label{def:SIR1R2epsilon}
S_{I,R_{1},R_{2},\varepsilon }:=\{\left( x,\eta,G,\phi \right) :x\in \left[
R_{1},R_{2}\right] ,\phi \in \mathbb{S}^{1},G=\mathcal{G}_{I}\left( x,\eta,\phi
\right) ,\eta\in \lbrack -\varepsilon ,\varepsilon ]\}.
\end{equation}

\begin{remark}
The earlier considered sector $S_{I,L,R}^{u}$ is expressed in the
coordinates $(q,p,G,\theta)$ and is connected with $\Lambda _{I}$. The set $%
S_{I,R_{1},R_{2},\varepsilon }$ is in coordinates $(x,\eta,G,\phi )$ and is
separated from $\Lambda _{I}$ since points in this set satisfy $x\geq R_{1}$.
\end{remark}

We now give a theorem which  ensures that the unstable manifold passes
through $S_{I,R_{1},R_{2},\varepsilon}$.

\begin{theorem}
\label{th:extended-sector}Let all assumptions of Theorem \ref%
{th:manifold-in-sector-3bp} be fullfilled. Let%
\begin{align*}
F_{x}\left( x,\eta,\phi \right) & :=-\frac{1}{4}x^{3}\left( \gamma _{I}\left(
x\right) +\eta\right) , \\
F_{\eta}\left( x,\eta,\phi \right) & :=\frac{1}{8}x^{6}\left( \mathcal{G}%
_{I}\left( x,\eta,\phi \right) \right) ^{2}-\frac{x^{3}}{4}\frac{\partial
\mathcal{U}}{\partial x}\left( x,\phi \right)  \\
& \quad +\frac{\partial \gamma _{I}}{\partial x}\left( x\right) \frac{1}{4}%
x^{3}\left( \gamma _{I}(x)+\eta\right) ,
\end{align*}%
Assume that there exist  $\varepsilon >0$ and $R_{1},R_{2}\in \mathbb{R}$%
, $0<R_{1}<R_{2}$, for which the following conditions hold:

\begin{enumerate}
\item The set $\overline{S_{I,L,R}^{u}}\cap \left\{ q=R\right\} $ is a
subset of $S_{I,R_{1},R_{2},\varepsilon }$; i.e. if $\left( q,p,G,\theta
\right) \in \overline{S_{I,L,R}^{u}}$ and $q=R$ then%
\begin{equation*}
\left( p+q,p-q-\gamma_I \left( p+q\right) ,G,\theta -G\left( p-q\right)
\right) \in S_{I,R_{1},R_{2},\varepsilon }.
\end{equation*}
\item We have
\begin{align*}
F_{\eta}\left( x,\varepsilon,G,\phi\right) & <0\qquad\text{for every }\left(
x,\varepsilon,G,\phi\right) \in S_{I,R_{1},R_{2},\varepsilon}, \\
F_{\eta}\left( x,-\varepsilon,G,\phi\right) & >0\qquad\text{for every }\left(
x,-\varepsilon,G,\phi\right) \in S_{I,R_{1},R_{2},\varepsilon}.
\end{align*}
\item For every $\left( x,\eta,G,\phi\right) \in S_{I,R_{1},R_{2},\varepsilon}$,%
\begin{equation*}
F_{x}\left( x,\eta,G,\phi\right) >0.
\end{equation*}
\end{enumerate}
Then, every point from $W_{\Lambda_{I}}^{u}\cap\left\{ q=R\right\} $ is
contained in $S_{I,R_{1},R_{2},\varepsilon}$ and the flow starting from such
point exits $S_{I,R_{1},R_{2},\varepsilon}$ through $\left\{ x=R_{2}\right\}
$.
\end{theorem}

\begin{proof}
Since $S_{I,L,R}^{u}$ is an unstable sector, by Theorem \ref%
{th:manifold-in-sector-3bp} it contains the unstable manifold $W^{u}.
$ The first condition therefore ensures that all points on $W^{u}\cap
\left\{ q=R\right\} $ are inside of $S_{I,R_{1},R_{2},\varepsilon }$. The $%
F_{\eta}$ is the the vector field on the coordinate $\eta$, so the second
condition ensures that trajectories starting from $W^{u}\cap \left\{
q=R\right\} $ cannot exit $S_{I,R_{1},R_{2},\varepsilon }$ through $\left\{
\eta=\varepsilon \right\} $ or $\left\{ \eta=-\varepsilon \right\} $. Since $F_{x}$
is the  vector field on the coordinate $x$, the third condition  ensures that
the coordinate $x$  increases along the flow. This means that a trajectory
starting from $W^{u}\cap \left\{ q=R\right\} $ has to exit $%
S_{I,R_{1},R_{2},\varepsilon }$ through $\left\{ x=R_{2}\right\} $.
\end{proof}

We have used Theorem \ref{th:extended-sector} to validate the following
result:

\begin{theorem}
\label{th:manifold-for-given-I-extended}Let%
\begin{equation*}
I=-1,\qquad L=4\cdot10^{-9},\qquad R=10^{-4},\qquad R_{1}=\frac{R}{2},\qquad
R_{2}=0.4,\qquad\varepsilon=2\cdot10^{-5}.
\end{equation*}
Then, every point from $W_{\Lambda_{I}}^{u}\cap\left\{ q=R\right\} $ is
contained in $S_{I,R_{1},R_{2},\varepsilon}$ and the flow starting from such
point exits $S_{I,R_{1},R_{2},\varepsilon}$ through $\left\{ x=R_{2}\right\}
$.
\end{theorem}
\begin{proof} The conditions to apply Theorem \eqref{th:extended-sector} can be validated by subdividing $S_{I,R_{1},R_{2},\varepsilon}$ into $10^5$ pieces along the coordinate $x$ and into $50$ pieces along the coordinate $\phi$ (in total $5\cdot 10^6$ pieces) and by checking the required conditions on each of the ``cubes'' separately. Such computation took under two minutes on a standard laptop.
\end{proof}


\section{Proof of the main Theorem\label{sec:oscillatory-proof}}

Here we shall construct oscillatory motions and, therefore, prove Theorem %
\ref{th:main}.

Throughout this section we will be working under the assumption that we have
a setting in which there is an unstable sector $S_{I,L,R}^{u}$ which is
established by means of Lemma \ref{th:invariant-sector}. Moreover, we shall
assume that the unstable manifold within it is established by means of
Theorem \ref{th:manifold-in-sector-3bp}. We will also assume that the bound
on the manifold is extended to $S_{I,R_{1},R_{2},\varepsilon }$ by means of
Theorem \ref{th:extended-sector}.

Since $S_{I,L,R}^{u}\cup S_{I,R_{1},R_{2},\varepsilon }\subset \{x>0,y<0\}$
by (\ref{eq:dot-x1}) we see that on $S_{I,L,R}^{u}\cup
S_{I,R_{1},R_{2},\varepsilon }$ we have $\dot x>0$. This means that if there
is a point in $S_{I,L,R}^{u}\cup S_{I,R_{1},R_{2},\varepsilon }$ whose
backward trajectory remains in this set, then it has to converge to $\Lambda
_{I}$ (See Figure \ref{fig:extension}).

\begin{remark}
For the energy level $H=I$ with  $I=-1$ we have validated the existence of $%
S_{I,L,R}^{u}$ and of $S_{I,R_{1},R_{2},\varepsilon }$ by means of a
computer assisted (see Theorems \ref{th:manifold-for-given-I} and \ref%
{th:manifold-for-given-I-extended}).
\end{remark}

Due to the symmetry property of system (\ref{eq:symmetry-flow}) from the
bound on the unstable manifold we automatically obtain the bound $\mathcal{S}%
\left( S_{I,L,R}^{u}\cup S_{I,R_{1},R_{2},\varepsilon }\right) $ on the
stable manifold.

The oscillatory motions will be established by the method of covering
relations \cite{MR2060531,MR2276430}. We start by recalling the method, and
then apply it to our problem in the subsequent section.

\subsection{Covering relations}

We restrict to the case of covering relations for maps on $\mathbb{R}^{2}$,
since this is sufficient for our application. This allows us to simplify
some of the introduced tools and notions. The methods from \cite%
{MR2060531,MR2276430} though, are general and can be applied to carry out
analogous constructions in higher dimensional settings.

\begin{definition}
An \emph{h-set} is a pair $\left( N,c_{N}\right) $ where $N\subset\mathbb{R}%
^{2}$ and $c_{N}:\mathbb{R}^{2}\rightarrow\mathbb{R}^{2} $ is a
homeomorphism such that $c_{N}\left( N\right) =\left[ -1,1\right] \times%
\left[ -1,1\right] $.
\end{definition}


For simplicity, when referring to an h-set we will sometimes write only the
set $N$, always assuming implicitly that we have an associated homeomorphism
$c_{N}$ with it.

For an h-set $\left( N,c_{N}\right) $ we define%
\begin{align*}
N_{c} & =\left[ -1,1\right] \times\left[ -1,1\right] , \\
N_{c}^{l} & =\left\{ -1\right\} \times\left[ -1,1\right] , \\
N_{c}^{r} & =\left\{ 1\right\} \times\left[ -1,1\right] , \\
N_{c}^{-} & =N_{c}^{l}\cup N_{c}^{r}, \\
N_{c}^{+} & =\left( \left[ -1,1\right] \times\left\{ -1\right\} \right)
\cup\left( \left[ -1,1\right] \times\left\{ 1\right\} \right) ,
\end{align*}
and%
\begin{equation*}
N^{l}=c_{N}^{-1}\left( N_{c}^{l}\right) ,\qquad N^{r}=c_{N}^{-1}\left(
N_{c}^{r}\right) ,\qquad N^{-}=c_{N}^{-1}\left( N_{c}^{-}\right) ,\qquad
N^{+}=c_{N}^{-1}\left( N_{c}^{+}\right) .
\end{equation*}

In this section we use the notation $\left( \mathrm{x},\mathrm{y}\right) \in
N_{c}\subset\mathbb{R}^{2}$ for the coordinates. We write $\pi_{\mathrm{x}}$
and $\pi_{\mathrm{y}}$ for the projections onto the $\mathrm{x}$ and $%
\mathrm{y}$ coordinates, respectively. Note that we use a different font $%
\mathrm{x}, \mathrm{y}$ in order to distinguish these with the coordinates $%
x,y$ of the PCR3BP.

\begin{figure}[tbp]
\begin{center}
\includegraphics[width=10cm]{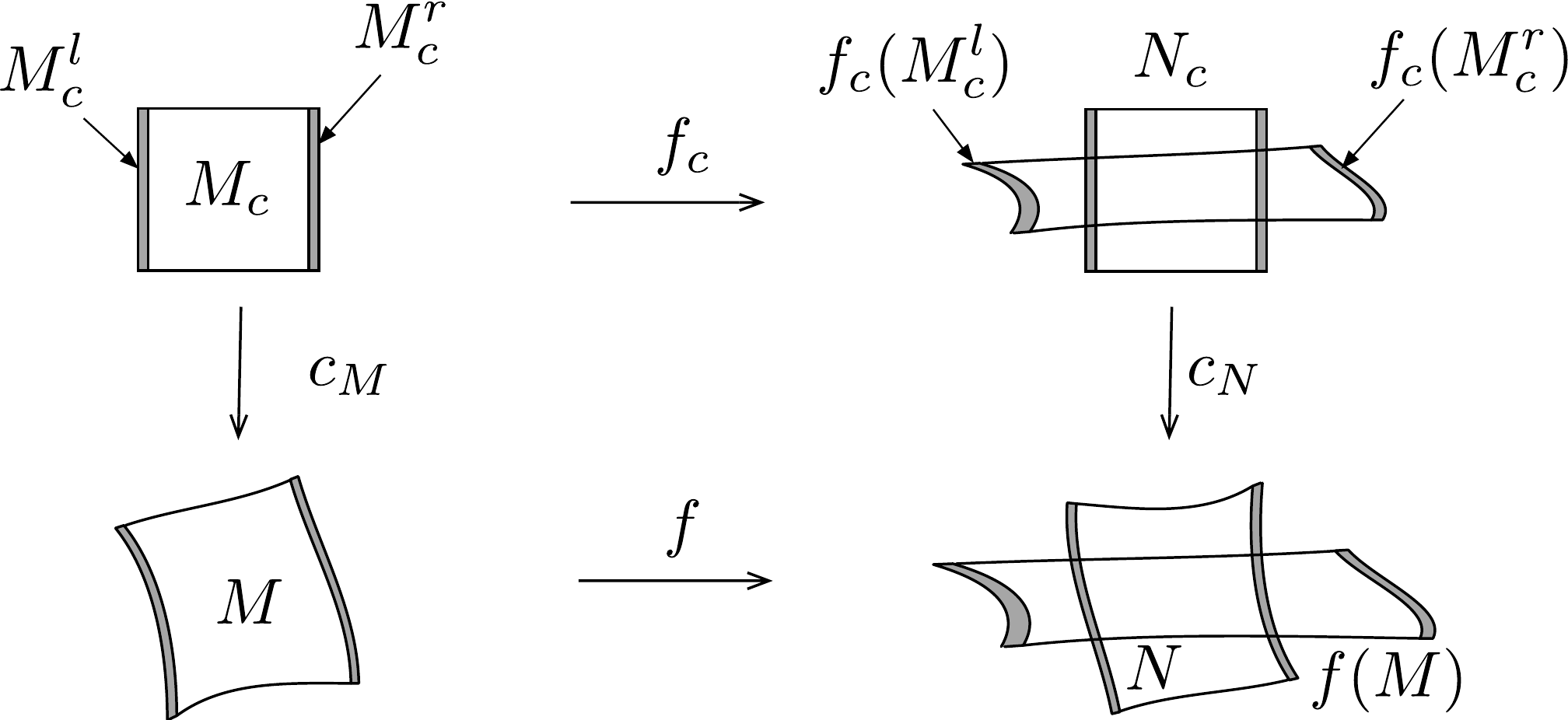}
\end{center}
\caption{The covering $M\protect\overset{f}{\Longrightarrow }N$.}
\end{figure}

\begin{definition}
Let $\left( M,c_{M}\right) $ and $\left( N,c_{N}\right) $ be two h-sets. Let
$f:\mathbb{R}^{2}\rightarrow \mathbb{R}^{2}$ be a continuous map and let $%
f_{c}:=c_{N}\circ f\circ c_{M}^{-1}$. We say that $\left( M,c_{M}\right) $ $%
f $-covers $\left( N,c_{N}\right) $, which we denote as%
\begin{equation*}
M\overset{f}{\Longrightarrow }N,
\end{equation*}%
iff%
\begin{equation}
\pi _{\mathrm{x}}f_{c}\left( M_{c}^{l}\right) <-1\qquad \text{and}\qquad \pi
_{\mathrm{x}}f_{c}\left( M_{c}^{r}\right) >1,
\label{eq:covering-exit-cond-1}
\end{equation}%
or%
\begin{equation}
\pi _{\mathrm{x}}f_{c}\left( M_{c}^{r}\right) <-1\qquad \text{and}\qquad \pi
_{\mathrm{x}}f_{c}\left( M_{c}^{l}\right) >1,
\label{eq:covering-exit-cond-2}
\end{equation}%
and%
\begin{equation}
\pi _{\mathrm{y}}f_{c}\left( M_{c}\right) \cap \left( \left[ -1,1\right]
\times \left( \mathbb{R\setminus }\left( -1,1\right) \right) \right)
=\emptyset .  \label{eq:covering-entry-cond}
\end{equation}
\end{definition}

The coordinate $\mathrm{x}$ plays the role of a local coordinate along which
we have a topological expansion, and $\mathrm{y}$ plays the role of a
coordinate along which we have a topological contraction.

We now introduce the notion of back-covering.

\begin{definition}
\label{def:symmetry-for-back-cover} Let $T:\mathbb{R}^{2}\rightarrow \mathbb{%
R}^{2}$ be defined as $T\left( \mathrm{x},\mathrm{y}\right) =\left( \mathrm{y%
},\mathrm{x}\right) $. For an h-set $\left( N,c_{N}\right) $ we define an
h-set $\left( N^{T},c_{N}^{T}\right) $ as
\begin{align*}
N^{T} & =N, \\
c_{N}^{T} & =T\circ c_{N}.
\end{align*}
\end{definition}

\begin{figure}[tbp]
\begin{center}
\includegraphics[width=7cm]{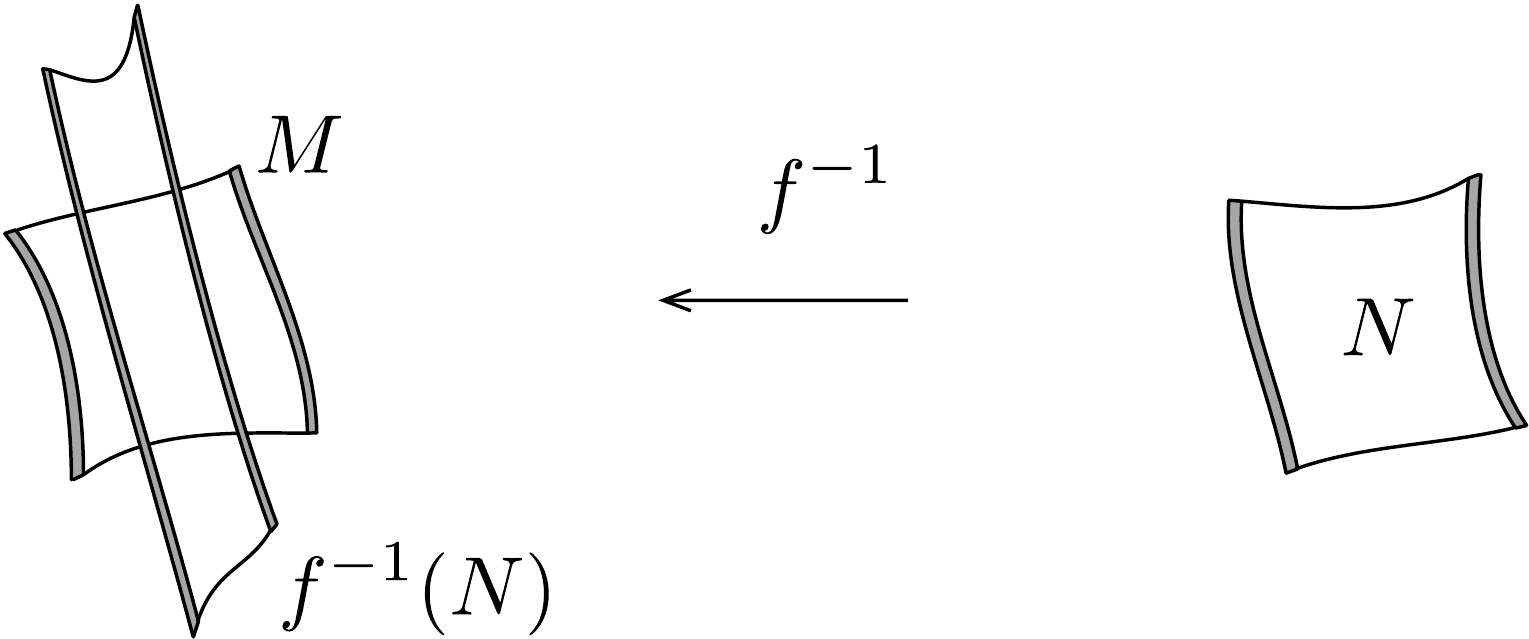}
\end{center}
\caption{The back-covering $M\protect\overset{f}{\Longleftarrow }N$.}
\end{figure}

\begin{definition}
\label{def:NT} Let $\left( M,c_{M}\right) $ and $\left( N,c_{N}\right) $ be
two h-sets. Let $f:\mathbb{R}^{2}\rightarrow\mathbb{R}^{2}$ be a continuous
and such that $f^{-1}:N\rightarrow\mathbb{R}^{2}$ is well defined. We say
that $\left( M,c_{M}\right) $ $f$-back-covers $\left( N,c_{N}\right) $,
which we denote as%
\begin{equation*}
M\overset{f}{\Longleftarrow}N,
\end{equation*}
iff
\begin{equation*}
N^{T}\overset{f^{-1}}{\Longrightarrow}M^{T}.
\end{equation*}
\end{definition}

For our shadowing theorem we also need the following notions.

\begin{definition}
\label{def:hor} Let $N$ be an h-set. Let $h:[-1,1]\rightarrow N$ be
continuous and let $h_{c}=c_{N}\circ h$. We say that $h$ is a horizontal
disk in $N$ if
\begin{equation*}
\pi_{\mathrm{x}}h_{c}(\mathrm{x})=\mathrm{x}.
\end{equation*}
\end{definition}

\begin{definition}
\label{def:ver} Let $N$ be an h-set. Let $v:[-1,1]\rightarrow N$ be
continuous and let $v_{c}=c_{N}\circ v$. We say that $v$ is a vertical disk
in $N$ if
\begin{equation*}
\pi_{\mathrm{y}}v_{c}(\mathrm{y})=\mathrm{y}.
\end{equation*}
\end{definition}

\begin{definition}
Let $N$ be an h-set and $h$ and $v$ be horizontal and vertical discs in $N$,
respectively. By $\left\vert h\right\vert $ and $\left\vert v\right\vert $
we denote the image of $h$ and $v$, respectively.
\end{definition}

The theorem below is our main tool for establishing oscillatory motions.


\begin{theorem}
\cite[Thm. 4]{MR2060531}\cite[Thm. 4]{MR2276430} \label{th:covering} Let $k\geq 1$. Assume that $%
N_{i}$, $i=0,\dots ,k$, are $h$-sets and for each $i=1,\dots ,k$ we have
either
\begin{equation}
N_{i-1}\cover{f_i}N_{i}  \label{eq:tran-dirgcov}
\end{equation}%
or
\begin{equation}
N_{i-1}\invcover{f_i}N_{i}.  \label{eq:tran-invgcov}
\end{equation}
\begin{enumerate}
\item If $N_k=N_0$ then there exists an $x\in \mathrm{int} N_0$ such that 
\[f_i\circ f_{i-1}\circ \ldots \circ f_1(x) \in \mathrm{int} N_i \qquad \mbox{for every } i\in \{1,\ldots, k \}, \]
 and 
\[f_k \circ f_{k-1} \circ \ldots \circ f_1(x) =x. \]
\item Let $h$ be a horizontal disk in $N_{0}$ and $v$ a vertical disk in $N_{k}$.
Then there exists a point $z\in \inter N_{0}$, such that
\begin{align}
z& \in \left\vert h\right\vert ,  \notag \\
f_{i}\circ f_{i-1}\circ \cdots \circ f_{1}(z)& \in \inter N_{i},\quad
i=1,\dots ,k  \label{eq:topop-transv-result} \\
f_{k}\circ f_{k-1}\circ \cdots \circ f_{1}(z)& \in \left\vert v\right\vert .
\notag
\end{align}
\end{enumerate}
\end{theorem}

\begin{corollary}
\label{cor:covering}%
We have the following results for infinite sequences of
coverings:
\begin{enumerate}
\item If%
\begin{equation*}
N_{i-1}\cover{f_i}N_{i}\qquad \text{or\qquad }N_{i-1}\invcover{f_i}%
N_{i}\qquad \text{for }i=1,2\dots
\end{equation*}%
then, for every horizontal disc $h$ in $N_{0}$, there exists a forward
trajectory $z_{0},z_{1},\ldots $, $f_{i}\left( z_{i-1}\right) =z_{i}$, such
that $z_{i}\in \inter N_{i}$ for $i=1,2,\ldots $ and $z_{0}\in \left\vert
h\right\vert $.

\item If%
\begin{equation*}
N_{i-1}\cover{f_i}N_{i}\qquad \text{or\qquad }N_{i-1}\invcover{f_i}%
N_{i}\qquad \text{for }i=0,-1,-2,\dots
\end{equation*}%
then, for every vertical disc $v$ in $N_{0}$, there exists a backward
trajectory $\ldots ,z_{-2},z_{-1},z_{0}$, $f_{i}\left( z_{i-1}\right) =z_{i}$%
, such that $z_{i}\in \inter N_{i}$ for $i=-1,-2,\ldots $ and $z_{0}\in
\left\vert v\right\vert $.

\item If%
\begin{equation*}
N_{i-1}\cover{f_i}N_{i}\qquad \text{or\qquad }N_{i-1}\invcover{f_i}%
N_{i}\qquad \text{for }i\in \mathbb{Z},
\end{equation*}%
there exists full trajectory $\left( z_{i}\right) _{i\in \mathbb{Z}}$, $%
f_{i}\left( z_{i-1}\right) =z_{i}$, such that $z_{i}\in \inter N_{i}$ for $%
i\neq 0$ and $z_{0}\in N_{0}$.
\end{enumerate}
\end{corollary}

\begin{proof}
For item 1, from Theorem \ref{th:covering} and a finite sequence of
coverings%
\begin{equation*}
N_{i-1}\cover{f_i}N_{i}\qquad \text{or\qquad }N_{i-1}\invcover{f_i}%
N_{i}\qquad \text{for }i=0,1,\dots ,k,
\end{equation*}%
we obtain $x_{k}\in \left\vert h\right\vert $, such that a trajectory
starting from $x_{k}$ visits the successive h-sets $N_i$, $i=1\ldots k$. By compactness of $%
\left\vert h\right\vert $, there exists a convergent subsequence $x_{k_{i}}%
\overset{i\rightarrow \infty }{\rightarrow }z_{0}\in \left\vert h\right\vert
$, which proves our claim.

Items 2 and 3 follow from mirror arguments. Item 2 follows by
considering finite sequences with $i\in \left\{ -k,\ldots ,-1,0\right\} $
and Item 3 by considering finite sequences with $i\in \left\{ -k,\dots
,k\right\} $.
\end{proof}

\subsection{Overview of the proof}

The proof will be based on a construction of appropriate h-sets, which will
be positioned on two dimensional sections along the flow. Note that we fix the energy level at $H=I$, which makes the system three dimensional, so
sections transversal to the flow are of dimension two. We consider three
sections on which we will position our h-sets
\begin{eqnarray*}
\Sigma _{0} &=&\left\{ y=0\right\} , \\
\Sigma _{1} &=&\left\{ x=R_{2}\right\} , \\
\Sigma &=&\left\{ \phi =0\right\} ,
\end{eqnarray*}%
where $R_{2}\in \mathbb{R}$ is a fixed constant (as obtained in Theorem \ref%
{th:manifold-for-given-I-extended}). On these sections we consider several
types of h-sets:%
\begin{eqnarray*}
N_{0} &\subset &\Sigma _{0}, \\
N_{1} &\subset &\Sigma _{1}, \\
N\left( r\right) ,N\left( \bar{r},r\right) ,M\left( \alpha ,\beta \right)
&\subset &\Sigma .
\end{eqnarray*}%
(The $r,\bar{r},\alpha ,\beta $ are parameters. The set $N_{0}$ is defined
in (\ref{eq:N0}) and depicted in Figure \ref{fig:cover1}; the set $N_{1}$ is
defined in (\ref{eq:N1}); The set $N\left( r\right) $ is defined in (\ref%
{eq:N2}) and depicted in Figure \ref{fig:N2}; The set $N\left( \bar{r}%
,r\right) $ is defined in (\ref{eq:N3}) and depicted in Figure \ref{fig:N3};
The set $M\left( \alpha ,\beta \right) $ is defined in (\ref{eq:M-h-set})
and depicted in Figure \ref{fig:M-h-set}. The smaller the $r,\bar{r},\alpha
,\beta $ the closer the h-sets are to $\Lambda _{I}$.)

Set $N_{0}$ is self $\mathcal{S}$-symmetric (see Figure \ref{fig:cover1}).
Checking that%
\begin{equation*}
N_{1}\cover{}N_{0}
\end{equation*}%
involves computer assisted validation. The rest of coverings is proven by
analytic arguments. The smaller the parameter $r$, the closer the left side
of the exit set of $N\left( r\right) $ is to $\Lambda _{I}$; see Figure \ref%
{fig:N2}. This will allow us to prove in Lemma \ref{lem:from-N2-to-N1} that
if $r$ is sufficiently small then
\begin{equation*}
N\left( r\right) \cover{}N_{1}.
\end{equation*}%
In Lemma \ref{lem:from-N3-to-N2} we prove that%
\begin{equation*}
N\left( \bar{r},r\right) \cover{}N\left( r\right) ,
\end{equation*}%
provided that we choose appropriately small $\bar{r}$. We also show in
Lemmas \ref{lem:from-M-to-Nr}, \ref{lem:to-M-alpha-beta} that we have a
sequence of coverings%
\begin{equation*}
M\left( \bar{r},r\right) \cover{}M\left( \alpha ,\beta \right) \cover{}%
N\left( \bar{r},r\right) ,
\end{equation*}%
for arbitrarily small $\alpha ,\beta $. The smaller we choose $\alpha ,\beta
$ the closer to $\Lambda _{I}$ is the h-set $M\left( \alpha ,\beta \right) $%
. This means that we can obtain an orbit that passes through the above
sequence of coverings to approach arbitrarily close to $\Lambda _{I}$. The
set $M\left( \bar{r},r\right) $ is an $\mathcal{S}$-symmetric to $N\left(
\bar{r},r\right) $, this allows us to automatically obtain coverings from $%
N_{0}$ to $M\left( \bar{r},r\right) $, this way we obtain a sequence of
coverings%
\begin{equation}
N_{0}\cover{}\ldots \cover{}M\left( \alpha ,\beta \right) \cover{}\ldots %
\cover{}N_{0}.  \label{eq:full-connection}
\end{equation}%
The smaller the $\alpha ,\beta $ the closer is the approach to $\Lambda _{I}$
for orbits that pass through such sequences. We then choose a sequence $%
M\left( \alpha _{1},\beta _{1}\right) ,M\left( \alpha _{2},\beta _{2}\right)
,\ldots $ and glue sequences (\ref{eq:full-connection}) to obtain
oscillatory motions.

To prove bounded motions we will glue infinite sequences (\ref%
{eq:full-connection}) with fixed $\alpha =\bar{r}$ and $\beta =r$.

To prove parabolic and hyperbolic motions, we will choose appropriate
horizontal discs in $N\left( \bar{r},r\right) $ and vertical discs in $M\left( \bar{r},r\right) $ and  use the
connection%
\begin{equation*}
N\left( \bar{r},r\right) \cover{}\ldots \cover{}N_{0}\cover{}\ldots \cover{}%
M\left( \bar{r},r\right) .
\end{equation*}%
The discs will be chosen so that when an orbit passes through a given disc, it escapes to infinity (i.e. $x=0$ in McGehee variables \eqref{def:McGehee}) according to a given type of motion (hyperbolic or parabolic).

\subsection{Proof of the main theorem}

Throughout this section we consider
\begin{equation*}
I=-1,\qquad L=4\cdot10^{-9},\qquad R=10^{-4},\qquad R_{1}=\frac{R}{2},\qquad
R_{2}=0.4,\qquad\varepsilon=2\cdot10^{-5},
\end{equation*}
as in Theorems \ref{th:manifold-for-given-I} and \ref%
{th:manifold-for-given-I-extended}. We focus on the case $I=-1$, but the
method is general and can be applied to different energy levels (i.e. values of the Jacobi constant).


Our objective is to construct a sequence of h-sets, positioned on sections
along the flow of the PCR3BP. Our objective is to define these h-sets, so
that Theorem \ref{th:covering} will lead to the existence of orbits which
can approach arbitrarily close to infinity and come back to the regions of
the primaries.

Let $\mathcal{G}_{I}\left( x,y,\phi\right) $ be defined as%
\begin{equation}
\mathcal{G}_{I}\left( x,y,\phi\right) :=\frac{1-\sqrt{1-\frac{x^{4}}{2}%
\left( \frac{1}{2}y^{2}-\mathcal{U}(x,\phi)-I\right) }}{\frac{x^{4}}{4}}
\label{eq:G-implicit-xy}
\end{equation}
a solution of the quadratic equation%
\begin{equation*}
I=\frac{x^{4}G^{2}}{8}+\frac{1}{2}y^{2}-\mathcal{U}(x,\phi)-G.
\end{equation*}
Note that $\lim_{(x,y)\to (0,0)}\mathcal{G}_{I}(x,y,\phi)=-I$ and
\begin{equation}
\mathcal{G}_{I}\left( x,y,\phi\right) =\mathcal{G}_{I}\left( x,-y,-\phi
\right) .  \label{eq:G-implicit-symmetry}
\end{equation}

We consider the first h-set on the section $\Sigma_{0}=\left\{y=0\right\} $.
We consider the flow restricted to energy level $H=I$, which means that
$\Sigma_{0}$ is two dimensional. The points on this section are parametrised
by $x$ and $\phi$, since the coordinate $G$ is determined by (\ref%
{eq:G-implicit-xy}). We will therefore specify an h-set on $\Sigma_{0}$ in
coordinates $x,\phi$.

To define our first h-set $N_{0}$ we consider two matrices $S_{\alpha ,\beta
}$ and $\mathcal{R}$ as%
\begin{equation*}
S_{\alpha ,\beta }:=\mathrm{diag}\left( \alpha ,\beta \right) ,\qquad
\mathcal{R}:=\left(
\begin{array}{cc}
\cos \left( \frac{\pi }{4}\right) & -\sin \left( \frac{\pi }{4}\right) \\
\sin \left( \frac{\pi }{4}\right) & \cos \left( \frac{\pi }{4}\right)%
\end{array}%
\right)
\end{equation*}%
with which we define
\begin{align}
N_{0}& :=S_{\alpha ,\beta }\circ \mathcal{R}\left( \left[ -1,1\right] \times %
\left[ -1,1\right] \right) +\left( x_{0},\pi \right) ,  \label{eq:N0} \\
c_{N_{0}}\left( z\right) & :=\mathcal{R}^{-1}S_{\alpha ,\beta }^{-1}\left(
z-\left( x_{0},\pi \right) \right) ,  \notag
\end{align}%
where $\alpha ,\beta ,x_{0}$ are given parameters. We have found that good
choices of the parameters are
\begin{equation*}
\alpha =2\cdot 10^{-4},\quad \beta =10^{-1},\quad x_{0}=1.9999.
\end{equation*}

\begin{figure}[ptb]
\begin{center}
\includegraphics[width=2.5in]{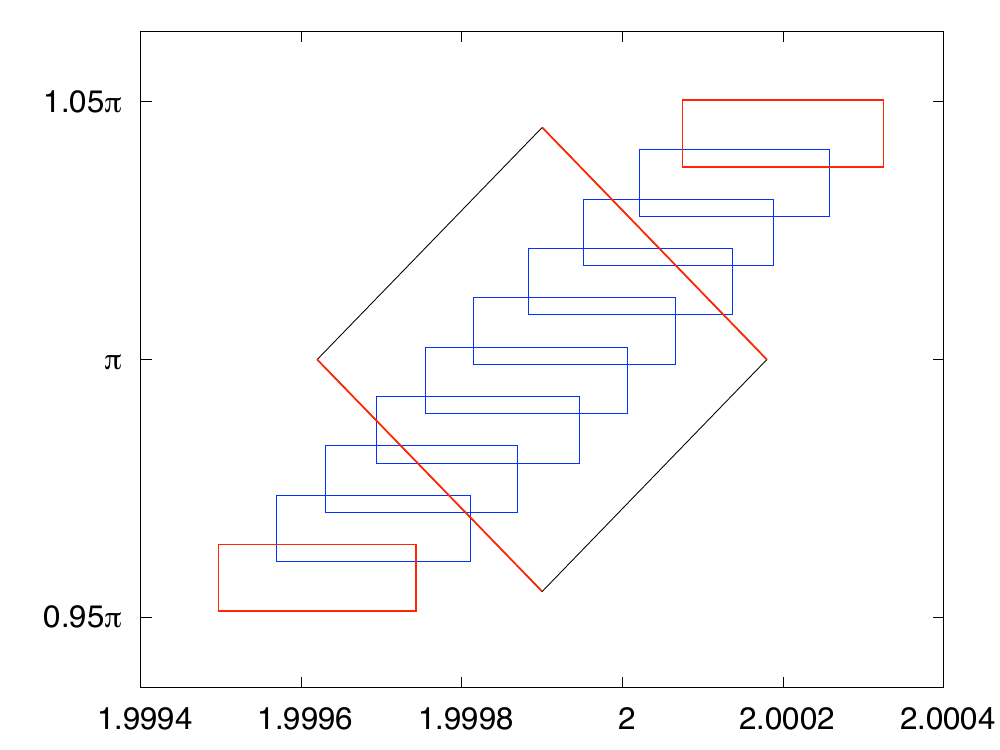}
\end{center}
\caption{The covering $N_{1} \protect\overset{f_{0}}{\Longrightarrow} N_{0}$%
. The set $N_{0}$ is the rhombus, with the exit set marked in red. The
computer assisted bound on the image of $N_{1}$ is depicted by the
diagonally placed rectangles, and the bound on the exit set is depicted in
red. The figure is in the $(x,\protect\phi)$ coordinates, with $x$ on the
horizontal axis and $\protect\phi$ on the vertical axis.}
\label{fig:cover1}
\end{figure}

The choice of $x_{0}$ is motivated by the fact that the point $\left(
x_{0},\pi \right) $ is (roughly) on the intersection of the stable and
unstable manifolds at $\Sigma _{0}$. The set $N_{0}$ is constructed by
rotating counterclockwise the square $[-1,1]^{2}$ by the angle $\frac{\pi }{4%
}$, rescaling it, and shifting it to be centered at $\left( x_{0},\pi
\right) $. This way we obtain a rhombus depicted in Figure \ref{fig:cover1}.
The scaling coefficients $\alpha ,\beta $ are chosen so that the edges $%
N^{+} $ are (roughly) aligned with the intersection of the unstable manifold
with $\Sigma _{0}$. In Figure \ref{fig:cover1} the (rigorous,
computer-assisted) bound on this intersection is depicted in blue.

\begin{remark}
The set $N_{0}$ is $\mathcal{S}$ symmetric. To be more precise, if $\left(
x,y=0,\phi,G\right) \in\Sigma_{0}$ is a point in $N_{0}$ then%
\begin{equation*}
\mathcal{S}\left( x,0,\phi,G\right) =\left( x,0,-\phi,G\right)
\end{equation*}
lies in $\left\{ y=0\right\} $. The point $\left( x,-\phi\right) $ lies in
the rhombus, since $-\phi$ is identified with $-\phi+2\pi$, and the rhombus
is symmetric with respect to the line $\phi=\pi$ (see Figure \ref{fig:cover1}%
), and from (\ref{eq:G-implicit-symmetry}) we see that $\mathcal{S}\left(
x,0,\phi,G\right) \in N_{0}\subset\Sigma_{0}$.
\end{remark}



Our second h-set $N_{1}$ is contained in the section%
\begin{equation*}
\Sigma _{1}=\left\{ x=R_{2}\right\} .
\end{equation*}%
The points on $\Sigma _{1}$ are parameterised by $\phi ,y$, since $G$ can be
computed as $G=\mathcal{G}_{I}\left( R_{2},y,\phi \right) .$ On the section $%
\Sigma _{1}$ the role of the exit coordinate is played by $\phi $, and the
topologically entry coordinate is $y$. We consider the set $N_{1}\subset
\Sigma _{1}$, defined on the $\phi ,y$ coordinates as%
\begin{equation}
N_{1}=\left[ \phi _{1},\phi _{2}\right] \times \left[ \gamma _{I}\left(
R_{2}\right) -\varepsilon ,\gamma _{I}\left( R_{2}\right) +\varepsilon %
\right] ,  \label{eq:N1}
\end{equation}%
where good choices of parameters $\phi _{1},\phi _{2}$ are
\begin{equation*}
\phi _{1}=3.45\qquad \text{and\qquad }\phi _{2}=3.75.
\end{equation*}%
We consider $f_{0}:\Sigma _{1}\rightarrow \Sigma _{0}$ to be the section to
section map along the flow of the PCR3BP. With the aid of computer assisted
computation we validate the following result.

\begin{lemma}
\label{lem:from-N1-to-N0} For $I=-1$ and the above defined h-sets $%
N_{0},N_{1}$ and the section to section map $f_{0}:\Sigma_{1}\rightarrow%
\Sigma_{0}$ we have%
\begin{equation*}
N_{1}\overset{f_{0}}{\Longrightarrow}N_{0}.
\end{equation*}
\end{lemma}

\begin{proof}
The computer assisted bounds that validate Lemma \ref{lem:from-N1-to-N0} are
depicted in Figure \ref{fig:cover1}. The computation took under eight
seconds on a standard laptop.
\end{proof}

\begin{corollary}
\label{cor:from-N0-to-SN1}Let $\tilde{f}_{0}:\Sigma_{0}\rightarrow\Sigma_{1}$
be section to section map along the flow of the PCR3BP (for positive time).
From the time reversing symmetry (\ref{eq:symmetry-flow}) we obtain that
(recall that $N^{T}$ was defined in Definition \ref{def:NT})
\begin{equation*}
N_{0}\overset{\tilde{f}_{0}}{\Longleftarrow}\mathcal{S}\left(
N_{1}^{T}\right) .
\end{equation*}
\end{corollary}

\begin{remark}
The computation for the proof of Lemma \ref{lem:from-N1-to-N0}, see Figure %
\ref{fig:cover1}, shows that the stable and unstable manifolds of $\Lambda
_{I}$ intersect. This follows from the fact that these manifolds are $%
\mathcal{S}$-symmetric.
\end{remark}

\begin{figure}[tbp]
\begin{center}
\includegraphics[width=2.5in]{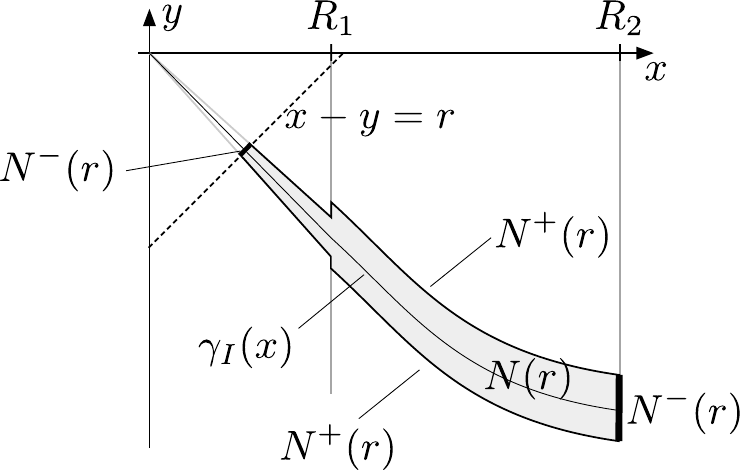}
\end{center}
\caption{The h-set $N(r)$.}
\label{fig:N2}
\end{figure}

We choose our next h-set on the section
\begin{equation*}
\Sigma =\left\{ \phi =0\right\} .
\end{equation*}%
The points on $\Sigma $ are parameterised by $x,y$, since $G$ can be
computed as $G=\mathcal{G}_{I}\left( x,y,0\right) $. On $\Sigma $ we can
therefore consider h-sets expressed in coordinates $x,y$. The coordinate $y$
is entry direction and $x$ is the exit direction. We define the following
h-set, $N\left( r\right) $, where $r\in \left( 0,R\right) $ and (for the
intuition behind the definition see Figure \ref{fig:N2} and the explanation
that comes in the next paragraph)
\begin{align}
N\left( r\right) & :=\left\{ \left( x,y\right) :\text{ }x\in \left[
R_{1},R_{2}\right] \text{ and }y\in \left[ \gamma _{I}\left( x\right)
-\varepsilon ,\gamma _{I}\left( x\right) +\varepsilon \right] \right\}
\label{eq:N2} \\
& \quad \cup \left\{ \left( x,y\right) :x<R_{1},\text{ }x-y\geq r\text{ and }%
\left\vert x+y\right\vert \leq L\left( x-y\right) \right\} .  \notag
\end{align}

The reason for above definition of $N(r)$ is the following. For $x\in \left[
R_{1},R_{2}\right] $, from the definition of $S_{I,R_{1},R_{2},\varepsilon }$
we know that the unstable manifold on $\Sigma $ is enclosed in an $%
\varepsilon $ neighbourhood of the curve $\left( x,\gamma _{I}\left(
x\right) \right) $. We therefore position our h-set around this curve for $%
x\in \left[ R_{1},R_{2}\right] $. Recall also that the coordinates $q,p$
were defined as $q=\frac{1}{2}(x-y)$ and $p=\frac{1}{2}(x+y)$, so the
condition $\left\vert x+y\right\vert \leq L\left\vert x-y\right\vert $
ensures that $\left\vert p\right\vert \leq Lq$, so the points in $N(r)$ that
have $x<R_{1}$ lie in the sector $S_{I,L,R}^{u}$.

The sets $N^{\pm }(r)$ are depicted on Figure \ref{fig:N2}. The formal
definition is follows
\begin{align*}
N^{+}\left( r\right) & =\left\{ \left( x,y\right) :\text{ }x\in \left[
R_{1},R_{2}\right] \text{ and }y=\gamma _{I}\left( x\right) -\varepsilon
\right\} \\
& \quad \cup \left\{ \left( x,y\right) :\text{ }x\in \left[ R_{1},R_{2}%
\right] \text{ and }y=\gamma _{I}\left( x\right) +\varepsilon \right\} \\
& \quad \cup \left\{ \left( x,y\right) :\text{ }x=R_{1},\text{ }-x\left(
1-L\right) \left( 1+L\right) ^{-1}\leq y\text{ and }y\leq \gamma _{I}\left(
x\right) +\varepsilon \right\} \\
& \quad \cup \left\{ \left( x,y\right) :\text{ }x=R_{1},\text{ }\gamma
_{I}\left( x\right) -\varepsilon \leq y\text{ and }y\leq -x\left( 1+L\right)
\left( 1-L\right) ^{-1}\right\} \\
& \quad \cup \left\{ \left( x,y\right) :\text{ }x\leq R_{1},\text{ }x-y\geq
r,\text{ and }\left\vert x+y\right\vert =L\left( x-y\right) \right\} , \\
N^{-}\left( r\right) & =\left\{ \left( x,y\right) :\text{ }x=R_{2}\text{ and
}y\in \left[ \gamma _{I}\left( x\right) -\varepsilon ,\gamma _{I}\left(
x\right) +\varepsilon \right] \right\} \\
& \quad \cup \left\{ \left( x,y\right) :\text{ }x-y=r,\text{ and }\left\vert
x+y\right\vert \leq L\left( x-y\right) \right\} .
\end{align*}

Let us assume that for every point from the set
\begin{equation*}
\left\{ \left( x,y,\phi ,G\right) :\left( x,y\right) \in N\left( r\right)
,\phi \in \mathbb{S}^{1},G=\mathcal{G}_{I}(x,y,\phi )\right\}
\end{equation*}%
we have (see (\ref{eq:dot-x1}))%
\begin{equation}
\dot{\phi}<0.  \label{eq:phi-decreasing}
\end{equation}

\begin{remark}
For a given fixed $I$ (\ref{eq:phi-decreasing}) can be validated with the
aid of a computer and interval arithmetic.
\end{remark}

Let
\begin{equation*}
f_{1}:\Sigma \rightarrow \Sigma _{1}
\end{equation*}%
be the section to section map along the flow.

\begin{lemma}
\label{lem:from-N2-to-N1}There exists an $r_{0}>0$ such that for any $r\in
\left( 0,r_{0}\right) $%
\begin{equation*}
N\left( r\right) \overset{f_{1}}{\Longrightarrow }N_{1}.
\end{equation*}
\end{lemma}

\begin{proof}
The flow starting from any point in $N\left( r\right) $ will remain in the
interior of $S_{I,L,R}^{u}\cup S_{I,R_{1},R_{2},\varepsilon }$ until it
exits this set through $\Sigma _{1}=\left\{ x=R_{2}\right\}$.
This means
that the trajectory will not pass through $N_{1}^{+}=\{y=\gamma_I(R_2)\pm%
\varepsilon, \phi\in\left[ \phi _{1},\phi _{2}\right]\}$ (see \eqref{eq:N1}%
), so
\begin{equation*}
\pi _{\mathrm{y}}f_{1}\left( N(r)\right) \cap N_{1}^{+}=\emptyset ,
\end{equation*}%
which validates condition (\ref{eq:covering-entry-cond}) from the definition
of the covering relation.

Let us now observe that%
\begin{equation*}
N^{-}(r)\cap \left\{ x=R_{2}\right\} \in \Sigma _{1}.
\end{equation*}%
This means that
\begin{equation*}
f_{1}\left( p\right) =p,\qquad \text{for every }p\in N^{-}(r)\cap \left\{
x=R_{2}\right\} ,
\end{equation*}%
which implies%
\begin{equation}
\pi _{\phi }f_{1}\left( p\right) =0\qquad \text{for every }p\in N^{-}(r)\cap
\left\{ x=R_{2}\right\} .  \label{eq:N2-phi-0-image}
\end{equation}%
Observe that since $\phi _{1},\phi _{2}$ used to define $N_{1}$ are such
that zero is not within the interval $\left[ \phi _{1},\phi _{2}\right] $.
This means that from (\ref{eq:N2-phi-0-image}) we obtain
\begin{equation*}
f_{1}\left( N^{-}(r)\cap \left\{ x=R_{2}\right\} \right) \cap
N_{1}=\emptyset .
\end{equation*}

Let us now consider $\phi $ to be from $\mathbb{R}$ instead of $\mathbb{S}%
^{1}.$ In other words, we consider a lift of the angle to the real line.
From (\ref{eq:phi-decreasing}) we see that $\pi _{\phi }f_{1}\left(
N(r)\right) \in \mathbb{R}_{-}$. We can identify $N_{1}$ with $\left[ \phi
_{1}-2\pi ,\phi _{2}-2\pi \right] \times \lbrack \gamma
_{I}(R_{2})-\varepsilon ,\gamma _{I}(R_{2})+\varepsilon ]$ so the angles of $%
N_{1}$ are negative. From (\ref{eq:N2-phi-0-image}) we see that
\begin{equation*}
\pi _{\phi }f_{1}\left( N^{-}(r)\cap \left\{ x=R_{2}\right\} \right) >\pi
_{\phi }N_{1},
\end{equation*}%
which is the second inequality of (\ref{eq:covering-exit-cond-2}), from the
definition of the covering relation.

We now need to show that if we choose $r>0$ to be sufficiently small, then
we will obtain%
\begin{equation*}
\pi _{\phi }f_{1}\left( N^{-}(r)\setminus \left\{ x=R_{2}\right\} \right)
<\pi _{\phi }N_{1}.
\end{equation*}%
This establishes the second inequality from (\ref{eq:covering-exit-cond-2}).
The smaller the $r$ we choose, the closer to the origin, which by $\dot{\phi}
$ in (\ref{eq:dot-x1}) implies that we will have a larger change in $\phi $
until we reach $\left\{ x=R_{2}\right\} $ from $\left\{ x-y=r\right\} .$
Choosing small $r$ we can therefore obtain $\pi _{\phi }f_{1}\left(
N^{-}(r)\setminus \left\{ x=R_{2}\right\} \right) <\phi _{1}-2\pi $, which
concludes our proof.
\end{proof}

\begin{corollary}
\label{cor:from-SN1-to-SN2}Let $\tilde{f}_{1}:\Sigma _{1}\rightarrow \Sigma $
be the section to section map along the flow of the PCR3BP. From the time
reversing symmetry (\ref{eq:symmetry-flow}) and from Lemma \ref%
{lem:from-N2-to-N1} we obtain that
\begin{equation*}
\mathcal{S}\left( N_{1}^{T}\right) \overset{\tilde{f}_{1}}{\impliedby }%
\mathcal{S}\left( N^{T}\left( r\right) \right) .
\end{equation*}
\end{corollary}

\begin{figure}[tbp]
\begin{center}
\includegraphics[width=2in]{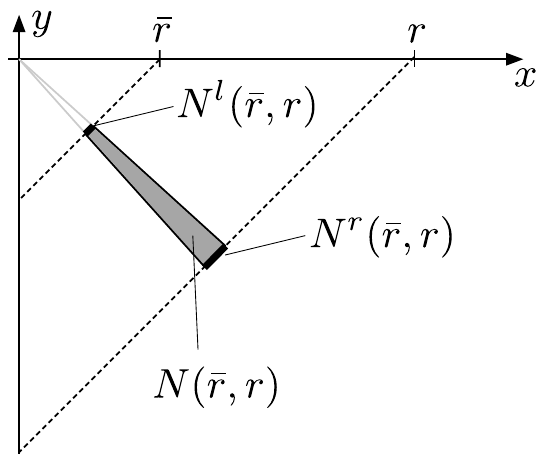}
\end{center}
\caption{The h-set $N(\bar{r},r)$.}
\label{fig:N3}
\end{figure}

Let $\bar{r}\in \left( 0,r\right) $. Recall that the section $\Sigma
=\left\{ \phi =0\right\} $ is parameterised by coordinates $x,y$. We define
the following h-set $N(\bar{r},r)$ in $\Sigma $ as (see Figure \ref{fig:N3}%
).
\begin{equation}
N(\bar{r},r)=\left\{ (x,y):x-y\in \left[ \bar{r},r\right] \text{ and }%
\left\vert x+y\right\vert \leq L\left( x-y\right) \right\} ,  \label{eq:N3}
\end{equation}%
and
\begin{align*}
N^{+}(\bar{r},r)& =N(\bar{r},r)\cap \left\{ (x,y):\left\vert x+y\right\vert
=L\left( x-y\right) \right\} , \\
N^{l}& =N(\bar{r},r)\cap \left\{ x-y=\bar{r}\right\} , \\
N^{r}& =N(\bar{r},r)\cap \left\{ x-y=r\right\} .
\end{align*}

Our objective now will be to show that if $\bar{r}$ is chosen to be
sufficiently small, then we can construct a covering%
\begin{equation*}
N\left( \bar{r},r\right) \Longrightarrow N\left( r\right) .
\end{equation*}%
The above statement is vague, since we have not specified which map we
consider for the covering. We make this more precise in the discussion that
follows.

We first note that for every point in $S_{I,R_{1},R_{2},\varepsilon}$ we
have $\dot x>0$. This means that there exists a $\delta>0$, such that a
trajectory that exits $S_{I,R_{1},R_{2},\varepsilon}$ will reach $\left\{
x=R_{2}+\delta\right\} $ without re-entering $S_{I,R_{1},R_{2},\varepsilon} $%
. Let such $\delta$ be fixed from now on.

Let $\tau :\mathbb{R}^{2}\times \mathbb{S}^{1}\times \mathbb{R\rightarrow R}$
be defined as
\begin{equation*}
\tau \left( \mathbf{x}\right) =\min \left( \inf \left\{ t>0:\pi _{\phi }\Phi
_{t}(\mathbf{x})=0\right\} ,\inf \left\{ t\geq 0:\pi _{x}\Phi _{t}(\mathbf{x}%
)=R_{2}+\delta \right\} \right) ,
\end{equation*}%
and let $f:\Sigma \rightarrow \Sigma $ be defined as
\begin{equation*}
f(\mathbf{x}):=\left( \pi _{x,y}\Phi _{\tau (\mathbf{x})}(\mathbf{x}),0,%
\mathcal{G}_{I}\left( \pi _{x,y}\Phi _{\tau (\mathbf{x})}(\mathbf{x}%
),0\right) \right) .
\end{equation*}

\begin{remark}
\label{rem:map-gives-true-trajectory}From the definition of $\tau$ we see
that if $\pi_{x}f(\mathbf{x})<R_{2}+\delta$, then $f(\mathbf{x})=\Phi _{\tau(%
\mathbf{x})}(\mathbf{x})$. This means that every point $\mathbf{x}$ in $%
\left\{ x<R_{2}+\delta\right\} $ whose image $f(\mathbf{x})$ lands in $%
\left\{ x<R_{2}+\delta\right\},$ the map $f$ corresponds to a true
trajectory of the PCR3BP.
\end{remark}

The definition of $f$ is somewhat artificial, but there is a reason for
which we choose it this way. When we take $\mathbf{x}\in S_{I,L,R}^{u}\cup
S_{I,R_{1},R_{2},\varepsilon }$, then from the way we have defined $f$ each
iterate $f^{n}(\mathbf{x})$ is well defined for every $n\in \mathbb{N}$.
Moreover, as long as $\pi _{x}f^{i}(\mathbf{x})\leq R_{2}$, for $i=0,\ldots
,n$, by Remark \ref{rem:map-gives-true-trajectory} the points $f^{i}(\mathbf{%
x})$ lie on a trajectory of the PCR3BP.

\begin{lemma}
\label{lem:from-N3-to-N2}For every small enough $r>0$ there exists a $k\in
\mathbb{N}$ and $\bar{r}\in \mathbb{R}$ such that
\begin{equation}
N\left( \bar{r},r\right) \overset{f^{k}}{\Longrightarrow }N\left( r\right) .
\label{eq:covering-f3}
\end{equation}
\end{lemma}

\begin{proof}
All trajectories which start from $S_{I,R,L}\cup
S_{I,R_{1},R_{2},\varepsilon }$ will exit this set through $%
S_{I,R_{1},R_{2},\varepsilon }\cap \{x=R_{2}\}$. This means that there
exists a $k$ such that $\pi _{x}f^{k}(N^{r}(\bar{r},r))>R_{2}$. This means
that the first inequality from the condition (\ref{eq:covering-exit-cond-1})
needed for (\ref{eq:covering-f3}) will be satisfied. The smaller the $\bar{r}
$, the closer the points from $N^{l}(\bar{r},r)$ will be to the origin,
where the dynamics is slow. This means that by choosing $\bar{r}$
sufficiently small, we will obtain $\pi _{x}f^{k}(N^{l}(\bar{r},r))<\pi
_{x}N(r)$, which means that the second inequality from the condition (\ref%
{eq:covering-exit-cond-1}) is also satisfied.

We need to show that $f^{k}(N(\bar{r},r))\cap N^{+}\left( r\right)
=\emptyset .$ This follows from the fact that points can exit $S_{I,R,L}\cup
S_{I,R_{1},R_{2},\varepsilon }$ only through $S_{I,R_{1},R_{2},\varepsilon
}\cap \{x=R_{2}\}$. No trajectory which starts in $S_{I,R,L}\cup
S_{I,R_{1},R_{2},\varepsilon }$ can therefore pass through $N^{+}\left(
r\right) $.
\end{proof}

Let $\tilde{\tau}:\mathbb{R}^{2}\times\mathbb{S}^{1}\times\mathbb{%
R\rightarrow R}$ be defined as
\begin{equation*}
\tilde{\tau}\left( \mathbf{x}\right) =\max\left( \sup\left\{ t<0:\pi
_{\phi}\Phi_{t}(\mathbf{x})=0\right\} ,\sup\left\{ t\leq0:\pi_{x}\Phi _{t}(%
\mathbf{x})=R_{2}+\delta\right\} \right) ,
\end{equation*}
and let $\tilde{f}:\Sigma\rightarrow\Sigma$ be defined as%
\begin{equation*}
\tilde{f}(\mathbf{x}):=\left( \pi_{x,y}\Phi_{\tilde{\tau}(\mathbf{x})}(%
\mathbf{x}),0,\mathcal{G}_{I}\left( \pi_{x,y}\Phi_{\tilde{\tau}(\mathbf{x})}(%
\mathbf{x}),0\right) \right) .
\end{equation*}

\begin{corollary}
\label%
{cor:from-SN2-to-SN3}From the time reversing symmetry (\ref{eq:symmetry-flow}%
) and Lemma \ref{lem:from-N3-to-N2} we obtain
\begin{equation*}
\mathcal{S}\left( N^{T}\left( r\right) \right) \overset{\tilde{f}^{k}}{%
\Longleftarrow }\mathcal{S}\left( N^{T}\left( \bar{r},r\right) \right) .
\end{equation*}
\end{corollary}

\begin{figure}[tbp]
\begin{center}
\includegraphics[width=2in]{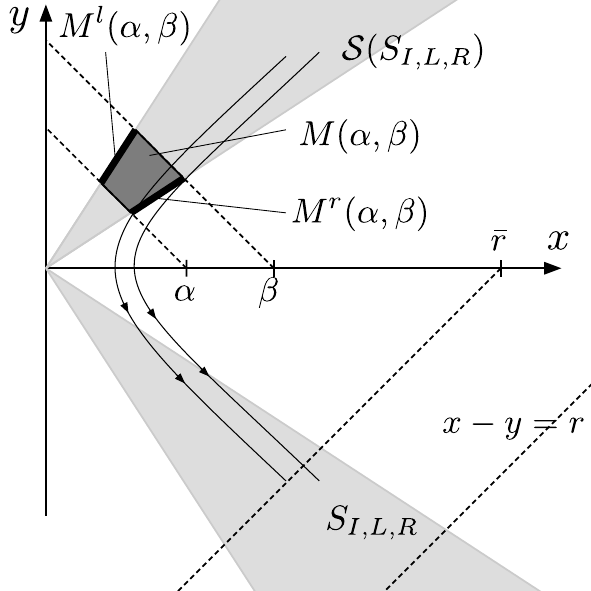}\qquad %
\includegraphics[width=2in]{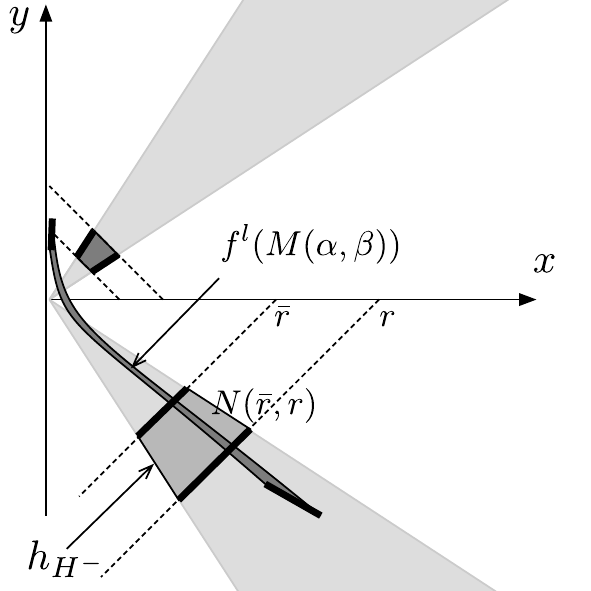}
\end{center}
\caption{The h-set $M(\protect\alpha ,\protect\beta )$ is in dark grey on
the left hand side plot. We see that the exit set $M^{r}(\protect\alpha ,%
\protect\beta )$ will enter the sector $S_{I,R,L}$. On the right is a sketch
of the covering of $N\left( \bar{r},r\right) $ by $M(\protect\alpha ,\protect%
\beta ).$}
\label{fig:M-h-set}
\end{figure}
Let us now consider h-sets on $\Sigma =\left\{ \phi =0\right\} $ of the form
(see Figure \ref{fig:M-h-set})%
\begin{align}
M(\alpha ,\beta )& =\left\{ (x,y):\left\vert x-y\right\vert \leq L\left(
x+y\right) ,x+y\in \left[ \alpha ,\beta \right] \right\}  \label{eq:M-h-set}
\\
M^{r}(\alpha ,\beta )& =M(\alpha ,\beta )\cap \left\{ (x,y):x-y=L\left(
x+y\right) \right\} ,  \notag \\
M^{l}(\alpha ,\beta )& =M(\alpha ,\beta )\cap \left\{ (x,y):x-y=-L\left(
x+y\right) \right\} ,  \notag \\
M^{+}(\alpha ,\beta )& =M(\alpha ,\beta )\cap \left( \{(x,y):x+y=\alpha
\}\cup \{(x,y):x+y=\beta \}\right) .  \notag
\end{align}

\begin{lemma}
\label{lem:from-M-to-Nr}If $B\in \mathbb{R}$, $B>0$, is sufficiently close
to zero, then for any $\alpha <\beta \leq B$ there exists an $\ell $
(depending on the choice of $\alpha ,\beta $) such that%
\begin{equation*}
M(\alpha ,\beta )\overset{f^{\ell }}{\Longrightarrow }N\left( \bar{r}%
,r\right) .
\end{equation*}
\end{lemma}

\begin{proof}
By Proposition \ref{th:orbits-enter-Scu}, if $\beta <r_{\ast }$ (see
statement of Proposition \ref{th:orbits-enter-Scu} for the constant $r_{\ast
}$) a trajectory starting from $M^{r}(\alpha ,\beta )$ will enter $%
S_{I,L,R}^{u}$; see Figure \ref{fig:M-h-set}-left. If $r_{\ast }$ is
sufficiently small, then the trajectory will enter $S_{I,L,R}^{u}\cap
\left\{ x-y<\bar{r}\right\} $. Moreover, such trajectory will exit $%
S_{I,L,R}^{u}\cap \left\{ x-y\leq r\right\} $ through $S_{I,L,R}^{u}\cap
\left\{ x-y=r\right\} .$ This means that there exists an $\ell >0$ such that
$f^{\ell }\left( M^{r}(\alpha ,\beta )\right) \cap N\left( \bar{r},r\right)
=\emptyset $. Moreover, $f^{\ell }\left( M^{r}(\alpha ,\beta )\right) $ will
be `to the right', along the coordinate $x$, of the set $N\left( \bar{r}%
,r\right) $ so the second condition from (\ref{eq:covering-exit-cond-1}) in
the definition of the covering will be fulfilled.

The set $M^{l}(\alpha ,\beta )$ gets mapped outside and above (with respect
to $y$) the set $\mathcal{S}\left( S_{I,L,R}^{u}\right) $; see Figure \ref%
{fig:M-h-set}. By the symmetry of the system a trajectory which starts from $%
M^{l}(\alpha ,\beta )$ will never re-enter $\mathcal{S}\left(
S_{I,L,R}^{u}\cup S_{I,R_{1},R_{2},\varepsilon }\right) $, so $f^{\ell
}\left( M^{l}(\alpha ,\beta )\right) $ will always have the $y$ coordinate
bigger than zero. This means that topologically the set $f^{\ell }\left(
M^{l}(\alpha ,\beta )\right) $ is to the left of $N\left( \bar{r},r\right) $
along the $x$ coordinate (see Figure \ref{fig:M-h-set}-right). Therefore the
first condition from (\ref{eq:covering-exit-cond-1}) in the definition of
the covering will be fulfilled.

Any point from $M(\alpha ,\beta )$ that enters $S_{I,L,R}^{u}$ does so
through $\left\{ x-y<\bar{r}\right\} $. Once a trajectory enters, it can not
exit through $\partial S_{I,L,R}^{u}\setminus \{q=R\}$, so
\begin{equation*}
f^{\ell }\left( M^{r}(\alpha ,\beta )\right) \cap N^{+}\left( \bar{r}%
,r\right) =\emptyset ,
\end{equation*}%
which means that condition (\ref{eq:covering-entry-cond}) definition of the
covering is fulfilled.
\end{proof}

Observe that%
\begin{equation*}
M(\bar{r},r)=\mathcal{S}\left( N^{T}(\bar{r},r)\right) .
\end{equation*}

\begin{lemma}
\label{lem:to-M-alpha-beta}Let $\bar{r},r>0$ be fixed. For every $\beta>0$
there exist $m>0$ and $\alpha>0$ (the $m$ and $\alpha$ depend on the choice
of $\bar{r},r,\beta$) such that
\begin{equation}
M(\bar{r},r)\overset{f^{m}}{\implies}M\left( \alpha,\beta\right)
\label{eq:final-covering-step}
\end{equation}
\end{lemma}

\begin{proof}
From the symmetry of the system, any point in $M^{-}(\bar{r},r)=M^l(\bar{r}%
,r)\cup M^r(\bar{r},r)$ will flow out of $\mathcal{S}\left(
S_{I,L,R}^{u}\right) $. From the definition of $f$, if $\mathbf{x}\in
M^{-}\left( \bar{r},r\right) $ then $f^{i}\left( \mathbf{x}\right) $ will
not return to $\mathcal{S}\left( S_{I,L,R}^{u}\right) $. Moreover, the
images of the components of $M^{-}(\bar{r},r)$ will remain on different
sides of $\mathcal{S}\left( S_{I,L,R}^{u}\cup S_{I,R_{1},R_{2},\varepsilon
}\right) .$ This implies (\ref{eq:covering-exit-cond-1}) for any choice of $%
\alpha ,\beta $.

As long as the flow which starts from $\mathbf{x}\in M(\bar{r},r)$ remains
in $\mathcal{S}\left( S_{I,L,R}^{u}\right) $, it approaches $\Lambda _{I}$
(This follows from the outflowing property in $S_{I,L,R}^{u}$ and the
symmetry of the system). Therefore, by compactness of $M(\bar{r},r)$, one
can find a sufficiently large $m$ so that if $\mathbf{x}\in M(\bar{r},r)$
and $f^{m}\left( \mathbf{x}\right) \in \mathcal{S}\left(
S_{I,L,R}^{u}\right) $ then $f^{m}\left( \mathbf{x}\right) \in \left\{
x+y<\beta \right\} .$ Let us fix such $m$. We can now choose $\alpha >0$
sufficiently small so that if $\mathbf{x}\in M(\bar{r},r)$ and $f^{m}\left(
\mathbf{x}\right) \in \mathcal{S}\left( S_{I,L,R}^{u}\right) $ then $%
f^{m}\left( \mathbf{x}\right) \in \left\{ x+y>\alpha \right\} $. We have
thus established (\ref{eq:final-covering-step}).
\end{proof}

We are now ready to prove Theorem \ref{th:main}.

\begin{proof}[Proof of Theorem \protect\ref{th:main}]
Let us start by choosing sufficiently small $0<\bar{r}<r$ so that we have%
\begin{equation}
N\left( \bar{r},r\right) \overset{f^{k}}{\implies }N\left( r\right) \overset{%
f_{1}}{\implies }N_{1}\overset{f_{0}}{\implies }N_{0}.
\label{eq:cover-seq-part-1}
\end{equation}%
This can be done by Lemmas \ref{lem:from-N1-to-N0}, \ref{lem:from-N2-to-N1}, %
\ref{lem:from-N3-to-N2}. We choose $r$ and $\bar{r}$ to be small enough, so that $\pi
_{q,p}N\left( \bar{r},r\right) \subset \lbrack -\sqrt{r^{\ast }},\sqrt{%
r^{\ast }}]^{2}$ where $r^{\ast }$ is the constant from Proposition \ref%
{th:orbits-enter-Scu}. We also choose $r$ to be small enough so that
\begin{equation*}
h_{P^{-}}:=W_{\Lambda _{I}}^{u}\cap \Sigma \cap N\left( \bar{r},r\right)
\end{equation*}%
is a horizontal disc in $N\left( \bar{r},r\right) $ and that
\begin{equation*}
v_{P^{+}}:=W_{\Lambda _{I}}^{s}\cap \Sigma \cap M\left( \bar{r},r\right)
\end{equation*}%
is a vertical disc in $M\left( \bar{r},r\right) $; see Figure \ref%
{fig:discs_H_P}. The $\bar{r}$ and $r$ will remain fixed throughout
the proof.

 We also define a horizontal disc $h_{H^{-}}$ in $N\left(
\bar{r},r\right) $ to be the lower boundary of $N\left( \bar{r},r\right) $
and define the vertical disc $v_{H^{+}}$ in $M\left( \bar{r},r\right) $ to
be the left boundary of $M\left( \bar{r},r\right) $; see Figure~\ref%
{fig:discs_H_P}.

All orbits which pass through $|h_{P^{-}}|$ converge backwards in time to $\Lambda_I$, hence they belong to $P^{-}$. Similarly, orbits which
pass through $|v_{P^{+}}|$ belong to $P^{+}.$ By Item 2 of Proposition \ref%
{th:orbits-enter-Scu}, orbits which pass through $|v_{H^{+}}|$ belong to
$H^{+}.$ By the symmetry of the system, orbits that pass through $|h_{H^{-}}|$
belong to $H^{-}.$ We will use the discs $v_{H^{+}},$ $h_{H^{-}},$ $v_{P^{+}}
$ and $h_{P^{-}}$ to obtain orbits from $H^{-}$ or $P^-$ to $H^+$ or $P^{+}$.
\begin{figure}[tbp]
\begin{center}
\includegraphics[width=2in]{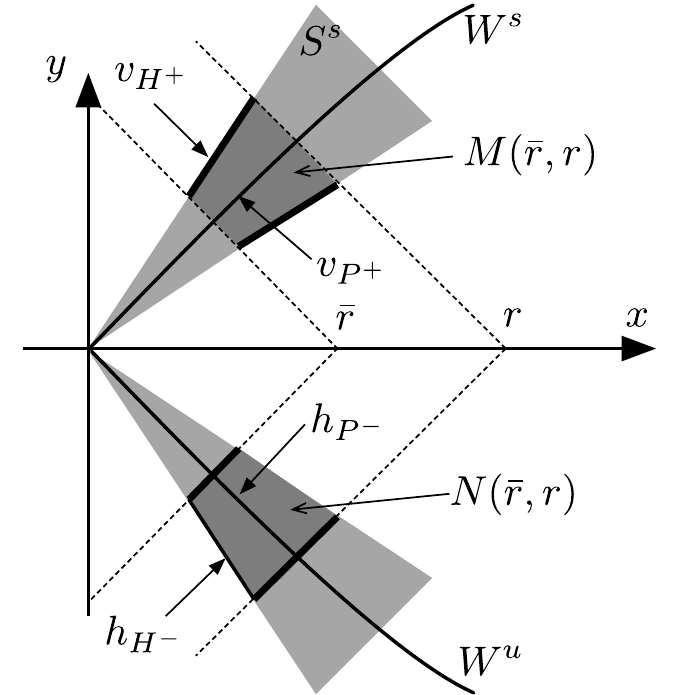}
\end{center}
\caption{The horizontal discs $h_{P^{-}},$ $h_{H^{-}}$ and the vertical
discs $v_{P^{+}}$ and $v_{H^{+}}.$}
\label{fig:discs_H_P}
\end{figure}

The idea for the proof of $OS^{\pm }$ is to consider a sequence $\left\{
\beta _{i}\right\} $ with $\beta _{i}>0$, and to construct a sequence of
coverings, which will successively link $N_{0}$ with $M\left( \alpha
_{i},\beta _{i}\right) $ (each time returning back to $N_{0}$), for suitably chosen $%
\alpha _{i}$. Below we will show how such link can be constructed for a fixed $%
\beta $. Later on in the proof this construction will be repeated for each $%
\beta _{i}$ from $\left\{ \beta _{i}\right\} $.

By Corollaries \ref{cor:from-N0-to-SN1}, \ref{cor:from-SN1-to-SN2}, \ref%
{cor:from-SN2-to-SN3}, there exists
\begin{equation}
N_{0}\overset{\tilde{f}_{0}}{\impliedby }\mathcal{S}\left( N_{1}^{T}\right)
\overset{\tilde{f}_{1}}{\impliedby }\mathcal{S}\left( N^{T}\left( r\right)
\right) \overset{\tilde{f}^{k}}{\impliedby }\mathcal{S}\left( N^{T}(\bar{r}%
,r)\right) =M(\bar{r},r).  \label{eq:cover-seq-part-2}
\end{equation}%
For a given small $\beta $, by Lemma \ref{lem:to-M-alpha-beta} there exists
an $\alpha =\alpha \left( \beta \right) >0$ and $m=m\left( \beta \right) \in
\mathbb{N}$ such that
\begin{equation}
M(\bar{r},r)\overset{f^{m}}{\implies }M\left( \alpha ,\beta \right) .
\label{eq:cover-seq-part-3}
\end{equation}%
By Lemma \ref{lem:from-M-to-Nr}, there exists a $\ell =\ell \left( \alpha
,\beta \right) $ such that%
\begin{equation}
M\left( \alpha ,\beta \right) \overset{f^{\ell }}{\implies }N\left( \bar{r}%
,r\right) .  \label{eq:cover-seq-part-4}
\end{equation}%
For the above choice of $\alpha ,m$ and $\ell $, combining (\ref%
{eq:cover-seq-part-1}--\ref{eq:cover-seq-part-4}) we obtain
\begin{align}
& N_{0}\overset{\tilde{f}_{0}}{\impliedby }\mathcal{S}\left(
N_{1}^{T}\right) \overset{\tilde{f}_{1}}{\impliedby }\mathcal{S}\left(
N^{T}\left( r\right) \right) \overset{\tilde{f}^{k}}{\impliedby }M(\bar{r},r)%
\overset{f^{m}}{\implies }  \label{eq:coverings-oscillatory-link-1} \\
& \qquad \qquad \qquad \overset{f^{m}}{\implies }M\left( \alpha ,\beta
\right) \overset{f^{\ell }}{\implies }N\left( \bar{r},r\right) \overset{f^{k}%
}{\implies }N\left( r\right) \overset{f_{1}}{\implies }N_{1}\overset{f_{0}}{%
\implies }N_{0}.  \notag
\end{align}%
To simplify the notation, let us denote such sequence of coverings by
\begin{equation}
N_{0}\underset{\beta }{\Longrightarrow }N_{0}.
\label{eq:coverings-oscillatory-link-2}
\end{equation}%
In the notation, by writing the $\beta $ we emphasize that (\ref%
{eq:coverings-oscillatory-link-1}) includes the set $M\left( \alpha ,\beta
\right) $. This will play an important role in our arguments. Whenever we
write (\ref{eq:coverings-oscillatory-link-2}), we will understand this as
choosing the $\beta $ first, and then the $\alpha ,m$ and $\ell $ in (\ref%
{eq:coverings-oscillatory-link-1}) are chosen so that the sequence of
coverings is ensured. The important issue is that we can construct such
sequence for an arbitrarily small $\beta $.

Let us now introduce the following sequences of coverings:

\begin{description}
\item[($H^{-}$)] The sequence that will lead to hyperbolic
motions in backward time, which we shall denote as $\left( H^{-}\right) $, is%
\begin{equation}
N\left( \bar{r},r\right) \overset{f^{k}}{\implies }N\left( r\right) \overset{%
f_{1}}{\implies }N_{1}\overset{f_{0}}{\implies }N_{0}.
\label{eq:Hminus-sequence}
\end{equation}%
\item[($P^{-}$)] We will also use the sequence \eqref%
{eq:Hminus-sequence} for the proof of parabolic motions in
backward time, which we denote by $\left(P^{-}\right)$.

\item[($H^{+}$)] The sequence that will lead to hyperbolic
motions in forward time, which we shall denote as $\left( H^{+}\right) $,%
\begin{equation}
N_{0}\overset{\tilde{f}_{0}}{\impliedby }\mathcal{S}\left( N_{1}^{T}\right)
\overset{\tilde{f}_{1}}{\impliedby }\mathcal{S}\left( N^{T}\left( r\right)
\right) \overset{\tilde{f}^{k}}{\impliedby }M(\bar{r},r).
\label{eq:Hplus-sequence}
\end{equation}%
\item[($P^{+}$)] We will also use the sequence \eqref{eq:Hplus-sequence} for the proof of parabolic motions in forward
time, denoted by $\left( P^{+}\right)$.

\item[($B^{-}$)] The sequence that will lead to bounded motions in backward
time%
\begin{equation}
\cdots \underset{r}{\implies }N_{0}\underset{r}{\implies }N_{0}\underset{r}{%
\implies }N_{0}.  \label{eq:Bminus-sequence}
\end{equation}%
(Above coverings are expressed in our simplified notation (\ref%
{eq:coverings-oscillatory-link-2}) with $\beta =r$.)

\item[($B^{+}$)] The sequence that will lead to bounded motions in forward
time:%
\begin{equation}
N_{0}\underset{r}{\implies }N_{0}\underset{r}{\implies }N_{0}\underset{r}{%
\implies }\cdots .  \label{eq:Bplus-sequence}
\end{equation}

\item[($OS^{-}$)] For a fixed sequence $\ldots ,\beta _{-3},\beta
_{-2},\beta _{-1}\in \mathbb{R}$ we consider a sequence of coverings%
\begin{equation*}
\cdots \underset{\beta _{-3}}{\implies }N_{0}\underset{\beta _{-2}}{\implies
}N_{0}\underset{\beta _{-1}}{\implies }N_{0}.
\end{equation*}%
(These coverings are expressed in our simplified notation (\ref%
{eq:coverings-oscillatory-link-2}).)

\item[($OS^{+}$)] For a fixed sequence $\beta _{0},\beta _{1},\beta
_{2},\ldots \in \mathbb{R}$ we consider a sequence of coverings%
\begin{equation*}
N_{0}\underset{\beta _{0}}{\implies }N_{0}\underset{\beta _{1}}{\implies }%
N_{0}\underset{\beta _{2}}{\implies }\cdots
\end{equation*}
\end{description}

To obtain orbits from $X^{-}\cap Y^{+}$ for $X,Y\in \left\{ H,P,B,OS\right\}
$ we combine sequences $\left( X^{-}\right) $ with $\left( X^{+}\right) $.

For example, to obtain an orbit from $H^{-}\cap P^{+}$ we glue $\left(
H^{-}\right) $ with $\left( P^{+}\right) $ which gives%
\begin{equation*}
N\left( \bar{r},r\right) \overset{f^{k}}{\implies }N\left( r\right) \overset{%
f_{1}}{\implies }N_{1}\overset{f_{0}}{\implies }%
\NoStart%
\overset{\tilde{f}_{0}}{\impliedby }\mathcal{S}\left( N_{1}^{T}\right)
\overset{\tilde{f}_{1}}{\impliedby }\mathcal{S}\left( N^{T}\left( r\right)
\right) \overset{\tilde{f}^{k}}{\impliedby }M(\bar{r},r).
\end{equation*}%
The circle indicates where the sequences are glued. Now, by Theorem \ref%
{th:covering} we obtain an orbit starting from $|h_{H^{-}}|$ which goes to $|v_{P^{+}}|$, which proves that $H^{-}\cap P^{+}\neq \emptyset $.

As another example we will show how to obtain an orbit from $B^{-}\cap OS^{+}
$. By gluing $\left( B^{-}\right) $ with $\left( OS^{+}\right) $, we obtain%
\begin{equation}
\cdots \underset{r}{\implies }N_{0}\underset{r}{\implies }N_{0}\underset{r}{%
\implies }%
\NoStart%
\underset{\beta _{0}}{\implies }N_{0}\underset{\beta _{1}}{\implies }N_{0}%
\underset{\beta _{2}}{\implies }\cdots .  \label{eq:example-2-covering}
\end{equation}%
We can take the sequence $\{\beta _{i}\}$ used in $\left( OS^{+}\right) $ to
converge to zero. By Item 3 of Corollary~\ref{cor:covering}, from (\ref%
{eq:example-2-covering}) we obtain an orbit which passes through this
sequence. We clearly have bounded motions in backward time. The trajectory
going forwards in time will be making alternating visits between $N_{0}$ and
$M\left( \alpha _{i},\beta _{i}\right) $ for $i\in \mathbb{N}$. The smaller
the $\beta _{i}$, the closer are the sets $M\left( \alpha _{i},\beta
_{i}\right) $  to $\Lambda _{I}$; which means that the further they are
from the origin for the original system associated to \eqref{eq:Jacobi-integral}. This means that the
orbit belongs to $OS^{+}$.

All other types of motions follow from mirror arguments: gluing of sequences
$(X^{-})$ and $\left( Y^{+}\right) $, for $X,Y\in \left\{ H,P,B,OS\right\} $
and using Theorem \ref{th:covering} or Corollary \ref{cor:covering}.

The bound for $r_0$ follows from computing \[\{r=2/x^2\,|\, x\in \pi_x N_0\}=[0.4999086, 0.5001915].\]

By Theorem \ref{th:covering} we obtain a periodic orbit passing through a sequence $N_0 \underset{r}{\implies } N_0$. Such orbit passes through the set $N(\bar r,r)$. The smaller the $r$ the further (in the original coordinates of the system) is such set from the origin. Since we can take $r$ as small as we wish, we can obtain periodic orbits which reach as far from the origin as we wish. The orbits also passes through $N_0$ so they pass $r_0$ close to the origin.  

The total computation time for the computer assisted part of the proof (i.e. for the validation of the sector, the extended sector and for the validation of $N_{1}\overset{f_{0}}{\Longrightarrow}N_{0}$) was $155$ seconds, running on a single thread on a standard laptop.

This concludes our proof.
\end{proof}

\section{Appendix}

Here we write the properties of the local Brouwer degree \cite{MR0493564}. The
\emph{local Brouwer degree} of a continuous map $f:\mathbb{R}^{n}%
\rightarrow\mathbb{R}^{n}$ at some point $c\in\mathbb{R}^{n}$, $n>0,$ in a set
$U\subset\mathbb{R}^{n}$ is a certain number. Suppose that
\begin{equation}
\mbox{the set}\quad f^{-1}(c)\cap U\quad\mbox{is compact.} \label{eq:zcomp}%
\end{equation}
Then \emph{the local Brouwer degree of $f$ at $c$ in the set $U$} is well
defined. We denote it by $\deg(f,U,c)$.

If $\overline{U}\subset\mbox{dom}(f)$ and $\overline{U}$ is compact, then
(\ref{eq:zcomp}) follows from the condition
\begin{equation}
c\notin f(\partial U). \label{eq:znoboundary}%
\end{equation}

Let us summarize the properties of the local Brouwer degree.

\noindent\textbf{Degree is an integer.}
\begin{equation}
\deg(f,U,c)\in\mathbb{Z}. \label{eq:degint}%
\end{equation}

\noindent\textbf{Solution property.}
\begin{equation}
\mbox{If}\quad\deg(f,U,c)\neq0,\quad\mbox{then there exists }x\in
U\mbox{ with }f(x)=c.
\end{equation}

\noindent\textbf{Homotopy property.} Let $H:[0,1]\times U\rightarrow
\mathbb{R}^{n}$ be continuous. Suppose that
\begin{equation}
\bigcup_{\lambda\in\lbrack0,1]}H_{\lambda}^{-1}(c)\cap U\quad\mbox{is
compact}. \label{eq:zhomcomp}%
\end{equation}
Then
\begin{equation}
\forall\lambda\in\lbrack0,1]\quad\deg(H_{\lambda},U,c)=\deg(H_{0},U,c).
\label{eq:deghom}%
\end{equation}
If $[0,1]\times\overline{U}\subset\mbox{dom}(H)$ and ${\overline{U}}$ is
compact, then (\ref{eq:zhomcomp}) follows from the following condition
\begin{equation}
c\notin H([0,1],\partial U). \label{eq:fhcb}%
\end{equation}

\noindent\textbf{Local degree for affine maps.} Suppose that $f(x)=A(x-x_{0}%
)+c$, where $A$ is a linear map and $x_{0}\in\mathbb{R}^{n}$. If the equation
$A(x)=0$ has no nontrivial solutions (i.e. if $Ax=0$, then $x=0$) and
$x_{0}\in U$, then
\begin{equation}
\deg(f,U,c)=\mbox{sgn}(\det A). \label{eq:indA}%
\end{equation}

\bibliographystyle{plain}
\bibliography{bib}
\end{document}